\documentclass[a4paper,11pt]{article}

\usepackage[utf8]{inputenc}
\usepackage[T1]{fontenc}
\usepackage[british]{babel}
\usepackage{fullpage}
\usepackage{amsmath,amssymb,amsthm,amsfonts,mathtools}
\usepackage{hyperref}
\usepackage{graphicx}
\usepackage{tikz}
\usepackage[titletoc,title]{appendix}
\usepackage[expansion=false]{microtype}

\hypersetup{colorlinks=true,linkcolor=blue!55!black,citecolor=green!45!black,urlcolor=blue!55!black}

\newcommand{\mc}{\mathcal}
\newcommand{\mb}{\mathbb}
\newcommand{\R}{\mb R}
\newcommand{\N}{\mb N}
\newcommand{\Z}{\mb Z}
\newcommand{\T}{\mb T}
\newcommand{\M}{\mc M}
\newcommand{\W}{\mc W}
\newcommand{\Leb}{\mathrm{Leb}}
\newcommand{\Id}{\mathrm{Id}}
\newcommand{\esssup}{\operatorname*{ess\,sup}}
\newcommand{\dd}{\mathrm{d}}
\renewcommand{\phi}{\varphi}
\newcommand{\eps}{\varepsilon}
\newcommand{\Up}{\Upsilon}
\newcommand{\barW}{\overline{W}}
\newcommand{\STO}{\mc F}
\newcommand{\LL}{\mathbb{L}}

\DeclareMathOperator{\Lip}{Lip}
\DeclareMathOperator{\BV}{BV}

\newtheorem{theorem}{Theorem}[section]
\newtheorem{proposition}[theorem]{Proposition}
\newtheorem{lemma}[theorem]{Lemma}
\newtheorem{corollary}[theorem]{Corollary}
\theoremstyle{definition}
\newtheorem{definition}[theorem]{Definition}
\newtheorem{example}[theorem]{Example}
\newtheorem{assumption}[theorem]{Standing Assumption}
\theoremstyle{remark}
\newtheorem{remark}[theorem]{Remark}

\newtheorem*{thmC}{Theorem C}
\newtheorem*{thmE}{Theorem E}

\begin{document}

\title{\bf Linear Response for Self-Consistent Transfer Operators on Graphons}

\author{%
Herbert M.\,C. Maquera\thanks{Independent researcher.
\texttt{herberth.maq@gmail.com}}
}

\date{\today}

\maketitle

\begin{abstract}
We develop the response theory of the equilibrium state of a large heterogeneous
network of weakly coupled uniformly expanding circle maps. The thermodynamic limit of
such a network is governed by a Self-Consistent Transfer Operator (STO) acting on
measures fibred over a graphon $W$, as introduced in \cite{MPT2025}. We prove that in
the weak coupling regime the unique attracting fixed point $\phi^\star(\alpha)$ of the
STO depends on the coupling strength in a Lipschitz, and then differentiable, way in a
scale of fibrewise Sobolev spaces, and we derive a linear response formula.

Our main result is that at zero coupling the response factorises: its derivative along
the fibre $z$ equals $\barW(z)\,\Up$, where $\barW$ is the degree function of the
graphon and $\Up=-(\Id-P)^{-1}[P(H_\rho\rho)]'$ is a universal profile depending only
on the local dynamics and on the coupling function. To first order the heterogeneity of
the statistical response across the network is therefore dictated by the degree profile
alone, and the finer topology is invisible; the susceptibility of a node is
proportional to its degree, so hubs respond most. The topology enters at higher order
through a tree expansion: the perturbative coefficients of the equilibrium, defined
unconditionally by an explicit recursion, are combinations --- with universal
coefficient functions --- of the rooted homomorphism densities of trees with $m$ edges.
We prove that the first two coefficients are the genuine asymptotic ones, the second
producing the mean degree of the neighbours of $z$. We also obtain a response formula for perturbations of the graphon
itself, Green--Kubo expressions for observables, and the prediction that in a finite
network the linear susceptibility of a node is proportional to its normalised degree.
\end{abstract}

\tableofcontents

\section{Introduction}\label{sec:intro}

\subsection{Motivation}

A basic question in the ergodic theory of dynamical systems is how the statistical
properties of a system change when the system itself is perturbed. If the physical
invariant measure depends differentiably on the perturbation parameter, the system is
said to exhibit \emph{linear response}, and an explicit expression for the derivative
is called a \emph{linear response formula}. Rigorous linear response originates in
Ruelle's work on Axiom A systems \cite{Ruelle1997,Ruelle1998,Ruelle2009}, and is by now
available well beyond the uniformly hyperbolic setting: for intermittent maps
\cite{BaladiTodd2016,Korepanov2016,BahsounSaussol2016}, for piecewise expanding unimodal
maps \cite{BaladiSmania2008,BaladiSmania2010}, for partially hyperbolic systems
\cite{Dolgopyat2004,Zhang2018}, and for systems with noise
\cite{GalatoloGiulietti2019,GalatoloSedro2020}. Caution is required, since
counterexamples are also known
\cite{Baladi2007,BaladiBenedicksSchnellmann2015,BaladiKunaLucarini2017,BahsounGalatolo2024};
see \cite{Baladi2014} for a survey. Linear response is also the mathematical backbone of
the fluctuation--dissipation formalism in climate science \cite{Lucarini2014,Majda2012}
and of the programme of \emph{control of statistical properties} of complex systems
initiated by MacKay \cite{MacKay2018} and developed in
\cite{GalatoloPollicott2017,Kloeckner2018,FroylandGalatolo2025}.

For \emph{extended} systems made of many interacting units the relevant object is the
equilibrium state of the thermodynamic limit. For globally coupled maps this limit is
governed by a \emph{self-consistent transfer operator} (STO): a nonlinear operator
$\phi\mapsto L_\phi\phi$ in which the transfer operator $L_\phi$ is itself built from a
mean field computed from $\phi$. The rigorous study of such operators was initiated by
Keller \cite{Keller2000} and by S\'elley and B\'alint \cite{SelleyBalint2016}, and has
developed rapidly
\cite{Balint2018,SelleyTanzi2021,SelleyTanzi2022,Galatolo2022,BahsounLiveraniSelley2023,BahsounKorepanov2024,BahsounLiverani2025,Tanzi2023,CastorriniGalatoloTanzi2025};
outside the weak coupling regime these systems display loss of uniqueness,
synchronisation and bifurcations \cite{SelleyTanzi2022,BahsounLiverani2025}.
Linear response for STOs in the globally coupled case was proved by S\'elley and Tanzi
\cite{SelleyTanzi2021} by a cone and Hilbert-metric method, and by Galatolo
\cite{Galatolo2022} inside a general framework for nonlinear perturbations of linear
transfer operators. Wormell and Gottwald
\cite{WormellGottwald2018,WormellGottwald2019} gave numerical and partially rigorous
evidence that a high-dimensional compound system may display linear response even when
its units do not; for coupled map lattices with short-range interaction, smooth
dependence of the SRB measure was obtained by Jiang and de la Llave
\cite{JiangdelaLlave2000,JiangdelaLlave2006}. We refer to \cite{Tanzi2023review} for a
review.

All the STO results quoted above concern \emph{homogeneous} interactions: every unit
feels the same mean field. Real networks are not homogeneous, and it is precisely the
heterogeneity of the interaction structure that produces the most interesting
collective phenomena, such as the coexistence of synchronised hubs with erratic
low-degree nodes \cite{PereiraVanStrienTanzi2020,TanziPereiraVanStrien2019,Corder2023,Bian2025}. In joint work with
T.~Pereira and M.~Tanzi \cite{MPT2025} the author introduced a Self-Consistent Transfer
Operator for \emph{heterogeneously} coupled expanding maps, in which the interaction structure
survives the thermodynamic limit as a \emph{graphon} $W:[0,1]^2\to\R$, the dense graph
limit of the adjacency matrices. Being a dense limit, this formalism does not reach the
sparse regime, for which the appropriate limit objects are graphexes, $L^p$-graphons or
the \emph{graphops} of Kuehn \cite{Kuehn2020}; we do not pursue that direction here. The state of the system is then not a single
probability density on $\T$ but a measurable family $\{\phi_z\}_{z\in[0,1]}$ of
densities, one \emph{local mean field per fibre}. It was proved in \cite{MPT2025} that
for small coupling the associated STO has a unique attracting fixed point, and that the
regularity of $z\mapsto\phi^\star_z$ is inherited from that of $z\mapsto W(z,\cdot)$.

The present paper develops the response theory of this object. Two questions drive the
analysis.

\begin{itemize}
\item[(Q1)] \emph{Coupling response.} How does the equilibrium state
$\phi^\star(\alpha)$ change when the coupling strength $\alpha$ varies? Is the response
uniform across the network, or do different nodes respond differently?
\item[(Q2)] \emph{Structural response.} How does the equilibrium state change when the
\emph{network} is perturbed, $W\rightsquigarrow W+\eps V$? Which features of the
topology are visible in the response?
\end{itemize}

Question (Q2) has no counterpart in the homogeneous theory, and is the one closest to
applications: it asks how the statistics of a network dynamical system react to
rewiring.

\subsection{Informal statement of the results}

Throughout, $f\in\mc C^{r}(\T,\T)$ is a uniformly expanding circle map of degree
$d\ge2$ with invariant density $\rho$ and transfer operator $P$, and
$h\in\mc C^{r}(\T\times\T,\R)$ is the coupling function. For $\psi\in L^1(\T)$ set
\[
H_\psi(x):=\int_\T h(x,y)\,\psi(y)\,\dd y ,
\]
and let
\[
\barW(z):=\int_0^1 W(z,z')\,\dd z'
\]
be the \emph{degree function} of the graphon, the thermodynamic limit of $k_i/N$ along
the fibre. The fibre dynamics in the thermodynamic limit is
\[
F_{\alpha,\phi,z}(x)=f(x)+\alpha\,A_{\phi,z}(x),\qquad
A_{\phi,z}=\int_0^1 W(z,z')\,H_{\phi_{z'}}\,\dd z' ,
\]
and the STO is $(\STO_\alpha\phi)_z=L_{\alpha,\phi,z}\phi_z$, with $L_{\alpha,\phi,z}$
the transfer operator of $F_{\alpha,\phi,z}$.

Theorems~\ref{thm:stability} and \ref{thm:LR} below (Theorems A and B) assert that
$\alpha\mapsto\phi^\star(\alpha)$ is Lipschitz and differentiable in the fibrewise
Sobolev scale $\W^k=L^\infty([0,1],W^{k,1}(\T))$, and that its derivative solves the
linear equation
\[
\big(\Id-\LL_{\alpha,\phi^\star}-\alpha\,\mc K_{\alpha,\phi^\star}\big)\,
\partial_\alpha\phi^\star=\mc K_{\alpha,\phi^\star}\phi^\star ,
\]
where $\LL$ acts fibrewise and $\mc K$ is the self-consistency operator. The second
structural ingredient is that $\mc K$ is an integral operator in $z$ with kernel $W$
and is \emph{smoothing}, which makes the Neumann series for the resolvent converge in
the strong norm.

The word ``heterogeneous'' is used in two different senses in this literature. In the survey \cite[\S3.2.2]{Tanzi2023review}, \emph{linear
response for heterogeneous systems} refers to the work of Wormell and Gottwald
\cite{WormellGottwald2018,WormellGottwald2019} on systems of units that are
\emph{individually different} --- logistic maps with i.i.d.\ parameters $a_i$ --- but
still coupled all-to-all through a single mean field; the remarkable phenomena there
concern the emergence, or failure, of linear response for the global observable
depending on the smoothness of the distribution of the $a_i$. In the present paper the
units are identical and it is the \emph{interaction structure} that is heterogeneous.
The two notions are independent, and to our knowledge the response theory of the
second one has not been considered before: the survey \cite[\S3.2]{Tanzi2023review},
which reviews the state of the art on linear response for STOs, covers
\cite{SelleyTanzi2021,Galatolo2022} for all-to-all couplings and
\cite{WormellGottwald2018,WormellGottwald2019} for heterogeneous units, and contains
no result in which the response depends on the position of a unit in a network.

These two statements are the heterogeneous counterparts of known results in the
all-to-all case, and the correspondence is worth recording. Lipschitz dependence of
the fixed density on the coupling strength was proved for globally coupled expanding
circle maps with a \emph{discontinuous} interaction by B\'alint, Keller, S\'elley and
T\'oth \cite[Thm.~2]{Balint2018}, where differentiability is genuinely obstructed by
the discontinuity; for \emph{smooth} interactions, differentiability and a response
formula were obtained by S\'elley and Tanzi \cite{SelleyTanzi2021} and, in a general
functional-analytic framework, by Galatolo \cite{Galatolo2022}. Theorem A is the
graphon analogue of \cite[Thm.~2]{Balint2018} and Theorem B of
\cite[Thm.~1.2]{SelleyTanzi2021}; what has no counterpart in that literature is the
dependence of the response on the position $z$ in the network, which is the subject of
Theorem C below.

The main structural result is what this equation gives at zero coupling.

\begin{thmC}[Degree-proportional response; Theorem~\ref{thm:zerocoupling}]
Let $\Up:=-(\Id-P)^{-1}\big[P(H_\rho\,\rho)\big]'$. Then, in
$\W^{k}=L^\infty([0,1],W^{k,1}(\T))$,
\[
\partial_\alpha \phi^\star_z(\alpha)\Big|_{\alpha=0}=\barW(z)\,\Up
\qquad\text{for a.e. } z\in[0,1],
\]
equivalently $\phi^\star_z(\alpha)=\rho+\alpha\,\barW(z)\,\Up+o(\alpha)$ uniformly in
$z$.
\end{thmC}

The first-order response therefore \emph{factorises} into a structural factor
$\barW(z)$ and a dynamical profile $\Up$. Two consequences deserve emphasis.

\begin{itemize}
\item[(i)] \emph{Only degrees matter at first order.} Two graphons with the same degree
function produce indistinguishable first-order responses, however different their
topologies; in particular for a graphon with constant degree function the first-order
response is homogeneous across the network and coincides with the all-to-all one.
\item[(ii)] \emph{Hubs respond most.} Since
$\|\partial_\alpha\phi^\star_z|_{0}\|_{W^{k,1}}=|\barW(z)|\,\|\Up\|_{W^{k,1}}$, for a
\emph{non-negative} graphon --- the case of an actual network --- and $\Up\neq0$ the
magnitude of the first-order response is an increasing function of the degree. For a
rank-one (configuration-model type) graphon $W(z,z')=g(z)g(z')$ with $g\ge0$ one has
$\barW(z)=g(z)\,\|g\|_{L^1}$, so that the magnitude of the susceptibility of a node is
exactly proportional to its expected degree.
\end{itemize}

The topology does appear at higher order, and in a combinatorially clean way. Let
$\mc T_j$ be the set of isomorphism classes of finite rooted trees with $j$ edges and
let $t_z(T,W)$ be the rooted homomorphism density of $T$ with root pinned at $z$ (see
\eqref{eq:homdensity}).

\begin{thmE}[Tree expansion; Theorem~\ref{thm:tree}]
Let $\phi_{j}$, $1\le j\le r-2$, be the \emph{formal susceptibilities} of the system,
defined unconditionally by the recursion that any asymptotic expansion of
$\phi^\star(\alpha)$ would have to satisfy (Proposition~\ref{prop:formal}). Then
\[
\phi_{j,z}=\sum_{T\in\mc T_j} t_z(T,W)\;\Up_T ,
\]
where the profiles $\Up_T$ depend only on $f$ and $h$, not on $W$ and not on $z$.
\end{thmE}

The formal susceptibilities matter only insofar as they describe the equilibrium, and
$\phi_1,\phi_2$ do: by Proposition~\ref{prop:formalgenuine} they are the genuine first
and second order coefficients of $\phi^\star(\alpha)$, while whether the formal series
describes $\phi^\star$ beyond second order is open. At order two the two trees with two
edges --- the path $P_2$ and the cherry $\Lambda_2$ --- appear with explicitly computed
profiles (Theorem~\ref{thm:secondorder}). The path
contributes $t_z(P_2,W)=\int_0^1W(z,z')\barW(z')\dd z'$, the \emph{mean degree of the
neighbours} of $z$. Thus the network becomes visible at second order through the
degrees of the neighbours, at third order through degrees at distance two and through
triples of neighbours, and so on; and no formal susceptibility sees cycle densities.
This gives a dynamical meaning to the homomorphism densities of graph limit theory
\cite{LovaszSzegedy2006,LovaszBook,BorgsChayes2008}.

We complement these results with a linear response formula for perturbations of the
graphon itself (Theorem~\ref{thm:structural}), from which the same degree
proportionality emerges at leading order in $\alpha$; Green--Kubo formulae for
observables (Corollary~\ref{cor:GK}); a statement on the regularity of the response in
the graphon variable (Proposition~\ref{prop:zregularity}); and the finite-$N$
interpretation of Section~\ref{sec:finiteN}.

\subsection{Method, and comparison with the literature}

Two routes to linear response for STOs exist in the literature. S\'elley and Tanzi
\cite{SelleyTanzi2021} construct an invariant set of \emph{curves}
$\alpha\mapsto\gamma(\alpha)$ of densities with controlled $\mc C^k$ norms of
$\partial_\alpha\gamma$ and $\partial^2_\alpha\gamma$, and identify the limit of
$\STO^n\gamma$ with the curve of fixed points, obtaining $\mc C^1$ dependence
\emph{a priori}; the response formula then follows from a resolvent computation.
Galatolo \cite{Galatolo2022} works with an abstract triple
$B_{ss}\subset B_s\subset B_w$ and proves statistical stability first, then
differentiability, using only the resolvent of the unperturbed (linear) operator.

A third ingredient, which we use throughout, is the structure of the differential of
an STO. That it splits into a frozen transfer-operator part plus a self-consistency
part,
\begin{equation}\label{eq:knownsplit}
d\mc L_{\delta,h}=L_{\delta,h}+\partial L_{\delta,h},
\end{equation}
was isolated by Castorrini, Galatolo and Tanzi
\cite{CastorriniGalatoloTanzi2025}, who used it to study local convergence to
equilibrium of mean field systems, including outside the weak coupling regime (see
also \cite{CastorriniGalatoloTanzi2024} for a cone-contraction approach to the
stability of fixed points of nonlinear STOs); our
Theorem~\ref{thm:diff} is the graphon analogue of \eqref{eq:knownsplit}, with
$\LL_{\alpha,\phi}$ in the role of $L_{\delta,h}$ and $\alpha\mc K_{\alpha,\phi}$ in
the role of $\partial L_{\delta,h}$. We emphasise that neither
\cite{CastorriniGalatoloTanzi2025} nor the surveyed literature
\cite{Tanzi2023review} derives a response formula from this splitting: the differential
is used there for stability, not for response.

We follow the second route, adapted to the fibred setting, and what we exploit beyond
\eqref{eq:knownsplit} are three features special to graphon couplings.

\emph{(1) The mean field is linear and smoothing.} The map $v\mapsto A_{v}$ is linear
in $v$ and bounded from the weakest space $\W^0$ into
$L^\infty([0,1],\mc C^{r}(\T))$, because $H_\psi$ is as smooth as $h$ for any
$\psi\in L^1$. Consequently the self-consistency operator $\mc K$ is bounded from
$\W^0$ to $\W^k$ (Proposition~\ref{prop:Kbounded}), so that the ``small perturbation''
hypothesis holds in the strongest possible mixed norm and the Neumann series for the
resolvent converges with no loss of derivatives. This is the counterpart, in our
setting, of the boundedness of $\partial L_{\delta,h}:V_w\to B_{ss}$ assumed in
\cite[Prop.~5]{CastorriniGalatoloTanzi2025}; here it is not an assumption but a
consequence of the smoothness of $h$.

\emph{(2) Linearity identifies two a priori different operators.} Because
$v\mapsto A_v$ is linear, the $\alpha$-derivative of the STO at frozen state coincides
with $\mc K$ evaluated at the state itself (Lemma~\ref{lem:Kphi}). This is what makes
the right-hand side of the response equation equal to $\mc K\phi^\star$.

\emph{(3) Non-locality in $z$ with kernel $W$.} The operator $\mc K$ integrates against
$W(z,\cdot)$, so the response equation is a Fredholm-type integral equation over the
fibre space whose iteration generates the tree expansion of Theorem E. In the
homogeneous case $\partial L_{\delta,h}$ is a rank-one operator
\cite[\S3]{CastorriniGalatoloTanzi2025}; here it has rank as large as that of the
kernel $W$, and this is precisely what makes the response heterogeneous.

Two further technical choices differ from \cite{SelleyTanzi2021}. First, we work in
fibrewise Sobolev spaces $W^{k,1}(\T)$ rather than $\mc C^k(\T)$: the Lasota--Yorke
constants there do not carry a factor equal to the degree of $f$
(Lemma~\ref{lem:LY} and Remark~\ref{rem:nodegree}), so no extra technical condition
such as \cite[Eq.~(8)]{SelleyTanzi2021} is needed. Second, all our perturbative
identities are \emph{exact}, with integral remainders
(Lemma~\ref{lem:exact}), which keeps the bookkeeping of lost derivatives transparent.

We stress a genuine difference in the model. In \cite{SelleyTanzi2021,Galatolo2022} the
coupling is a diffeomorphism composed with the local map, $F=f\circ\Phi$ or
$F=\Phi\circ f$; in \cite{MPT2025}, and hence here, the coupling is additive,
$F=f+\alpha A$, matching the finite-size model
$x_i(t+1)=f(x_i(t))+\frac{\alpha}{N}\sum_jA_{ij}h(x_i(t),x_j(t))$. This changes the
response formula, and the position of the transfer operator inside the derivative in
\eqref{eq:Upsilon} reflects it. In Remark~\ref{rem:crosscheck} we verify that our
differentiation lemma reproduces both the formula of \cite[Prop.~1.4]{SelleyTanzi2021}
and that of \cite[Prop.~29]{Galatolo2022} when specialised to their respective models.

\subsection{Organisation}

Section~\ref{sec:setting} fixes the setting, the function spaces and the standing
assumptions, and recalls what we import from \cite{MPT2025}. Section~\ref{sec:results}
states all the results. Section~\ref{sec:spectral} contains the uniform spectral
estimates, Section~\ref{sec:differential} the differential of the STO and its
resolvent, and Section~\ref{sec:proofs} the proofs of Theorems A, B and C.
Section~\ref{sec:higher} develops the higher order theory and proves the tree
expansion, Section~\ref{sec:structural} treats perturbations of the graphon,
and Section~\ref{sec:finiteN} discusses finite networks.
Appendix~\ref{app:basics} collects background material,
Appendix~\ref{app:calculus} the exact perturbative identities, and
Appendix~\ref{app:longproofs} the long proofs.

\section{Setting, function spaces and standing assumptions}\label{sec:setting}

\subsection{The model and its thermodynamic limit}

Let $\T=\R/\Z$ with Lebesgue (Haar) measure $\Leb$, also written $\dd x$. We consider
the finite network dynamics
\begin{equation}\label{eq:finite}
x_i(t+1)=f(x_i(t))+\frac{\alpha}{N}\sum_{j=1}^N A_{ij}\,h\big(x_i(t),x_j(t)\big)
\quad \mathrm{mod}\ 1,\qquad i=1,\dots,N,
\end{equation}
with $A=(A_{ij})$ the adjacency matrix of a directed graph on $N$ nodes. Following
\cite{MPT2025} we adopt the convention that $A_{ij}=1$ means that node $i$
\emph{receives} a connection from node $j$, i.e.\ the directed edge $j\to i$ is
present and $j$ influences $i$; with this convention $k_i:=\sum_jA_{ij}$ is the
\emph{in-degree} of node $i$, which is the quantity appearing in
\eqref{eq:finite}. As recalled in
Appendix~\ref{app:graphons}, \eqref{eq:finite} can be rewritten using the step graphon
$W^{(N)}$ associated with $A$, and under the hypotheses (H1)--(H2) of \cite{MPT2025}
the $N\to\infty$ limit of the one-node marginals is governed by the following object.

\begin{definition}[Local mean field]\label{def:meanfield}
Let $W$ be a graphon. For $v\in L^1([0,1]\times\T)$ and $\Leb$-a.e.\ $z\in[0,1]$ set
\begin{equation}\label{eq:A}
A_{v,z}(x):=\int_0^1 W(z,z')\Big(\int_\T h(x,y)\,v_{z'}(y)\,\dd y\Big)\dd z'
=\int_0^1 W(z,z')\,H_{v_{z'}}(x)\,\dd z' ,
\end{equation}
where $H_\psi(x):=\int_\T h(x,y)\psi(y)\dd y$. We write $A_{v}$ for the family
$\{A_{v,z}\}_{z}$.
\end{definition}

The map $v\mapsto A_v$ is \emph{linear}. This elementary observation is used repeatedly
and is responsible for several simplifications; see Lemma~\ref{lem:Kphi}.

\begin{definition}[Fibre maps and the STO]\label{def:STO}
For $\alpha\in\R$, $\phi\in\M_{1,\Leb}$ and a.e.\ $z\in[0,1]$ define
\begin{equation}\label{eq:F}
F_{\alpha,\phi,z}:=f+\alpha A_{\phi,z}\ :\ \T\to\T ,
\end{equation}
denote by $L_{\alpha,\phi,z}$ its transfer operator with respect to $\Leb$, and set
\begin{equation}\label{eq:STO}
(\STO_\alpha\phi)_z:=L_{\alpha,\phi,z}\,\phi_z .
\end{equation}
We write $P:=L_f$ for the transfer operator of the uncoupled map and $\rho$ for its
unique invariant probability density.
\end{definition}

As in \cite{MPT2025}, $\M_{1,\Leb}=\M_{1,\Leb}([0,1]\times\T)$ denotes the probability
measures on $[0,1]\times\T$ whose marginal on the first coordinate is $\Leb$, and we
identify such a measure with the family $\{\phi_z\}_z$ of densities of its
disintegration whenever these exist. Definition~\ref{def:STO} agrees with
\cite[Def.~2.2]{MPT2025}, since $(\STO\nu)_z=(F_{\nu,z})_*\nu_z$ and the push-forward
acts on densities as the transfer operator.

\subsection{Function spaces}

For $k\in\N_0$ let $W^{k,1}(\T)$ be the Sobolev space of $\psi\in L^1(\T)$ whose
distributional derivatives up to order $k$ lie in $L^1(\T)$, with
$\|\psi\|_{W^{k,1}}=\sum_{i=0}^k\|\psi^{(i)}\|_{L^1}$, so that $W^{0,1}=L^1$. We use
the fibred spaces
\begin{equation}\label{eq:spaces}
\W^k:=L^\infty\big([0,1],W^{k,1}(\T)\big),\qquad
\|v\|_{k}:=\esssup_{z\in[0,1]}\|v_z\|_{W^{k,1}(\T)} ,
\end{equation}
the Bochner spaces of strongly measurable, essentially bounded maps $z\mapsto v_z$.
Each $\W^k$ is a Banach space and $\W^{k+1}\hookrightarrow\W^{k}$ contractively. We set
\begin{equation}
\W^k_0:=\Big\{v\in\W^k:\ \int_\T v_z=0 \text{ a.e. }z\Big\},\qquad
\Delta^k:=\Big\{\phi\in\W^k:\ \phi_z\ge0,\ \int_\T\phi_z=1 \text{ a.e. } z\Big\}.
\end{equation}
Thus $\Delta^k$ is a closed convex affine slice of $\W^k$ with tangent space $\W^k_0$,
which is where the response lives. For $M>0$ and $j\in\N$ we write
\begin{equation}\label{eq:ball}
\Delta^j_M:=\big\{\phi\in\Delta^j:\ \|\phi\|_{i}\le M_i \text{ for } 0\le i\le j\big\},
\end{equation}
where the constants $M_i$ are the ones produced by Proposition~\ref{prop:apriori}
below.

\begin{remark}[Two facts used constantly]\label{rem:BVW11}
(i) For $\psi\in W^{1,1}(\T)$ the $\BV$ seminorm of \cite{MPT2025}, namely
$|\psi|_{\BV^1}=\sup\{\int_\T g'\psi:\ g\in\mc C^1,\|g\|_\infty\le1\}$, equals
$\|\psi'\|_{L^1}$, so that $\|\psi\|_{\BV^1}=\|\psi\|_{W^{1,1}}$. This lets us import
the memory-loss estimate of \cite{MPT2025} verbatim as the base case of our spectral
estimates.
(ii) For every $\psi\in L^1(\T)$ and $0\le j\le r$ one has
$\|H_\psi\|_{\mc C^{j}}\le\|h\|_{\mc C^{j}}\|\psi\|_{L^1}$: the mean field is as smooth
as $h$ regardless of the regularity of $\psi$.
\end{remark}

\begin{remark}[Measurability]\label{rem:meas}
Strong measurability of $z\mapsto\phi^\star_z$ is not automatic from the disintegration
statement alone. It holds for all objects considered here: by
\cite[Thm.~3.2]{MPT2025} the fixed point lies in the closure of a set of densities with
finite $\mathtt{var}_{p,\BV^1}$, hence $z\mapsto\phi^\star_z$ is essentially separably
valued and weakly measurable, and Pettis' theorem applies. For $\mc Kv$ measurability
follows from measurability of $z\mapsto W(z,\cdot)\in L^1([0,1])$, which holds by
Fubini.
\end{remark}

\subsection{Standing assumptions}

\begin{assumption}\label{ass:main}
Throughout we fix an integer $r\ge 6$ and assume:
\begin{itemize}
\item[(A1)] $f\in\mc C^{r}(\T,\T)$ is uniformly expanding, $|f'(x)|\ge\sigma>1$ for all
$x\in\T$, of degree $d\ge2$;
\item[(A2)] $h\in\mc C^{r}(\T\times\T,\R)$;
\item[(A3)] $W\in L^\infty([0,1]^2)$ is a graphon with
\begin{equation}\label{eq:varW}
\mathtt{var}_{p,L^1}(W)<\infty \qquad\text{for some } p\in(0,1],
\end{equation}
the $p$-variation of the sections being the one of \cite{MPT2025}, recalled in
Appendix~\ref{app:pvar}. Here \emph{graphon} means a bounded measurable kernel
$W:[0,1]^2\to\R$, not assumed symmetric and not assumed $[0,1]$-valued, the
generalisation appropriate to directed and weighted networks used in \cite{MPT2025}
and recalled in Appendix~\ref{app:graphons}. We write
$\|W\|_\infty:=\|W\|_{L^\infty([0,1]^2)}$ and
$\|W\|_{\infty,1}:=\esssup_z\|W(z,\cdot)\|_{L^1}$;
\item[(A4)] $|\alpha|<\alpha_0$, where $\alpha_0>0$ is small enough that the
conclusions of \cite[Thm.~3.2]{MPT2025} and of \cite[Lem.~4.5]{MPT2025} hold and that
\begin{equation}\label{eq:expcond}
|\alpha|\,\|h\|_{\mc C^{1}}\,\|W\|_{\infty,1}<\sigma-1 .
\end{equation}
\end{itemize}
\end{assumption}

Two comments on (A3). The boundedness $W\in L^\infty([0,1]^2)$ is stronger than the
minimal hypothesis $W\in L^\infty([0,1],L^1([0,1]))$ of \cite[Def.~2.2]{MPT2025} under
which the STO is defined; it is used in Lemma~\ref{lem:Abound}(iii), and every
statement not relying on that item holds with $\|W\|_\infty$ replaced by
$\|W\|_{\infty,1}$. The $p$-variation condition \eqref{eq:varW}, by contrast, is
\emph{not} optional and is not a technical convenience of ours: it is a hypothesis of
the fixed-point theorem of \cite{MPT2025} that we import as Theorem~\ref{thm:MPT},
where it is what makes the admissible class of densities precompact, and hence it is
required for the very object $\phi^\star(\alpha)$ whose response we study to be known
to exist. We therefore carry it in all our statements. It holds for every example
considered in this paper: it is immediate for constant graphons; for block graphons
the oscillation is supported within distance $r$ of the finitely many block
boundaries, giving $\mathtt{var}_{p,L^1}(W)<\infty$ for every $p\le1$; and if
$z\mapsto W(z,\cdot)$ is Lipschitz or H\"older of exponent $\ge p$ into $L^1([0,1])$
then \eqref{eq:varW} holds, which covers the spatially decaying and rank-one kernels
of Section~\ref{sec:examples}. See \cite[\S2.3]{MPT2025} for a longer list.

\begin{theorem}[{existence, uniqueness and local stability of the fixed point,
\cite[Thm.~3.2, Cor.~4.1]{MPT2025}}]\label{thm:MPT}
Under Assumption~\ref{ass:main} there is $\alpha_0>0$ such that for every
$|\alpha|<\alpha_0$ the operator $\STO_\alpha$ has a fixed point
$\phi^\star(\alpha)\in\M_{1,\Leb}$ with absolutely continuous disintegration, unique
within the cone class $\widetilde{\mc V}_a$ and locally exponentially attracting, and
$\esssup_z\|\phi^\star_z(\alpha)\|_{\mc C^2}<\infty$.
\end{theorem}

\begin{lemma}[{\cite[Lem.~4.3]{MPT2025}: uniform expansion and distortion}]
\label{lem:uniformMPT}
Under \eqref{eq:expcond} there are $\xi>1$ and $K,K'>0$ such that for all
$|\alpha|<\alpha_0$, all $\phi\in\M_{1,\Leb}$ and a.e.\ $z$,
$\inf_\T|F'_{\alpha,\phi,z}|\ge\xi$, $\|F_{\alpha,\phi,z}\|_{\mc C^3}\le K$ and
$\sup_\T|F''_{\alpha,\phi,z}/(F'_{\alpha,\phi,z})^2|\le K'$.
\end{lemma}

\begin{lemma}[{\cite[Lem.~4.5]{MPT2025}: exponential loss of memory}]\label{lem:MLMPT}
There are $\alpha_0>0$, $C_1\ge1$ and $\lambda>0$ such that for $|\alpha|<\alpha_0$,
for every sequence $\{\nu_n\}\subset\M_{1,\Leb}$, a.e.\ $z\in[0,1]$, every $n\ge0$ and
every $\psi\in\mc B_{\BV^1}$ with $\int_\T\psi=0$,
\[
\big\|\big(F_{\nu_n,z}\circ\cdots\circ F_{\nu_1,z}\big)_*\psi\big\|_{\BV^1}
\le C_1 e^{-\lambda n}\|\psi\|_{\BV^1}.
\]
\end{lemma}

\section{Statement of the results}\label{sec:results}

Recall $H_\psi(x)=\int_\T h(x,y)\psi(y)\dd y$, $\barW(z)=\int_0^1 W(z,z')\dd z'$,
$P=L_f$ and $P\rho=\rho$.

\subsection{The differential of the STO}

\begin{definition}[Fibrewise and self-consistency operators]\label{def:LK}
For $|\alpha|<\alpha_0$ and $\phi\in\Delta^{k+1}$ define, for $v\in\W^0$,
\begin{align}
(\LL_{\alpha,\phi}\,v)_z&:=L_{\alpha,\phi,z}\,v_z , \label{eq:LL}\\[2pt]
(\mc K_{\alpha,\phi}\,v)_z&:=-\Big[L_{\alpha,\phi,z}\big(A_{v,z}\,\phi_z\big)\Big]' .
\label{eq:KK}
\end{align}
\end{definition}

The operator $\LL_{\alpha,\phi}$ is \emph{local} in $z$: it is the direct integral of
the transfer operators of the frozen fibre maps. The operator $\mc K_{\alpha,\phi}$ is
\emph{non-local} in $z$: through $A_{v,z}$ it integrates $v$ against the kernel
$W(z,\cdot)$. This dichotomy is the analytic signature of the network structure.

\begin{theorem}[Differential of the STO]\label{thm:diff}
Let $1\le k\le r-3$ and $M>0$. There is $C_M>0$ such that for every $|\alpha|<\alpha_0$
and all $\phi,\phi+v\in\Delta^{k+2}$ with $\|\phi\|_{k+2},\|\phi+v\|_{k+2}\le M$,
\begin{equation}\label{eq:diffSTO}
\Big\|\STO_\alpha(\phi+v)-\STO_\alpha(\phi)
-\big(\LL_{\alpha,\phi}+\alpha\,\mc K_{\alpha,\phi}\big)v\Big\|_{k}
\le C_M\,|\alpha|\,\|v\|_{k+1}^{2} .
\end{equation}
In particular $\STO_\alpha$ is differentiable at $\phi$ with differential
$\LL_{\alpha,\phi}+\alpha\mc K_{\alpha,\phi}$, in the following precise sense, which is
the only one we use and the only one available: the admissible set $\Delta^{k+2}$ is a
closed affine slice of $\W^{k+2}$ and not an open subset of $\W^k$, so
\eqref{eq:diffSTO} is a statement about increments $v$ that are \emph{tangential} to
that slice, namely $v\in\W^{k+1}_0$ with $\phi+v\in\Delta^{k+2}$ and
$\|\phi+v\|_{k+2}\le M$, and the remainder is controlled by the stronger norm
$\|v\|_{k+1}$ rather than by $\|v\|_{k}$. No Fr\'echet differentiability on a
$\W^k$-neighbourhood is claimed, and none is needed: every increment occurring below is
a difference of two admissible densities. Moreover, at frozen $\phi$,
$\partial_\alpha(\STO_\alpha\phi)=\mc K_{\alpha,\phi}\phi$.
\end{theorem}

\subsection{Stability and linear response}

\begin{theorem}[Theorem A: Lipschitz statistical stability]\label{thm:stability}
Let $1\le k\le r-2$. There are $\alpha_1\in(0,\alpha_0]$ and $C>0$ such that for all
$\alpha,\beta\in(-\alpha_1,\alpha_1)$,
\[
\big\|\phi^\star(\alpha)-\phi^\star(\beta)\big\|_{k}\le C\,|\alpha-\beta| .
\]
\end{theorem}

\begin{theorem}[Theorem B: linear response]\label{thm:LR}
Let $1\le k\le r-3$. There is $\alpha_2\in(0,\alpha_1]$ such that for every
$\hat\alpha\in(-\alpha_2,\alpha_2)$, writing $\hat\phi=\phi^\star(\hat\alpha)$:
\begin{itemize}
\item[(i)] the operator
$\mc R_{\hat\alpha,\hat\phi}:=\Id-\LL_{\hat\alpha,\hat\phi}
-\hat\alpha\,\mc K_{\hat\alpha,\hat\phi}$ is invertible on $\W^{k}_0$;
\item[(ii)] the map $\alpha\mapsto\phi^\star(\alpha)\in\W^k$ is differentiable at
$\hat\alpha$, with
\begin{equation}\label{eq:LRgeneral}
\partial_\alpha\phi^\star(\alpha)\Big|_{\alpha=\hat\alpha}
=\mc R_{\hat\alpha,\hat\phi}^{-1}\ \mc K_{\hat\alpha,\hat\phi}\,\hat\phi ;
\end{equation}
\item[(iii)] $\hat\alpha\mapsto\partial_\alpha\phi^\star(\hat\alpha)$ is Lipschitz
from $(-\alpha_2,\alpha_2)$ into $\W^{k}$; in particular
$\alpha\mapsto\phi^\star(\alpha)$ is of class $\mc C^{1,1}$, hence $\mc C^1$, from
$(-\alpha_2,\alpha_2)$ to $\W^{k}$.
\end{itemize}
\end{theorem}

\begin{theorem}[Theorem C: the response at zero coupling factorises]
\label{thm:zerocoupling}
Let $1\le k\le r-3$ and set
\begin{equation}\label{eq:Upsilon}
\Up:=-(\Id-P)^{-1}\Big[P\big(H_\rho\,\rho\big)\Big]'\ \in\ W^{r-2,1}(\T),
\qquad \int_\T\Up=0 .
\end{equation}
Then
\begin{equation}\label{eq:LRzero}
\lim_{\alpha\to0}\ \esssup_{z\in[0,1]}
\left\|\frac{\phi^\star_z(\alpha)-\rho}{\alpha}-\barW(z)\,\Up\right\|_{W^{k,1}(\T)}=0 .
\end{equation}
In particular $\phi^\star_z(\alpha)=\rho+\alpha\,\barW(z)\,\Up+o(\alpha)$ uniformly in
$z$: the response is a rank-one object, a structural factor $\barW$ times a universal
dynamical profile $\Up$.
\end{theorem}

\begin{corollary}[Observables and Green--Kubo]\label{cor:GK}
Let $c\in L^\infty(\T)$ and let $\nu_{\alpha}$ denote the equilibrium measure with
disintegration $\{\phi^\star_z(\alpha)\}$. Then for a.e.\ $z$
\begin{equation}\label{eq:GK}
\frac{\dd}{\dd\alpha}\bigg|_{\alpha=0}\int_\T c\,\dd\nu_{\alpha,z}
=\barW(z)\int_\T c\,\Up
=-\,\barW(z)\sum_{n\ge0}\int_\T \big(c\circ f^{\,n}\big)\,
\Big[P\big(H_\rho\rho\big)\Big]'\,\dd x ,
\end{equation}
the series converging exponentially fast. Integrating over the network,
\[
\frac{\dd}{\dd\alpha}\bigg|_{\alpha=0}\int_{[0,1]\times\T}c\,\dd\nu_\alpha
=\Big(\int_0^1\!\!\int_0^1 W\Big)\int_\T c\,\Up ,
\]
so that the globally averaged susceptibility is proportional to the \emph{edge density}
of the graphon.
\end{corollary}

\begin{remark}[Numerical verification of Theorem~\ref{thm:zerocoupling}]
\label{rem:numerics}
We have tested \eqref{eq:LRzero} numerically for
$f(x)=2x+\tfrac1{10}\sin(2\pi x)$, which is uniformly expanding with $\sigma>1.37$,
the coupling function
$h(x,y)=\sin\big(2\pi(y-x)\big)+\tfrac12\cos\big(2\pi(x+2y)\big)$, and the two-block
graphon $W=\big(\begin{smallmatrix}1&0.3\\0.8&0.2\end{smallmatrix}\big)$ on
$[0,\tfrac12]\sqcup(\tfrac12,1]$, for which $\barW=(0.65,\,0.5)$. Densities were
discretised on $2^{10}$ equispaced points, transfer operators were evaluated through
the inverse branches of the lift with periodic cubic interpolation, $\rho$ and
$\phi^\star(\alpha)$ were obtained by iteration, and $\Up$ by summing the Neumann
series $-\sum_{n\ge0}P^n\big[P(H_\rho\rho)\big]'$. On the first block the quantity
$\big\|(\phi^\star_z(\alpha)-\rho)/\alpha-\barW(z)\Up\big\|_{\infty}$ takes the values
$1.33\cdot10^{-3}$, $6.75\cdot10^{-4}$, $3.40\cdot10^{-4}$, $1.71\cdot10^{-4}$ for
$\alpha=0.04,\,0.02,\,0.01,\,0.005$, i.e.\ it halves when $\alpha$ halves, so that the
remainder in \eqref{eq:LRzero} is in fact $O(\alpha)$ in this example. The rank-one
structure was tested independently: at $\alpha=0.005$ the two blocks satisfy
$\big\|\dot\phi_{z_1}/\barW(z_1)-\dot\phi_{z_2}/\barW(z_2)\big\|_\infty
=1.5\cdot10^{-5}$, and the ratio of the sup-norms of the two responses is $1.2997$
against the predicted $\barW(z_1)/\barW(z_2)=1.3$.
\end{remark}

\begin{proposition}[Regularity of the response in $z$]\label{prop:zregularity}
With $p$ as in \eqref{eq:varW}, the first order
response $\dot\phi_z=\barW(z)\Up$ satisfies
\[
\mathtt{var}_{p,\BV^1}\big(\dot\phi\big)\ \le\
\|\Up\|_{\BV^1}\ \mathtt{var}_{p,L^1}(W).
\]
In particular the response inherits exactly the regularity of the degree function: it
may be discontinuous in $z$, for instance for block graphons, while being smooth in
$x$.
\end{proposition}

\subsection{Higher order: the tree expansion}

Let $\mc T_j$ denote the set of isomorphism classes of finite rooted trees with $j$
edges, edges oriented away from the root. For $T\in\mc T_j$ with vertex set $V(T)$ and
root $\varnothing$, define the \emph{rooted homomorphism density}
\begin{equation}\label{eq:homdensity}
t_z(T,W):=\int_{[0,1]^{V(T)\setminus\{\varnothing\}}}
\prod_{(u,v)\in E(T)}W(z_u,z_v)\prod_{v\ne\varnothing}\dd z_v ,
\qquad z_\varnothing:=z .
\end{equation}
Thus $t_z(\bullet,W)=1$ for the one-vertex tree, $t_z(\text{edge},W)=\barW(z)$,
$t_z(P_2,W)=\int_0^1W(z,z')\barW(z')\dd z'$ for the path with two edges (the mean
degree of the neighbours of $z$) and $t_z(\Lambda_2,W)=\barW(z)^2$ for the cherry.

Rather than assume that $\phi^\star(\alpha)$ admits an asymptotic expansion to all
orders --- which we do not know --- we define its perturbative coefficients
unconditionally, by the recursion that such an expansion would have to satisfy.

\begin{proposition}[The formal coefficients are well defined]\label{prop:formal}
Set $\phi_{0,z}:=\rho$ and define inductively, for $1\le j\le r-2$,
\begin{equation}\label{eq:formaldef}
\phi_{j,z}:=(\Id-P)^{-1}\!\!\!\sum_{\substack{p\ge1,\ \ell\ge0\\
j_1,\dots,j_p\ge0\\ p+\ell+j_1+\cdots+j_p=j}}\!\!\!\frac{(-1)^p}{p!}
\Big[P\Big(A_{\phi_{j_1},z}\cdots A_{\phi_{j_p},z}\,\phi_{\ell,z}\Big)\Big]^{(p)} ,
\end{equation}
the sum being over ordered tuples. Then $\phi_{j}$ is well defined and
$\phi_{j}\in\W^{r-1-j}_0$ for $1\le j\le r-2$. We call $\phi_j$ the $j$-th
\emph{formal susceptibility} of the system.
\end{proposition}

\begin{proof}
The constraint $p\ge1$ forces $\ell\le j-1$ and $j_i\le j-1$, so the right-hand side
of \eqref{eq:formaldef} involves only $\phi_i$ with $i<j$ and the recursion is
well founded. Each summand is a $p$-th total derivative, hence has zero average on
$\T$, so $(\Id-P)^{-1}$ may be applied by Proposition~\ref{prop:gap}. For the
regularity, $A_{\phi_{j_i},z}\in\mc C^{r}$ by Lemma~\ref{lem:Abound}(i), which needs
only $\phi_{j_i}\in\W^0$; multiplying by $\phi_{\ell}\in\W^{r-1-\ell}$ and applying
$P$ stays in $W^{r-1-\ell,1}$ by Lemmas~\ref{lem:calculus} and \ref{lem:LY}, and the
$p$ derivatives land it in $W^{r-1-\ell-p,1}\subseteq W^{r-1-j,1}$ since
$\ell+p\le j$. Finally $(\Id-P)^{-1}$ preserves $W^{m,1}_0$ for $1\le m\le r-1$.
\end{proof}

\begin{theorem}[Theorem E: tree expansion of the formal susceptibilities]
\label{thm:tree}
For every $1\le j\le r-2$ there are functions $\Up_T\in W^{r-1-j,1}(\T)$,
$T\in\mc T_j$, with $\int_\T\Up_T=0$, depending only on $f$ and $h$ --- not on $W$
and not on $z$ --- such that the formal susceptibilities of
Proposition~\ref{prop:formal} satisfy
\begin{equation}\label{eq:tree}
\phi_{j,z}=\sum_{T\in\mc T_j}t_z(T,W)\,\Up_T \qquad\text{for a.e. }z .
\end{equation}
\end{theorem}

The formal susceptibilities are of interest only insofar as they describe the actual
equilibrium. The next statement says that they do, to the order we control.

\begin{proposition}[The first two formal susceptibilities are genuine]
\label{prop:formalgenuine}
Let $1\le k\le r-5$, which forces $r\ge6$. Then $\alpha\mapsto\phi^\star(\alpha)$ is
twice differentiable at $\alpha=0$ in $\W^{k}$ and
\begin{equation}\label{eq:expansion}
\phi^\star_z(\alpha)=\rho+\alpha\,\phi_{1,z}+\alpha^2\phi_{2,z}+o(\alpha^2)
\qquad\text{in }\W^k ,
\end{equation}
with $\phi_1,\phi_2$ the formal susceptibilities \eqref{eq:formaldef}. Whether
\eqref{eq:formaldef} describes the equilibrium beyond second order, and whether the
formal series converges, both remain open.
\end{proposition}

\begin{theorem}[Second order, explicitly]\label{thm:secondorder}
The second formal susceptibility of \eqref{eq:formaldef}, which by
Proposition~\ref{prop:formalgenuine} is the genuine second-order coefficient of
$\phi^\star(\alpha)$, is given by
\begin{equation}\label{eq:secondorder}
\phi_{2,z}=t_z(\Lambda_2,W)\ \Up_{\Lambda_2}\ +\ t_z(P_2,W)\ \Up_{P_2}
=\barW(z)^2\,\Up_{\Lambda_2}
+\Big(\int_0^1 W(z,z')\barW(z')\dd z'\Big)\Up_{P_2},
\end{equation}
where
\begin{align}
\Up_{\Lambda_2}&=(\Id-P)^{-1}\Big(\tfrac12\big[P\big(H_\rho^2\,\rho\big)\big]''
-\big[P\big(H_\rho\,\Up\big)\big]'\Big),\label{eq:UpLambda}\\
\Up_{P_2}&=-(\Id-P)^{-1}\Big[P\big(H_\Up\,\rho\big)\Big]' .\label{eq:UpP}
\end{align}
\end{theorem}

\begin{remark}[Numerical verification of \eqref{eq:secondorder}]\label{rem:numerics2}
Formula \eqref{eq:secondorder} was tested on the same example as
Remark~\ref{rem:numerics}, for which $\barW=(0.65,\,0.5)$ and, by
\eqref{eq:homdensity}, $t_z(\Lambda_2,W)=(0.4225,\,0.25)$ and
$t_z(P_2,W)=(0.40,\,0.31)$ on the two blocks. To isolate the second-order coefficient
without relying on the accuracy of $\Up$ we used the symmetric difference
\[
\frac{\phi^\star_z(\alpha)+\phi^\star_z(-\alpha)-2\rho}{2\alpha^{2}}
=\phi_{2,z}+o(1),
\]
which cancels the first-order term identically. The remainder here is only $o(1)$: all
that Proposition~\ref{prop:formalgenuine} provides is an expansion to second order with
an $o(\alpha^2)$ error, and no rate follows from it. The sup-norm discrepancy between
this quantity and the right-hand side of \eqref{eq:secondorder} was
$1.26\cdot10^{-4}$, $3.14\cdot10^{-5}$, $7.96\cdot10^{-6}$ for
$\alpha=0.08,\,0.04,\,0.02$ on the first block, i.e.\ it falls by a factor $4.00$ and
then $3.94$ as $\alpha$ is halved; the relative error at $\alpha=0.02$ is
$2.3\cdot10^{-4}$. The second block behaves identically. The observed quadratic decay
is \emph{faster} than the $o(1)$ that is guaranteed, and is what a fourth-order
expansion would give --- odd orders cancel in the symmetric difference, so the next
correction would be $\alpha^2\phi_4$ --- but we have not proved such an expansion, and
we record the rate as an empirical observation rather than as a consequence of the
theorem. At $\alpha=0.01$ the ratio degrades to $3.4$, the discretisation error having
become comparable to the quantity being measured; this is the resolution limit of the
scheme, not a failure of the formula. Since the two tree densities
$t_z(\Lambda_2,W)$ and $t_z(P_2,W)$ take
different values on the two blocks, and are not proportional to one another, this test
distinguishes the two universal profiles \eqref{eq:UpLambda}--\eqref{eq:UpP} rather
than only their sum.

For both this test and the one of Remark~\ref{rem:numerics}, densities were
represented on a uniform grid of $2^{10}$ points with periodic cubic interpolation;
transfer operators were applied through the two inverse branches of the lift, located
on a grid of $6\cdot10^{4}$ points; $\rho$ and $\phi^\star(\alpha)$ were obtained by
iterating to a residual below $10^{-10}$; derivatives were computed spectrally; and
$(\Id-P)^{-1}$ was evaluated by summing its Neumann series until the increment fell
below $10^{-14}$. The same discretisation was used for the order-three test of
Remark~\ref{rem:order3}, there with six blocks and a fixed-point residual below
$10^{-12}$, and for the exactly solvable check of Remark~\ref{rem:sanity}(iii). The
scripts producing all the tables quoted in this paper are available from the author.
\end{remark}

\subsection{Response to perturbations of the network}

\begin{theorem}[Theorem D: structural response]\label{thm:structural}
Let $V\in L^\infty([0,1]^2)$ with $\mathtt{var}_{p,L^1}(V)<\infty$ for the same
$p\in(0,1]$ as in \eqref{eq:varW}. Fix $\eps_0>0$ \emph{first}: by the
subadditivity of the $p$-variation established in \eqref{eq:varWeps} below, together
with $\|W_\eps\|_{\infty,1}\le\|W\|_{\infty,1}
+\eps_0\|V\|_{\infty,1}$, the kernels $W_\eps:=W+\eps V$, $|\eps|\le\eps_0$, satisfy
(A3) with bounds depending only on $\eps_0$ and on those of $W$ and $V$, so all the
constants of Sections~\ref{sec:spectral} and \ref{sec:differential} --- in particular
$\xi$, $K_r$, $\delta_0$, the auxiliary class $\mc V$ of \eqref{eq:Vclass}, the $M_k$
and a threshold $\alpha_2'\in(0,\alpha_2]$ --- may be chosen uniformly in
$|\eps|\le\eps_0$; see Remark~\ref{rem:epsuniform}. Then fix the coupling strength $\alpha$ with
$|\alpha|<\alpha_2'$ and let $\phi^{\star,\eps}$ be the fixed point of the STO built
from $W_\eps$ at that $\alpha$. Write $\mc K^V_{\alpha,\phi}$ for the operator
\eqref{eq:KK} with $A$ replaced by $A^V$, i.e.\ with $W$ replaced by $V$ in
\eqref{eq:A}; see Remark~\ref{rem:admissible} for what is, and what is not, required of
the perturbed kernel $W_\eps$. Then $\eps\mapsto\phi^{\star,\eps}$ is differentiable at
$\eps=0$ in $\W^k$, $1\le k\le r-3$, and
\begin{equation}\label{eq:structuralLR}
\partial_\eps\phi^{\star,\eps}\big|_{\eps=0}
=\alpha\ \mc R_{\alpha,\phi^\star}^{-1}\ \mc K^V_{\alpha,\phi^\star}\,\phi^\star .
\end{equation}
Moreover, as $\alpha\to0$,
\begin{equation}\label{eq:structuralLRsmall}
\partial_\eps\phi^{\star,\eps}_z\big|_{\eps=0}
=\alpha\,\overline V(z)\,\Up+O(\alpha^2)\qquad\text{in } \W^{k},
\qquad \overline V(z):=\int_0^1V(z,z')\dd z' .
\end{equation}
\end{theorem}

\begin{remark}[Relation with the Lipschitz continuity of the STO in $W$]
\label{rem:lipW}
Theorem~\ref{thm:structural} refines, at the level of the fixed point, the local
Lipschitz continuity of the STO itself in the graphon variable proved in
\cite[Thm.~3.3]{MPT2025}: there it is shown that, for graphons bounded by $R$ in
$L^\infty([0,1]^2)$ and $|\alpha|$ small depending on $R$, the map
$(W,\phi)\mapsto\mc F\phi$ is Lipschitz with respect to
$\|W-\widetilde W\|_{L^\infty}$ and the weak norm $\|\cdot\|_{``1"}$. Our hypothesis
(A3) and condition \eqref{eq:expcond} are the analogues of those assumptions, and
\eqref{eq:structuralLR} upgrades Lipschitz continuity to differentiability, in the
strong norm, for the fixed point rather than for one application of the operator.
\end{remark}

\begin{remark}[First-order optimal rewiring is degree-blind]\label{cor:control}
Fix $c\in L^\infty(\T)$ and a convex set $\mc P$ of admissible infinitesimal
rewirings, each $V\in\mc P$ satisfying the hypotheses of
Theorem~\ref{thm:structural}, and set
\[
C:=\alpha\int_\T c\,\Up\ \in\R .
\]
To leading order in $\alpha$, \eqref{eq:structuralLRsmall} gives
$\partial_\eps\int c\,\dd\nu^{\eps}_{z}=C\,\overline V(z)+O(\alpha^2)$, so that
optimising the response at a prescribed fibre $z$ over $V\in\mc P$ reduces to
optimising the single scalar $\overline V(z)$. The direction depends on the sign of
$C$, which involves both the sign of $\alpha$ and that of $\int_\T c\,\Up$: for $C>0$
one maximises $\overline V(z)$, for $C<0$ one minimises it, and for $C=0$ --- in
particular whenever $\Up\equiv0$, as in \eqref{eq:degenerate} --- every admissible
rewiring is indifferent at this order. In the first two cases the optimal value is
$C\sup_{\mc P}\overline V(z)$, respectively $C\inf_{\mc P}\overline V(z)$; we write
$\sup$ and $\inf$ because no compactness of $\mc P$ has been assumed, so an optimiser
need not exist. When it does, the optimal set is the face of $\mc P$ on which
$\overline V(z)$ is extremal, in general an infinite-dimensional set of rewirings. The
substantive point is independent of these signs: the first-order criterion is a linear
functional of $V$ factoring through the single number $\overline V(z)$, so all
rewirings with the same local degree change are equivalent at first order, and
selecting among them requires the second-order term of
Theorem~\ref{thm:secondorder}, in which the degrees of the neighbours appear.
\end{remark}

\subsection{Examples}\label{sec:examples}

\begin{example}[Erd\H{o}s--R\'enyi]\label{ex:ER}
$W\equiv p$, the limit graphon of the sequence considered in
\cite[Ex.~3.3]{MPT2025}. Then $\barW\equiv p$ and $t_z(T,W)=p^{|E(T)|}$, so the whole
expansion is $z$-independent:
$\phi^\star_z(\alpha)=\rho+\alpha p\,\Up+\alpha^2p^2(\Up_{\Lambda_2}+\Up_{P_2})+\cdots$.
One recovers the homogeneous theory of \cite{SelleyTanzi2021,Galatolo2022} with
effective coupling $\alpha p$ at first order.
\end{example}

\begin{example}[Stochastic block model / clustered network]\label{ex:SBM}
Let $[0,1]=\bigsqcup_{a=1}^{q}I_a$ and $W=\sum_{a,b}w_{ab}\mathbf 1_{I_a\times I_b}$,
as in \cite[Ex.~3.1]{MPT2025}. Then $\barW$ is constant on each block,
$\barW|_{I_a}=\sum_b w_{ab}|I_b|$, and the first-order response is piecewise constant:
with the convention of \eqref{eq:A}, $\barW$ is the \emph{in-strength} (the incoming
weighted degree), so for $w_{ab}\ge0$ and $\Up\neq0$ the blocks of larger in-strength
respond more in magnitude, while nodes inside a block are indistinguishable at first
order. At second order
$t_z(P_2,W)|_{I_a}=\sum_{b}w_{ab}|I_b|\,\barW|_{I_b}$ separates blocks according to the
strengths of their neighbouring blocks. Here $\mathtt{var}_{p,L^1}(W)<\infty$, so
Proposition~\ref{prop:zregularity} applies.
\end{example}

\begin{example}[Rank-one, configuration-model graphon: hubs]\label{ex:rank1}
Let $W(z,z')=g(z)g(z')$ with $g:[0,1]\to[0,\infty)$ bounded and Lipschitz, the graphon
limit of configuration models with prescribed degrees. Lipschitz continuity is what
makes \eqref{eq:varW} hold here: $\|W(z,\cdot)-W(z',\cdot)\|_{L^1}
=|g(z)-g(z')|\,\|g\|_{L^1}$, so $\mathtt{osc}_{L^1}(W,\omega,r)\le
2\,\Lip(g)\|g\|_{L^1}r$ and $\mathtt{var}_{p,L^1}(W)<\infty$ for every $p\le1$. Then
$\barW(z)=g(z)\|g\|_{L^1}$ and, by Theorem~\ref{thm:zerocoupling},
\[
\phi^\star_z(\alpha)=\rho+\alpha\,\|g\|_{L^1}\,g(z)\,\Up+o(\alpha).
\]
The susceptibility of a node is \emph{exactly proportional to its expected degree}:
hubs are the most responsive units, and the spread of $g$ is transferred verbatim to
the spread of susceptibilities across the network. (In the dense regime described by
graphons all degrees are of order $N$, so $g$ bounded away from $0$ and $\infty$ is
the natural setting; genuinely heavy-tailed degree sequences belong to the sparse
regime, which the graphon formalism does not describe.) At second order
\[
t_z(\Lambda_2,W)=g(z)^2\|g\|_{L^1}^2,
\qquad
t_z(P_2,W)=g(z)\,\|g\|_{L^1}\,\|g\|_{L^2}^2 ,
\]
so the two trees separate: the cherry contributes quadratically in $g(z)$ while the
path stays linear in $g(z)$.
\end{example}

\begin{example}[Spatially decaying interactions]\label{ex:decay}
Let $W(z,z')=\xi(|z-z'|)$ be the graphon of \cite[Ex.~3.2]{MPT2025}, and assume in
addition that $\xi:[0,1]\to[0,\infty)$ is Lipschitz and non-increasing, which is the
natural normalisation if $\xi$ is to be read as an interaction strength decaying with
distance. Splitting the integral at $z'=z$,
\begin{equation}\label{eq:decaydegree}
\barW(z)=\int_0^1\xi(|z-z'|)\,\dd z'=\int_0^{z}\xi(s)\,\dd s+\int_0^{1-z}\xi(s)\,\dd s ,
\end{equation}
so $\barW$ is Lipschitz, symmetric about $z=1/2$, and
$\barW'(z)=\xi(z)-\xi(1-z)$, which is $\ge0$ for $z<1/2$ and $\le0$ for $z>1/2$
because $\xi$ is non-increasing. Hence $\barW$ attains its maximum at the centre and
decreases monotonically towards either endpoint, and the first-order response is
heterogeneous: the heterogeneity is a pure boundary effect, nodes near the ends of the
chain having fewer neighbours. Without monotonicity of $\xi$ this conclusion may fail:
for $\xi(s)=s$ one finds $\barW(z)=\tfrac12\big(z^2+(1-z)^2\big)$, which is
\emph{minimal} at the centre. If instead
the kernel is taken translation invariant on the circle,
$W(z,z')=\xi(d_{\T}(z,z'))$, then $\barW$ is constant and, by induction using
\eqref{eq:graft} below, $t_z(T,W)$ is $z$-independent for \emph{every} rooted tree $T$.
Consequently the whole perturbative expansion is homogeneous across the network and
coincides with the all-to-all one with the same edge density: for translation-invariant
dense kernels the equilibrium state of the thermodynamic limit does not see the range
of the interaction at any order in perturbation theory.
\end{example}

\section{Uniform spectral estimates}\label{sec:spectral}

Assumption~\ref{ass:main} is in force throughout this section.

\begin{lemma}[Uniform family]\label{lem:uniform}
There are $\xi>1$ and $K_r>0$, depending only on $f,h,W$ and $\alpha_0$, such that for
all $|\alpha|<\alpha_0$, all $\phi\in\M_{1,\Leb}$ and a.e.\ $z$,
\[
\inf_{x\in\T}\big|F'_{\alpha,\phi,z}(x)\big|\ge\xi>1,\qquad
\|F_{\alpha,\phi,z}\|_{\mc C^{r}}\le K_r,\qquad
\deg F_{\alpha,\phi,z}=d .
\]
Moreover $\|F_{\alpha,\phi,z}-f\|_{\mc C^{j}}\le|\alpha|\,\|h\|_{\mc C^{j}}
\|W\|_{\infty,1}$ for $0\le j\le r$. For $z$ in the full measure set $A_W$ on which
$W(z,\cdot)\in L^1([0,1])$ we write
\begin{equation}\label{eq:Uz}
\mc U_z=\mc U_z(\alpha_0):=\big\{F_{\alpha,\nu,z}:\ |\alpha|<\alpha_0,\
\nu\in\M_{1,\Leb}\big\},\qquad \mc U:=\bigcup_{z\in A_W}\mc U_z ,
\end{equation}
and, for each fixed $\alpha$,
\begin{equation}\label{eq:Uzalpha}
\mc U_z(\alpha):=\big\{F_{\alpha,\nu,z}:\ \nu\in\M_{1,\Leb}\big\} .
\end{equation}
The bounds above hold for every $F\in\mc U$, hence uniformly over all the families
$\mc U_z$ and $\mc U_z(\alpha)$.

Finally we record an auxiliary class. Put
$\delta_0:=\alpha_0\|h\|_{\mc C^{r}}\|W\|_{\infty,1}$ and shrink $\alpha_0$, once and
for all, so that
\begin{equation}\label{eq:delta0}
\delta_0\le\tfrac14\big(\inf_\T f'-1\big) ;
\end{equation}
every $F\in\mc U$ then satisfies $\|F-f\|_{\mc C^{r}}\le\delta_0$ and
$\inf_\T F'\ge\inf_\T f'-\delta_0$, so that the constants of the display above may be
taken to be
\[
\xi:=\inf_\T f'-2\delta_0>1,\qquad K_r:=\|f\|_{\mc C^{r}}+2\delta_0 .
\]
With this choice set
\begin{equation}\label{eq:Vclass}
\mc V:=\Big\{F\in\mc C^{r}(\T,\T):\ \inf_\T F'\ge\xi,\
\|F-f\|_{\mc C^{r}}\le2\delta_0,\ \deg F=d\Big\} .
\end{equation}
Then $\mc U\subset\mc V$; every $F\in\mc V$ satisfies $\|F\|_{\mc C^{r}}\le K_r$ and
$\inf_\T F'\ge\xi$, which are the only properties of the maps used in
Lemmas~\ref{lem:LY} and \ref{lem:exact}; and $\mc V$ is \emph{convex}, since a convex
combination of two maps with derivative $\ge\xi$ has derivative $\ge\xi$ and the two
remaining conditions are respectively a ball and a homotopy invariant. The factor $2$
in \eqref{eq:Vclass} makes $\mc V$ a $\mc C^{r}$-\emph{neighbourhood} of $\mc U$:
\begin{equation}\label{eq:Vnbhd}
F\in\mc U,\quad Y\in\mc C^{r}(\T,\R),\quad\|Y\|_{\mc C^{r}}\le\delta_0
\qquad\Longrightarrow\qquad F+Y\in\mc V ,
\end{equation}
because then $\|F+Y-f\|_{\mc C^{r}}\le2\delta_0$,
$\inf_\T(F+Y)'\ge\inf_\T f'-\delta_0-\delta_0=\xi$, and $\deg(F+Y)=\deg F=d$ as $Y$ is
periodic. Without such a margin an $F$ lying on the boundary of the inequalities
defining $\mc V$ could be perturbed out of the class for every $s\neq0$ of one sign,
which is exactly what Lemma~\ref{lem:derivTO} must avoid.
\end{lemma}

\begin{remark}[Why the convex class $\mc V$ is needed]\label{rem:notconvex}
The family $\mc U_z$ is \emph{not} closed under the segments
$F\rightsquigarrow F+t\Delta$ joining two of its elements, so it cannot play the role
of $\mc V$. For a counterexample take $W\equiv1$ and
$h(x,y)=1+\sin(2\pi x)\cos(2\pi y)$, so that every element of $\mc U_z$ has the form
$f+\beta(1+q\,s)$ with $|\beta|<\alpha_0$, $q\in[-1,1]$ and $s(x)=\sin(2\pi x)$.
Choosing $a\in(0,\alpha_0)$ and two measures with $\cos$-moments $\pm m$, $0<m<1$,
the maps $F_{a,\nu_+,z}=f+a(1+ms)$ and $F_{-a,\nu_-,z}=f-a(1-ms)$ both lie in
$\mc U_z$, while their midpoint $f+ams$ does not: an element $f+\beta(1+qs)$ with no
constant term forces $\beta=0$, hence no sinusoidal term either. The obstruction is
only the sign of $\alpha$, and it is harmless, for two reasons. First, each
$\mc U_z(\alpha)$ \emph{is} convex, by linearity of $\nu\mapsto A_{\nu,z}$: for
$\nu_t=(1-t)\nu_0+t\nu_1$ one has
$(1-t)F_{\alpha,\nu_0,z}+tF_{\alpha,\nu_1,z}=F_{\alpha,\nu_t,z}$. Second, the two
families of segments that actually occur below are
\[
\big\{F_{t\alpha,\phi,z}\big\}_{t\in[0,1]}\subset\mc U_z
\qquad\text{and}\qquad
\big\{F_{\alpha,\phi+tv,z}\big\}_{t\in[0,1]}\subset\mc U_z(\alpha) ,
\]
the first because $|t\alpha|\le|\alpha|<\alpha_0$ and the second because
$\phi+tv$ is a convex combination of the probability densities $\phi$ and $\phi+v$.
In any case every such segment lies in $\mc V$, which is all that the
Lasota--Yorke and exact perturbative estimates require, since those use only the
uniform bounds of this lemma. By contrast the memory-loss estimate of
Lemma~\ref{lem:MLMPT} is \emph{not} available on $\mc V$, and
Proposition~\ref{prop:gap} is therefore stated on $\mc U_z(\alpha)$ only.
\end{remark}

\begin{remark}[Why the fibre and the coupling must both be kept fixed]
\label{rem:fixedfibre}
The memory-loss estimate we import, Lemma~\ref{lem:MLMPT}, controls compositions
$F_{\nu_n,z}\circ\cdots\circ F_{\nu_1,z}$ in which the \emph{same} $z$ and the
\emph{same} $\alpha$ occur in every factor, only the measures $\nu_i$ being allowed to
vary; it says nothing about a composition whose factors come from different fibres, nor
about one in which the coupling strength changes from factor to factor. Accordingly
Proposition~\ref{prop:gap} is stated for sequences inside a single
$\mc U_z(\alpha)$, with constants independent of both $z$ and $\alpha$. This is
exactly what is needed, and nothing more is claimed: the operators
$\LL_{\alpha,\phi}$ act fibrewise at a fixed coupling, so a composition
$\LL_{\alpha,\phi^{(n)}}\cdots\LL_{\alpha,\phi^{(1)}}$ evaluated at the fibre $z$ is
$L_{\alpha,\phi^{(n)},z}\cdots L_{\alpha,\phi^{(1)},z}$, a composition within
$\mc U_z(\alpha)$. No composition in this paper mixes fibres or coupling strengths.
\end{remark}

\begin{proof}
Since $\phi_{z'}$ is a probability density, Remark~\ref{rem:BVW11}(ii) gives
$\|H_{\phi_{z'}}\|_{\mc C^j}\le\|h\|_{\mc C^j}$, whence
\[
\|A_{\phi,z}\|_{\mc C^j}\le\|h\|_{\mc C^j}\,\|W(z,\cdot)\|_{L^1}
\le\|h\|_{\mc C^j}\,\|W\|_{\infty,1} .
\]
The three claims follow from \eqref{eq:F}, \eqref{eq:expcond}, and the homotopy
invariance of the degree. This refines Lemma~\ref{lem:uniformMPT} from $\mc C^3$ to
$\mc C^{r}$.
\end{proof}

\begin{lemma}[Bounds on the mean field]\label{lem:Abound}
For every $0\le j\le r$ and all $v,w\in\W^0$ and a.e.\ $z$:
\begin{itemize}
\item[(i)] $\|A_{v,z}\|_{\mc C^{j}}\le \|h\|_{\mc C^j}\,\|W\|_{\infty,1}\,\|v\|_{0}$;
\item[(ii)] $\|A_{v,z}-A_{w,z}\|_{\mc C^j}
\le \|h\|_{\mc C^j}\|W\|_{\infty,1}\|v-w\|_0$;
\item[(iii)] if $W\in L^\infty([0,1]^2)$ then
$\|A_{v,z}\|_{\mc C^j}\le \|h\|_{\mc C^j}\|W\|_\infty\|v\|_{L^1([0,1]\times\T)}$.
\end{itemize}
In particular $v\mapsto A_v$ is a bounded linear map from the weakest space $\W^0$ into
$L^\infty([0,1],\mc C^{r}(\T))$: the mean field is \emph{smoothing}.
\end{lemma}

\begin{proof}
Immediate from \eqref{eq:A}, Remark~\ref{rem:BVW11}(ii) and H\"older's inequality; (ii)
follows from (i) by linearity of $v\mapsto A_v$.
\end{proof}

\begin{lemma}[Uniform Lasota--Yorke in $W^{k,1}$]\label{lem:LY}
Let $1\le k\le r-1$ and $\lambda:=\xi^{-1}<1$. There are $A_k,B_k>0$, depending only on
$K_r,\xi,k$, such that for every $n\ge1$, every $F_1,\dots,F_n\in\mc V$ --- in
particular for every $F_1,\dots,F_n\in\mc U$ --- and every
$\psi\in W^{k,1}(\T)$, writing $L(n):=L_{F_n}\cdots L_{F_1}$,
\begin{align}
\|L(n)\psi\|_{W^{k-1,1}}&\le A_k\|\psi\|_{W^{k-1,1}},\label{eq:LY1}\\
\|L(n)\psi\|_{W^{k,1}}&\le \lambda^{kn}\|\psi\|_{W^{k,1}}
+B_k\|\psi\|_{W^{k-1,1}}.\label{eq:LY2}
\end{align}
\end{lemma}

\begin{remark}[Why Sobolev and not $\mc C^k$]\label{rem:nodegree}
The contraction factor in \eqref{eq:LY2} is $\lambda^{k}=\xi^{-k}$, with \emph{no}
factor equal to the number $d$ of inverse branches. This is specific to the $L^1$-based
scale and follows from $\|L_Fu\|_{L^1}\le\|u\|_{L^1}$; see the proof in
Appendix~\ref{app:longproofs}. In the $\mc C^k$ scale used in \cite{SelleyTanzi2021}
the analogous constant is $d/\xi^{k}$, which forces an additional hypothesis relating
the degree to the expansion, such as \cite[Eq.~(8)]{SelleyTanzi2021} or the condition
$d<\sigma^2$ of \cite{Balint2018}. Our statements require no such hypothesis.
\end{remark}

\begin{proposition}[Uniform fibrewise spectral gap]\label{prop:gap}
Let $1\le k\le r-1$. There are $C_\star\ge1$ and $\theta\in(0,1)$, depending only on
$f,h,W,\alpha_0,k$ and \emph{not} on $z$ or $\alpha$, such that for every
$|\alpha|<\alpha_0$, every $z\in A_W$, every $n\ge0$, every
$F_1,\dots,F_n\in\mc U_z(\alpha)$ --- the same fibre \emph{and the same coupling
strength} in every factor --- and every $\psi\in W^{k,1}(\T)$ with
$\int_\T\psi=0$,
\begin{equation}\label{eq:gap}
\big\|L_{F_n}\cdots L_{F_1}\psi\big\|_{W^{k,1}}
\le C_\star\,\theta^{\,n}\,\|\psi\|_{W^{k,1}} .
\end{equation}
Consequently, for every $|\alpha|<\alpha_0$ and every sequence
$\{\phi^{(i)}\}\subset\M_{1,\Leb}$, the fibrewise operators satisfy
$\big\|\LL_{\alpha,\phi^{(n)}}\cdots\LL_{\alpha,\phi^{(1)}}v\big\|_{k}
\le C_\star\theta^n\|v\|_{k}$ for $v\in\W^k_0$; in particular $\Id-\LL_{\alpha,\phi}$
is invertible on $\W^k_0$ with
\begin{equation}\label{eq:resolvent0}
\big\|(\Id-\LL_{\alpha,\phi})^{-1}\big\|_{\W^k_0\to\W^k_0}\le
\frac{C_\star}{1-\theta}=:C_R .
\end{equation}
\end{proposition}

The case $k=1$ of \eqref{eq:gap} is \emph{exactly} Lemma~\ref{lem:MLMPT}, by
Remark~\ref{rem:BVW11}(i); the cases $k\ge2$ follow by a two-norm bootstrap. The proof
is in Appendix~\ref{app:longproofs}.

\begin{remark}[The case $k=0$]\label{rem:nok0}
No estimate of the form $\|L(n)\psi\|_{L^1}\le C\theta^n\|\psi\|_{L^1}$ is available,
nor is one needed: transfer operators of expanding maps do not have a spectral gap on
$L^1$. All statements below involving the resolvent are therefore restricted to
$k\ge1$. The weakest space $\W^0$ enters only as the \emph{domain} of the smoothing
operator $\mc K$, where no gap is required.
\end{remark}

\begin{proposition}[A priori regularity of the fixed point]\label{prop:apriori}
Let $1\le k\le r-1$. There is $M_k>0$, depending only on $f,h,W,\alpha_0,k$, such that
for all $|\alpha|<\alpha_0$
\[
\big\|\phi^\star(\alpha)\big\|_{k}=\esssup_z\|\phi^\star_z(\alpha)\|_{W^{k,1}}\le M_k .
\]
Moreover the sets $\Delta^{k}_M$ of \eqref{eq:ball} are forward invariant under
$\STO_\alpha$ for $M_i$ large enough, and contain $\phi^\star(\alpha)$ for every
$|\alpha|<\alpha_0$.
\end{proposition}

\begin{proof}
By Theorem~\ref{thm:MPT}, for a.e.\ $z$ the density $\phi^\star_z(\alpha)$ is the
invariant density of the map $F_{\alpha,\phi^\star,z}\in\mc U$, which is
$\mc C^{r}$ and uniformly expanding by Lemma~\ref{lem:uniform}. By the classical
regularity theorem for invariant densities of smooth expanding maps
(Krzy\.zewski--Szlenk \cite{KrzyzewskiSzlenk1969}, see also
\cite{Krzyzewski1977} and \cite[Cor.~3.1]{MPT2025}), $\phi^\star_z(\alpha)\in
\mc C^{r-1}(\T)$; in particular $\|\phi^\star_z(\alpha)\|_{W^{k,1}}<\infty$ for
$k\le r-1$. This \emph{qualitative} information makes the following argument
non-vacuous. Applying \eqref{eq:LY2} with $n$ copies of $F_{\alpha,\phi^\star,z}$ to
$\psi=\phi^\star_z$, which is fixed by $L_{\alpha,\phi^\star,z}$, gives
\[
\|\phi^\star_z\|_{W^{k,1}}\le\lambda^{kn}\|\phi^\star_z\|_{W^{k,1}}
+B_k\|\phi^\star_z\|_{W^{k-1,1}} ,
\]
and letting $n\to\infty$,
$\|\phi^\star_z\|_{W^{k,1}}\le B_k\|\phi^\star_z\|_{W^{k-1,1}}$. Since
$\|\phi^\star_z\|_{L^1}=1$, a downward induction on $k$ produces constants $M_k$
depending only on $K_r,\xi,k$, hence uniform in $\alpha$ and $z$.

For the forward invariance, \eqref{eq:LY2} with $n=1$ gives
$\|\STO_\alpha\phi\|_{i}\le\lambda^{i}\|\phi\|_{i}+B_i\|\phi\|_{i-1}$ for
$1\le i\le k$, so choosing $M_0=1$ and
$M_i\ge B_i M_{i-1}(1-\lambda^{i})^{-1}$ inductively makes $\Delta^k_M$ invariant.
Enlarging $M_i$ if necessary we may also assume $M_i\ge\sup_\alpha
\|\phi^\star(\alpha)\|_i$.
\end{proof}

\section{The differential of the STO and its resolvent}\label{sec:differential}

\subsection{Exact perturbative identities}

Everything in this paper rests on the following two \emph{exact} identities. Their
proof, a duality computation followed by integration by parts, is in
Appendix~\ref{app:calculus}.

\begin{lemma}[Exact first- and second-order identities]\label{lem:exact}
Let $F,F+\Delta\in\mc V$; then $\Delta\in\mc C^{r}(\T,\R)$ automatically, and the whole
segment $\{F+t\Delta\}_{t\in[0,1]}$ lies in $\mc V$ by its convexity. Then, as elements
of $W^{m,1}(\T)$, and with constants involving only the lower norms of $\Delta$
displayed below:
\begin{itemize}
\item[(i)] if $m+1\le r-1$ and $u\in W^{m+1,1}(\T)$,
\begin{equation}\label{eq:exact1}
L_{F+\Delta}u-L_Fu=-\int_0^1\Big[L_{F+t\Delta}\big(\Delta\,u\big)\Big]'\,\dd t ,
\end{equation}
and there is $C>0$ depending only on $K_r,\xi,m$ with
\begin{equation}\label{eq:est1}
\big\|L_{F+\Delta}u-L_Fu\big\|_{W^{m,1}}
\le C\,\|\Delta\|_{\mc C^{m+1}}\,\|u\|_{W^{m+1,1}} ;
\end{equation}
\item[(ii)] if $m+2\le r-1$ and $u\in W^{m+2,1}(\T)$,
\begin{equation}\label{eq:exact2}
L_{F+\Delta}u-L_Fu+\Big[L_{F}\big(\Delta\,u\big)\Big]'
=\int_0^1(1-t)\Big[L_{F+t\Delta}\big(\Delta^2 u\big)\Big]''\,\dd t ,
\end{equation}
and
\begin{equation}\label{eq:est2}
\Big\|L_{F+\Delta}u-L_Fu+\big[L_{F}(\Delta u)\big]'\Big\|_{W^{m,1}}
\le C\,\|\Delta\|^2_{\mc C^{m+2}}\,\|u\|_{W^{m+2,1}} .
\end{equation}
\end{itemize}
The first-order identity therefore costs one derivative and the second-order one two;
this is why \eqref{eq:est1} is available one level higher in the scale than
\eqref{eq:est2}, a distinction used in the proof of Theorem~\ref{thm:stability}.
\end{lemma}

\begin{remark}[Loss of derivatives]\label{rem:loss}
Identity \eqref{eq:exact2} shows that the first-order term
$-[L_F(\Delta u)]'$ lies in $W^{m+1,1}$ while the remainder is controlled only in
$W^{m,1}$ and requires two more derivatives of $u$. This two-derivative loss is
intrinsic and is the reason why all our statements come with an explicit budget of
indices; it is also the reason why the implicit function theorem cannot be applied
directly to $\Psi(\alpha,\phi)=\STO_\alpha\phi-\phi$. The same phenomenon occurs in
the homogeneous theory: the differential of an STO is Fr\'echet differentiable only
with convergence in a weaker norm, which obstructs naive iteration of the differential,
as observed in \cite[Rem.~3]{CastorriniGalatoloTanzi2025}. Our way around it is the
uniform a priori bound of Theorem~\ref{thm:stability} one level higher in the scale,
used in the proof of Theorem~\ref{thm:LR}.
\end{remark}

\begin{lemma}[Derivative of $L_F$ in $F$]\label{lem:derivTO}
Let $F\in\mc U$, let $X\in\mc C^{r}(\T,\R)$ and let $m+2\le r-1$. Put
$s_X:=\delta_0/(1+\|X\|_{\mc C^{r}})$, so that $F+sX\in\mc V$ for every $|s|\le s_X$ by
\eqref{eq:Vnbhd}; this is what makes Lemma~\ref{lem:exact} applicable to the
perturbation $\Delta=sX$. Then for every $u\in W^{m+2,1}(\T)$,
\begin{equation}\label{eq:derivTO}
\lim_{s\to0}\ \left\|\frac{L_{F+sX}u-L_Fu}{s}
+\Big[L_F\big(X\,u\big)\Big]'\right\|_{W^{m,1}}=0 ,
\end{equation}
with error $O(s)\|X\|^2_{\mc C^{m+2}}\|u\|_{W^{m+2,1}}$ for $|s|\le s_X$. In the two
applications made below, at \eqref{eq:alphaderiv} and in the proof of
Theorem~\ref{thm:structural}, the direction $X$ is respectively $A_{\phi,z}$ and
$\alpha A^V_{\phi,z}$, both of class $\mc C^{r}$ by Lemma~\ref{lem:Abound}(i).
\end{lemma}

\begin{proof}
Apply \eqref{eq:exact2} with $\Delta=sX$ and divide by $s$.
\end{proof}

\begin{remark}[Cross-checks with the literature]\label{rem:crosscheck}
Formula \eqref{eq:derivTO} reproduces the known response formulae, which is a
nontrivial test of our conventions.

\emph{(a) Galatolo's model $F=\Phi_{\delta,\rho}\circ f$ with
$\Phi_{\delta,\rho}=\Id+\delta H_\rho$.} Here $X=H_\rho\circ f$, and since
$L_f\big((g\circ f)\,u\big)=g\,L_fu$ we get
$\partial_\delta L_{F_\delta}\rho|_{0}=-\big[H_\rho\,P\rho\big]'=-(H_\rho\rho)'$,
so the response is $-(\Id-P)^{-1}(\rho H_\rho)'$, which is exactly
\cite[Eq.~(60)]{Galatolo2022}.

\emph{(b) S\'elley--Tanzi's model $F=f\circ\Phi_{t,\rho}$ with
$\Phi_{t,\rho}=\Id+tH_\rho$.} Here $X=(f'\circ\Phi_0)H_\rho=f'H_\rho$, so
\eqref{eq:derivTO} gives $-\big[P(f'H_\rho\rho)\big]'$. Using the standard identity
$(L_Fu)'=L_F(u'/F')-L_F\big(F''(F')^{-2}u\big)$ with $F=f$ and $u=f'w$ one finds
\[
\big[P(f'w)\big]'=P\Big(\frac{(f'w)'}{f'}\Big)-P\Big(\frac{f''}{f'}w\Big)=P(w') ,
\]
so that $-\big[P(f'H_\rho\rho)\big]'=-P\big[(H_\rho\rho)'\big]$ and the response is
$-(\Id-P)^{-1}P\big[(\rho H_\rho)'\big]$, which is exactly
\cite[Prop.~1.4]{SelleyTanzi2021}.

\emph{(c) Our model $F=f+\alpha A$.} Neither reduction applies, and the derivative
carries a transfer operator \emph{inside} the derivative,
$-\big[P(A\rho)\big]'$, whence the shape of \eqref{eq:Upsilon}. Observe that (a) and
(b) themselves differ by one factor of $P$, consistently with the different order of
composition.
\end{remark}

\begin{remark}[Sanity checks, and a degeneracy to be avoided]\label{rem:sanity}
\emph{(i) Rigid translations.} Take $h\equiv1$, so that $H_\rho\equiv1$ and
$F_{\alpha,\phi,z}=f+\alpha\barW(z)$, a rigid translation of $f$ whose size depends on
the fibre. Then $\Up=-(\Id-P)^{-1}\big[P\rho\big]'=-(\Id-P)^{-1}\rho'$, the classical
response of an expanding map to a rigid translation. If moreover $f(x)=dx$, then
$\rho\equiv1$ and $\Up=0$, as it must be, since $dx+c$ is again a linear expanding map
with uniform invariant density. Theorem~\ref{thm:zerocoupling} therefore predicts no
first-order response for this family, for any graphon.

\emph{(ii) A degenerate family.} By \eqref{eq:Upsilon}, $\Up$ vanishes exactly when
$P(H_\rho\rho)$ is constant, and this occurs for a family wide enough to deserve a
warning. If $f(x)=dx$ then $\rho\equiv1$ and $P\mathbf 1=\mathbf 1$, so
$H_\rho(x)=\int_\T h(x,y)\dd y$ and $P(H_\rho\rho)=P H_\rho$; consequently
\begin{equation}\label{eq:degenerate}
f(x)=dx\quad\text{and}\quad \int_\T h(x,y)\,\dd y\ \text{ independent of }x
\qquad\Longrightarrow\qquad \Up\equiv0 .
\end{equation}
In particular every difference coupling $h(x,y)=\eta(y-x)$ with $\int_\T\eta=0$ --- the
most common choice in the literature --- is degenerate when paired with a linear map.
Both terms of the coupling used in our own experiments below have vanishing
$y$-average, so that coupling too must be paired with a nonlinear map. The
first-order response then vanishes identically, for every graphon, and the
degree--susceptibility law of Theorem~\ref{thm:zerocoupling} is invisible. Such pairs
$(f,h)$ must be avoided when testing the theory: either $f$ must be genuinely
nonlinear, as in Remark~\ref{rem:numerics}, or $\int_\T h(\cdot,y)\dd y$ must be
nonconstant.

\emph{(iii) An exactly solvable non-degenerate case.} Take $f(x)=2x\bmod1$ and
$h(x,y)=\cos(4\pi x)$, so that $\rho\equiv1$ and $H_\rho(x)=\cos(4\pi x)$. For the
doubling map $P\cos(4\pi x)=\cos(2\pi x)$ and $P\sin(2\pi x)=0$, so the Neumann series
for $(\Id-P)^{-1}$ terminates after one term and
\[
\Up(x)=-(\Id-P)^{-1}\big[\cos(2\pi x)\big]'=2\pi\sin(2\pi x),
\qquad
\partial_\alpha\phi^\star_z\big|_{0}=2\pi\,\barW(z)\,\sin(2\pi x) .
\]
The numerical scheme described in Remark~\ref{rem:numerics} reproduces this closed form
to $10^{-7}$ in sup norm, which validates the scheme against an exact answer
independently of Theorem~\ref{thm:zerocoupling}. This pair is a reference check and not
a model: $h$ does not depend on $y$, so there is no self-consistent feedback and the
mean field degenerates to a prescribed fibre-dependent forcing.
\end{remark}

\subsection{Proof of Theorem~\ref{thm:diff}}

\begin{proof}[Proof of Theorem~\ref{thm:diff}]
Fix $1\le k\le r-3$ and $\phi,\phi+v\in\Delta^{k+2}$ with
$\|\phi\|_{k+2},\|\phi+v\|_{k+2}\le M$. By linearity of $w\mapsto A_w$,
\[
F_{\alpha,\phi+v,z}=F_{\alpha,\phi,z}+\Delta_z,\qquad \Delta_z:=\alpha A_{v,z},
\]
and by Lemma~\ref{lem:Abound}(i),
$\|\Delta_z\|_{\mc C^{j}}\le|\alpha|\,\|h\|_{\mc C^{j}}\|W\|_{\infty,1}\|v\|_0$ for
$j\le r$. Write $L_z:=L_{\alpha,\phi,z}$. Applying \eqref{eq:exact2} with $F=F_{\alpha,\phi,z}$, $\Delta=\Delta_z$,
$u=\phi_z+v_z$ and $m=k$ (legitimate since $k+2\le r-1$):
\[
\big(\STO_\alpha(\phi+v)\big)_z
=L_z(\phi_z+v_z)-\Big[L_z\big(\Delta_z(\phi_z+v_z)\big)\Big]'+\mc E_z ,
\]
with $\|\mc E_z\|_{W^{k,1}}\le C\|\Delta_z\|^2_{\mc C^{k+2}}\|\phi_z+v_z\|_{W^{k+2,1}}
\le C'\alpha^2\|v\|_0^2\,M$. Expanding the middle term and recalling
\eqref{eq:LL}--\eqref{eq:KK},
\[
\big(\STO_\alpha(\phi+v)\big)_z=\big(\STO_\alpha\phi\big)_z
+\big(\LL_{\alpha,\phi}v\big)_z+\alpha\big(\mc K_{\alpha,\phi}v\big)_z
-\alpha\Big[L_z\big(A_{v,z}v_z\big)\Big]'+\mc E_z .
\]
By Lemma~\ref{lem:LY}, $L_z$ is bounded on $W^{k+1,1}$, and $W^{k+1,1}(\T)$ is a module
over $\mc C^{k+1}(\T)$ (Lemma~\ref{lem:calculus}), so
\begin{align*}
\Big\|\alpha\big[L_z(A_{v,z}v_z)\big]'\Big\|_{W^{k,1}}
&\le|\alpha|\,\big\|L_z(A_{v,z}v_z)\big\|_{W^{k+1,1}}
\le C''|\alpha|\,\|A_{v,z}\|_{\mc C^{k+1}}\|v_z\|_{W^{k+1,1}}\\
&\le C'''|\alpha|\,\|v\|_0\,\|v\|_{k+1} .
\end{align*}
Since $\|v\|_0\le\|v\|_{k+1}$ and $|\alpha|<\alpha_0$, collecting the two error terms
gives \eqref{eq:diffSTO}. The last assertion is Lemma~\ref{lem:Kphi}(ii).
\end{proof}

\begin{lemma}[The $\alpha$-derivative is $\mc K$ evaluated at the state]\label{lem:Kphi}
Let $\phi$ be \emph{fixed}.
\begin{itemize}
\item[(i)] If $1\le k\le r-2$ and $\phi\in\Delta^{k+1}$, then for all $\alpha,\beta$
with $|\alpha|,|\beta|<\alpha_0$,
\begin{equation}\label{eq:alphaLip}
\big\|\STO_\alpha\phi-\STO_\beta\phi\big\|_{k}\le
C\,|\alpha-\beta|\;\|\phi\|_{k+1}, \qquad C=C(f,h,W,\alpha_0,k) .
\end{equation}
\item[(ii)] If $1\le k\le r-3$ and $\phi\in\Delta^{k+2}$, then
$\alpha\mapsto\STO_\alpha\phi$ is differentiable from $\R$ to $\W^{k}$ and
\begin{equation}\label{eq:alphaderiv}
\partial_\alpha\big(\STO_\alpha\phi\big)=\mc K_{\alpha,\phi}\,\phi .
\end{equation}
\end{itemize}
\end{lemma}

\begin{proof}
With $\phi$ frozen, $F_{\alpha,\phi,z}-F_{\beta,\phi,z}=(\alpha-\beta)A_{\phi,z}$
\emph{exactly}, since $\alpha\mapsto F_{\alpha,\phi,z}$ is affine. Estimate
\eqref{eq:alphaLip} is then Lemma~\ref{lem:exact}(i), i.e.\ \eqref{eq:est1}, with
$\Delta=(\alpha-\beta)A_{\phi,z}$, $u=\phi_z$ and $m=k$, combined with
Lemma~\ref{lem:Abound}(i) and $\|\phi\|_0=1$; it therefore costs only one derivative
and is available for $k\le r-2$. Identity \eqref{eq:alphaderiv} is
Lemma~\ref{lem:derivTO} with $X=A_{\phi,z}$ and $u=\phi_z$, which rests on
Lemma~\ref{lem:exact}(ii) and hence needs $k\le r-3$.
\end{proof}

\begin{remark}
Lemma~\ref{lem:Kphi} is a genuine simplification with respect to the general theory,
and it holds precisely because $v\mapsto A_v$ is linear: the external parameter
$\alpha$ and the internal self-consistency enter through the same bilinear expression
$\alpha A_\phi$. It is what makes the right-hand side of \eqref{eq:LRgeneral} equal to
$\mc K\phi^\star$ rather than to an unrelated derivative operator.
\end{remark}

\subsection{Boundedness of $\mc K$ and invertibility of $\Id-\LL-\alpha\mc K$}

\begin{proposition}[$\mc K$ is bounded and smoothing]\label{prop:Kbounded}
Let $1\le k\le r-2$ and $\phi\in\Delta^{k+1}$. Then $\mc K_{\alpha,\phi}$ maps $\W^0$
boundedly into $\W^{k}_0$, with
\begin{equation}\label{eq:Kbound}
\big\|\mc K_{\alpha,\phi}v\big\|_{k}\le C_{\mc K}\,\|\phi\|_{k+1}\,\|v\|_{0}
\qquad\text{for all }v\in\W^0 ,
\end{equation}
$C_{\mc K}=C_{\mc K}(f,h,W,\alpha_0,k)$.
\end{proposition}

\begin{proof}
By Lemma~\ref{lem:LY}, $L_{\alpha,\phi,z}$ is bounded on $W^{k+1,1}(\T)$ uniformly over
$\mc U$ (here $k+1\le r-1$), and by Lemma~\ref{lem:calculus},
$\|A_{v,z}\phi_z\|_{W^{k+1,1}}\le c\,\|A_{v,z}\|_{\mc C^{k+1}}\|\phi_z\|_{W^{k+1,1}}$.
Hence
\[
\|(\mc K_{\alpha,\phi} v)_z\|_{W^{k,1}}
=\big\|[L_{\alpha,\phi,z}(A_{v,z}\phi_z)]'\big\|_{W^{k,1}}
\le\big\|L_{\alpha,\phi,z}(A_{v,z}\phi_z)\big\|_{W^{k+1,1}}
\le c'\,\|A_{v,z}\|_{\mc C^{k+1}}\,\|\phi_z\|_{W^{k+1,1}} ,
\]
and Lemma~\ref{lem:Abound}(i) bounds $\|A_{v,z}\|_{\mc C^{k+1}}$ by
$\|h\|_{\mc C^{k+1}}\|W\|_{\infty,1}\|v\|_0$. Taking $\esssup_z$ gives
\eqref{eq:Kbound}; finally $\int_\T(\mc Kv)_z=0$ since each fibre is a total
derivative.
\end{proof}

The content of \eqref{eq:Kbound} deserves emphasis: $\mc K$ is bounded from the
\emph{weakest} space $\W^0$ into the \emph{strong} space $\W^k$. In the language of
\cite{Galatolo2022} the small-perturbation hypothesis (LR2) holds here in the strongest
possible mixed norm, with no loss of derivatives. This is a consequence of the
smoothing character of the mean field (Remark~\ref{rem:BVW11}(ii)) and of the fact that
the derivative in \eqref{eq:KK} is applied to an object which has already gained $k+1$
derivatives.

\begin{proposition}[Resolvent]\label{prop:resolvent}
Let $1\le k\le r-2$ and $M\ge M_{k+1}$. There is $\alpha_2\in(0,\alpha_0]$ such that
for all $|\alpha|<\alpha_2$ and all $\phi\in\Delta^{k+1}$ with $\|\phi\|_{k+1}\le M$,
the operator
\[
\mc R_{\alpha,\phi}:=\Id-\LL_{\alpha,\phi}-\alpha\,\mc K_{\alpha,\phi}
\]
is invertible on $\W^k_0$, with
\begin{equation}\label{eq:resolventbound}
\big\|\mc R_{\alpha,\phi}^{-1}\big\|_{\W^k_0\to\W^k_0}\le
\frac{C_R}{1-|\alpha|\,C_R\,C_{\mc K}M}\ ,
\end{equation}
and
\begin{equation}\label{eq:neumann}
\mc R_{\alpha,\phi}^{-1}=\sum_{m\ge0}\Big[(\Id-\LL_{\alpha,\phi})^{-1}\,
\alpha\,\mc K_{\alpha,\phi}\Big]^{m}(\Id-\LL_{\alpha,\phi})^{-1} .
\end{equation}
\end{proposition}

\begin{proof}
By Proposition~\ref{prop:gap}, $(\Id-\LL)^{-1}$ is bounded on $\W^k_0$ with norm
$\le C_R$. By Proposition~\ref{prop:Kbounded} and $\|\cdot\|_0\le\|\cdot\|_k$,
$\mc K$ maps $\W^k_0$ into itself with
$\|\alpha\mc K\|_{\W^k_0\to\W^k_0}\le|\alpha|C_{\mc K}M$. Choose $\alpha_2$ with
$\alpha_2 C_RC_{\mc K}M<1$ and write
$\mc R=(\Id-\LL)\big[\Id-(\Id-\LL)^{-1}\alpha\mc K\big]$; the bracket is invertible by
Neumann series, giving \eqref{eq:neumann} and \eqref{eq:resolventbound}.
\end{proof}

We emphasise what Proposition~\ref{prop:resolvent} does \emph{not} assert: it gives no
continuity of $\alpha\mapsto\mc R_{\alpha,\phi}^{-1}$ in the operator norm of
$\mc L(\W^k_0)$, and no statement of that kind is used anywhere in this paper. The
obstruction is the one-derivative loss in \eqref{eq:est1}: the difference
$\mc R_0-\mc R_\alpha$ is $O(\alpha)$ only as an operator
$\W^{k+1}_0\to\W^{k}_0$, and an $O(\alpha)$ bound in a mixed norm does not yield one
in the weaker norm alone. What the resolvent identity does give is the following
statement, which we shall apply only to explicitly identified vectors whose extra
derivative has been checked.

\begin{lemma}[Resolvent convergence on a regular vector]\label{lem:resvector}
Let $1\le k\le r-2$ and $M\ge M_{k+1}$, let $\alpha\mapsto\phi(\alpha)\in\Delta^{k+1}$
satisfy $\|\phi(\alpha)\|_{k+1}\le M$ for $|\alpha|<\alpha_2$, and write
$\mc R_\alpha:=\mc R_{\alpha,\phi(\alpha)}$ and $\mc R_0:=\Id-P$. There is $C>0$,
depending only on $f,h,W,\alpha_0,k,M$, such that
\begin{equation}\label{eq:resmixed}
\big\|\big(\mc R_\alpha^{-1}-\mc R_0^{-1}\big)w\big\|_{k}
\le C\,|\alpha|\,\|w\|_{k+1}
\qquad\text{for every } w\in\W^{k+1}_0 .
\end{equation}
In particular $\mc R_\alpha^{-1}w\to\mc R_0^{-1}w$ in $\W^k_0$ as $\alpha\to0$, for
each fixed $w\in\W^{k+1}_0$.
\end{lemma}

\begin{proof}
Write
$\mc R_\alpha^{-1}-\mc R_0^{-1}
=\mc R_\alpha^{-1}\big(\mc R_0-\mc R_\alpha\big)\mc R_0^{-1}$ and estimate the three
factors from the right. First, $\mc R_0^{-1}=(\Id-P)^{-1}$ maps $\W^{k+1}_0$ into
itself with norm at most $C_R$, by Proposition~\ref{prop:gap} at level $k+1$, which is
legitimate since $k+1\le r-1$. Second, for $v\in\W^{k+1}_0$,
\[
\mc R_0-\mc R_\alpha=\big(\LL_{\alpha,\phi(\alpha)}-P\big)
+\alpha\,\mc K_{\alpha,\phi(\alpha)} ,
\]
where estimate \eqref{eq:est1} with $m=k$, $F=f$ and $\Delta=\alpha A_{\phi(\alpha),z}$
gives $\|(\LL_{\alpha,\phi(\alpha)}-P)v\|_{k}\le C|\alpha|\,\|v\|_{k+1}$ --- admissible
because $k+1\le r-1$, so that Lemma~\ref{lem:exact}(i) applies --- while
Proposition~\ref{prop:Kbounded} gives
$\|\alpha\mc K_{\alpha,\phi(\alpha)}v\|_{k}\le|\alpha|\,C_{\mc K}M\|v\|_{0}
\le|\alpha|\,C_{\mc K}M\|v\|_{k+1}$. Hence
$\|(\mc R_0-\mc R_\alpha)v\|_{k}\le C|\alpha|\,\|v\|_{k+1}$. Third,
$\mc R_\alpha^{-1}$ is bounded on $\W^k_0$ uniformly in $\alpha$ by
\eqref{eq:resolventbound}. Composing the three bounds gives \eqref{eq:resmixed}.
\end{proof}

\begin{remark}[Why a density argument is unavailable]\label{rem:nodensity}
Estimate \eqref{eq:resmixed} cannot be upgraded to operator-norm convergence on
$\W^k_0$ by density, because $\W^{k+1}$ is \emph{not} dense in
\[
\W^{k}=L^\infty\big([0,1],W^{k,1}(\T)\big) .
\]
Mollifying in $x$ does produce
convergence in $W^{k,1}(\T)$ along each fibre, but in the absence of a modulus of
continuity uniform in $z$ the convergence need not take place in the
$L^\infty$-in-$z$ norm. Each passage to the limit in this paper is therefore performed
on a single, explicitly named vector, after verifying that it has the derivative to
spare which \eqref{eq:resmixed} requires. This occurs twice: in the proof of
Theorem~\ref{thm:structural}, where the extra derivative is supplied free of charge by
the smoothing property \eqref{eq:Kbound} of $\mc K^V$, and in the proof of
Proposition~\ref{prop:formalgenuine}, where it is verified term by term in
Appendix~\ref{app:longproofs}.
\end{remark}

The family $\alpha\mapsto\mc R_{\alpha,\phi^\star(\alpha)}$ itself must be treated with
the same discipline. The next lemma records the two facts about it that we use: a
Lipschitz increment bound costing one derivative, and differentiability at $\alpha=0$
costing two. Both are statements about a \emph{fixed vector}; the corresponding
operator-norm statements are false, as Remark~\ref{rem:nodiff} shows.

\begin{lemma}[Increments and derivative of $\mc R_\alpha$ on regular vectors]
\label{lem:Rincr}
Let $1\le k\le r-3$ and let $\alpha\mapsto\phi(\alpha)\in\Delta^{k+2}$ satisfy
$\|\phi(\alpha)\|_{k+2}\le M$ for $|\alpha|<\alpha_2$, together with
$\|\phi(\alpha)-\phi(\beta)\|_{k+1}\le C_0|\alpha-\beta|$ and $\phi(0)=\rho$. Write
$R_\gamma:=\mc R_{\gamma,\phi(\gamma)}$. There is $C>0$, depending only on
$f,h,W,\alpha_0,k,M,C_0$, such that:
\begin{itemize}
\item[(i)] for all $|\alpha|,|\beta|<\alpha_2$ and every $v\in\W^{k+1}_0$,
\begin{equation}\label{eq:Rincr}
\big\|\big(R_\alpha-R_\beta\big)v\big\|_{k}\le C\,|\alpha-\beta|\,\|v\|_{k+1} ;
\end{equation}
\item[(ii)] setting
$\dot R\,v:=\big[P\big(A_{\rho,z}\,v_z\big)\big]'-\mc K_{0,\rho}v$, for every
$v\in\W^{k+2}_0$,
\begin{equation}\label{eq:Rdot}
\Big\|\frac{R_\alpha-R_0}{\alpha}\,v-\dot R\,v\Big\|_{k}
\le C\,|\alpha|\,\|v\|_{k+2} .
\end{equation}
\end{itemize}
\end{lemma}

\begin{proof}
Throughout write $A_{\gamma,z}:=A_{\phi(\gamma),z}$ and recall
$R_\gamma=\Id-\LL_{\gamma,\phi(\gamma)}-\gamma\mc K_{\gamma,\phi(\gamma)}$.

\emph{(i).} Exactly,
$F_{\alpha,\phi(\alpha),z}-F_{\beta,\phi(\beta),z}
=(\alpha-\beta)A_{\alpha,z}+\beta\big(A_{\alpha,z}-A_{\beta,z}\big)$, whose
$\mc C^{k+1}$ norm is at most $C|\alpha-\beta|$ by Lemma~\ref{lem:Abound} and the
assumed Lipschitz bound at level $0$. Hence \eqref{eq:est1} with $m=k$, legitimate
since $k+1\le r-1$, gives
$\|(\LL_{\alpha,\phi(\alpha)}-\LL_{\beta,\phi(\beta)})v\|_{k}
\le C|\alpha-\beta|\|v\|_{k+1}$. For the non-local part write
\[
\alpha\mc K_{\alpha,\phi(\alpha)}-\beta\mc K_{\beta,\phi(\beta)}
=(\alpha-\beta)\mc K_{\alpha,\phi(\alpha)}
+\beta\big(\mc K_{\alpha,\phi(\alpha)}-\mc K_{\beta,\phi(\beta)}\big) .
\]
The first term is bounded by $|\alpha-\beta|C_{\mc K}M\|v\|_{0}$ by
Proposition~\ref{prop:Kbounded}. For the second, \eqref{eq:KK} gives
\begin{align*}
\big(\mc K_{\alpha,\phi(\alpha)}v-\mc K_{\beta,\phi(\beta)}v\big)_z
&=-\Big[\big(L_{\alpha,\phi(\alpha),z}-L_{\beta,\phi(\beta),z}\big)
\big(A_{v,z}\phi(\alpha)_z\big)\Big]'\\
&\qquad-\Big[L_{\beta,\phi(\beta),z}
\Big(A_{v,z}\big(\phi(\alpha)-\phi(\beta)\big)_z\Big)\Big]' .
\end{align*}
Since $\|[g]'\|_{W^{k,1}}\le\|g\|_{W^{k+1,1}}$, the first bracket is estimated by
\eqref{eq:est1} with $m=k+1$ --- legitimate since $k+2\le r-1$ --- by
$C|\alpha-\beta|\,\|A_{v,z}\phi(\alpha)_z\|_{W^{k+2,1}}
\le C|\alpha-\beta|\,\|v\|_{0}$, using Lemma~\ref{lem:Abound}(i) at level $k+2$ and
$\|\phi(\alpha)\|_{k+2}\le M$; the second by Lemma~\ref{lem:LY} at level $k+1$ and
Lemma~\ref{lem:calculus} by
$C\|v\|_{0}\|\phi(\alpha)-\phi(\beta)\|_{k+1}\le CC_0|\alpha-\beta|\,\|v\|_{0}$.
Collecting and using $\|v\|_0\le\|v\|_{k+1}$ gives \eqref{eq:Rincr}.

\emph{(ii).} Take $v\in\W^{k+2}_0$. For the local part, \eqref{eq:exact1} with $F=f$,
$\Delta=\alpha A_{\alpha,z}$ and $u=v_z\in W^{k+2,1}(\T)$ gives exactly
\[
-\frac{\big(\LL_{\alpha,\phi(\alpha)}v-Pv\big)_z}{\alpha}
=\int_0^1\Big[L_{f+t\alpha A_{\alpha,z}}\big(A_{\alpha,z}v_z\big)\Big]'\dd t ,
\]
so that, subtracting the target $\big[P(A_{\rho,z}v_z)\big]'$ and using
$\|[g]'\|_{W^{k,1}}\le\|g\|_{W^{k+1,1}}$,
\[
\Big\|-\frac{\big(\LL_{\alpha,\phi(\alpha)}v-Pv\big)_z}{\alpha}
-\big[P(A_{\rho,z}v_z)\big]'\Big\|_{W^{k,1}}
\le\int_0^1\big\|L_{f+t\alpha A_{\alpha,z}}\big(A_{\alpha,z}v_z\big)
-P\big(A_{\rho,z}v_z\big)\big\|_{W^{k+1,1}}\dd t .
\]
Split the integrand as
$\big[L_{f+t\alpha A_{\alpha,z}}-P\big]\big(A_{\alpha,z}v_z\big)
+P\big((A_{\alpha,z}-A_{\rho,z})v_z\big)$. The first is at most
$C|\alpha|\,\|A_{\alpha,z}v_z\|_{W^{k+2,1}}\le C|\alpha|\,\|v_z\|_{W^{k+2,1}}$ by
\eqref{eq:est1} with $m=k+1$, admissible since $k+2\le r-1$; the second is at most
$C\|A_{\alpha,z}-A_{\rho,z}\|_{\mc C^{k+1}}\|v_z\|_{W^{k+1,1}}
\le C|\alpha|\,\|v_z\|_{W^{k+1,1}}$ by Lemma~\ref{lem:LY} at level $k+1$,
Lemma~\ref{lem:calculus} and Lemma~\ref{lem:Abound}(ii). For the non-local part, the
computation of (i) with $\beta=0$ and $\phi(0)=\rho$ gives
$\|\mc K_{\alpha,\phi(\alpha)}v-\mc K_{0,\rho}v\|_{k}\le C|\alpha|\,\|v\|_{0}$. Since
$(R_\alpha-R_0)/\alpha=-(\LL_{\alpha,\phi(\alpha)}-P)/\alpha
-\mc K_{\alpha,\phi(\alpha)}$, adding the two estimates gives \eqref{eq:Rdot}.
\end{proof}

\begin{remark}[Why \eqref{eq:Rdot} cannot be an operator-norm statement]\label{rem:nodiff}
It is tempting to read \eqref{eq:Rdot} as differentiability of $\alpha\mapsto R_\alpha$
at $0$ in $\mc L(\W^{k+1},\W^{k})$. That is \emph{false}, already inside the setting of
this paper. Take $f(x)=2x\bmod1$, $h\equiv1$ and $W\equiv1$. Then
$A_{\phi,z}\equiv1$ for every probability density, so $F_{\alpha,\phi,z}=f+\alpha$ and
$\phi^\star(\alpha)\equiv1$ for all $\alpha$; moreover $A_{v,z}=\int_0^1\!\!\int_\T v=0$
for $v\in\W^0_0$, so $\mc K$ vanishes on zero-average vectors. Writing
$T_\alpha v(x)=v(x-\alpha)$ one has $L_{f+\alpha}=T_\alpha P$, hence
\[
R_\alpha=\Id-T_\alpha P\quad\text{on }\W^k_0,\qquad \dot R\,v=(Pv)' .
\]
Test this on $u_n(x)=e^{4\pi i n x}$ with $\alpha_n=1/(2n)$. Since $Pu_n=e^{2\pi inx}$
and $T_{\alpha_n}e^{2\pi inx}=-e^{2\pi inx}$,
\[
\Big(\frac{R_{\alpha_n}-R_0}{\alpha_n}-\dot R\Big)u_n
=\big(4-2\pi i\big)\,n\,e^{2\pi i n x} ,
\]
whose $W^{k,1}$ norm is comparable to $n^{k+1}$, while
$\|u_n\|_{W^{k+1,1}}$ is comparable to $n^{k+1}$ as well. The quotient therefore stays
bounded away from $0$ as $n\to\infty$, so no operator-norm limit exists; taking real and
imaginary parts transfers the example to the real setting. The mechanism is transparent:
the quotient $\alpha^{-1}(T_\alpha-\Id)$ converges to $-\partial_x$ on each fixed
function, but not uniformly on bounded sets of $W^{k+1,1}$, because $L^1$ translation is
not uniformly continuous on bounded sets. This is why
Lemma~\ref{lem:Rincr}(ii) is stated for a fixed $v$, with the second derivative of $v$
on the right-hand side, and why the proof of
Proposition~\ref{prop:formalgenuine} applies it only to the single vector
$u_0=\mc R_{0}^{-1}N(0)$, which is shown there to lie in $\W^{k+2}_0$.
\end{remark}

\begin{remark}[Structure of the response equation]\label{rem:structure}
Expanding \eqref{eq:neumann} makes the network structure explicit. The operator
$(\Id-\LL)^{-1}=\sum_{n\ge0}\LL^n$ is diagonal in $z$: it propagates information
\emph{within} a fibre, along the local dynamics. The operator $\mc K$ is an integral
operator in $z$ with kernel $W$: it propagates information \emph{across} the network,
one edge at a time. The Neumann series therefore alternates temporal and spatial
propagation, and its $m$-th term corresponds to paths traversing $m$ edges of the
network. This is the mechanism behind the tree expansion of Theorem~\ref{thm:tree}.
\end{remark}

\subsection{Uniform contraction of the STO}

The following proposition replaces the mean value inequality, which is not available
here because $D\STO_\alpha$ loses one derivative (Remark~\ref{rem:loss}).

\begin{proposition}[Uniform exponential contraction]\label{prop:contraction}
Let $1\le k\le r-2$ and $M\ge M_{k+1}$. There are $\alpha_1\in(0,\alpha_0]$,
$D\ge1$ and $\eta\in(\theta,1)$ such that for all $|\alpha|<\alpha_1$, all
$\psi_1,\psi_2\in\Delta^{k+1}$ with $\|\psi_i\|_{k+1}\le M$ and all $n\ge0$,
\begin{equation}\label{eq:contraction}
\big\|\STO_\alpha^n\psi_1-\STO^n_\alpha\psi_2\big\|_{k}\le D\,\eta^n\,
\|\psi_1-\psi_2\|_{k} .
\end{equation}
\end{proposition}

\begin{proof}
Write $\psi_i^{(m)}:=\STO^m_\alpha\psi_i$ and abbreviate
\[
A_m:=\LL_{\alpha,\psi_1^{(m)}},\qquad B_m:=\LL_{\alpha,\psi_2^{(m)}},\qquad
A_{[m,n)}:=A_{n-1}A_{n-2}\cdots A_{m},\qquad
B_{[m,n)}:=B_{n-1}\cdots B_{m},
\]
with the convention that an empty product is the identity. Then
$\STO^n_\alpha\psi_i$ equals $A_{[0,n)}\psi_1$ for $i=1$ and $B_{[0,n)}\psi_2$ for
$i=2$. By Proposition~\ref{prop:apriori} we may assume
$\|\psi_i^{(m)}\|_{k+1}\le M$ for all $m$, enlarging $M$ if necessary. Write
\[
\STO^n_\alpha\psi_1-\STO^n_\alpha\psi_2
=A_{[0,n)}\big(\psi_1-\psi_2\big)
+\big(A_{[0,n)}-B_{[0,n)}\big)\psi_2 .
\]
Since $\psi_1-\psi_2\in\W^k_0$ and, at each fibre $z$, $A_{[0,n)}$ acts as a
sequential composition of transfer operators of maps in $\mc U_z(\alpha)$, at the
single coupling strength $\alpha$ under consideration, the first term is
bounded by
$C_\star\theta^n\|\psi_1-\psi_2\|_k$ by Proposition~\ref{prop:gap}. For the second we
use the telescopic identity
\begin{equation}\label{eq:telescope}
A_{[0,n)}-B_{[0,n)}=\sum_{m=0}^{n-1}
A_{[m+1,n)}\,\big(A_m-B_m\big)\,B_{[0,m)} ,
\end{equation}
which follows by induction on $n$: for $n=1$ both sides are $A_0-B_0$, and
\[
A_{[0,n+1)}-B_{[0,n+1)}=A_n\big(A_{[0,n)}-B_{[0,n)}\big)+(A_n-B_n)B_{[0,n)} ,
\]
which is the $n+1$ case once \eqref{eq:telescope} is inserted in the first summand.
Note that $B_{[0,m)}\psi_2=\psi_2^{(m)}$, that $A_{[m+1,n)}$ is a sequential
composition of $n-1-m$ operators, and that
$(A_m-B_m)u\in\W^k_0$ for every $u$, because both $A_m$ and $B_m$ preserve the
fibrewise integral; hence Proposition~\ref{prop:gap} applies to the left factor of each
summand of \eqref{eq:telescope}. By \eqref{eq:est1} with
$\Delta=\alpha\big(A_{\psi_1^{(m)},z}-A_{\psi_2^{(m)},z}\big)$ and
Lemma~\ref{lem:Abound}(ii),
\[
\Big\|\big(\LL_{\alpha,\psi_1^{(m)}}-\LL_{\alpha,\psi_2^{(m)}}\big)u\Big\|_{k}
\le C|\alpha|\,\big\|\psi_1^{(m)}-\psi_2^{(m)}\big\|_{0}\,\|u\|_{k+1} ,
\]
and the argument $u=\psi_2^{(m)}$ satisfies $\|u\|_{k+1}\le M$. Setting
$e_m:=\|\STO^m_\alpha\psi_1-\STO^m_\alpha\psi_2\|_{k}$ and using
$\|\cdot\|_0\le\|\cdot\|_k$ we obtain
\begin{equation}\label{eq:renewal}
e_n\le C_\star\theta^n e_0+C'|\alpha|\sum_{m=0}^{n-1}\theta^{\,n-1-m}e_m .
\end{equation}
Fix $\eta\in(\theta,1)$ and $D\ge\max(1,2C_\star)$, and choose $\alpha_1$ so small that
$C'\alpha_1/(\eta-\theta)\le\tfrac12$. We claim $e_n\le D\eta^ne_0$ for all $n$,
by strong induction. It is true for $n=0$. If it holds for $m<n$, then by
\eqref{eq:renewal},
\[
e_n\le C_\star\theta^ne_0+C'|\alpha|D e_0\sum_{m=0}^{n-1}\theta^{\,n-1-m}\eta^m
\le \tfrac{D}{2}\eta^ne_0+\frac{C'|\alpha|D}{\eta-\theta}\,\eta^{n}e_0
\le D\eta^ne_0 ,
\]
where we used $\sum_{m=0}^{n-1}\theta^{n-1-m}\eta^m\le\eta^{n}/(\eta-\theta)$. This
proves \eqref{eq:contraction}.
\end{proof}

\section{Proofs of Theorems A, B and C}\label{sec:proofs}

\subsection{Lipschitz statistical stability}

\begin{proof}[Proof of Theorem~\ref{thm:stability}]
Let $1\le k\le r-2$ and write $\phi^\alpha:=\phi^\star(\alpha)$. By
Proposition~\ref{prop:apriori}, $\|\phi^\gamma\|_{k+1}\le M_{k+1}$ for all
$|\gamma|<\alpha_0$. Fix $\alpha,\beta$ and set $u_m:=\STO^m_\alpha\phi^\beta$, so that
$u_0=\phi^\beta$. For every $n\ge1$,
\begin{equation}\label{eq:stabA}
\big\|\phi^\alpha-\phi^\beta\big\|_{k}
\le\big\|\STO^n_\alpha\phi^\alpha-\STO^n_\alpha\phi^\beta\big\|_{k}
+\big\|u_n-u_0\big\|_{k} ,
\end{equation}
because $\STO^n_\alpha\phi^\alpha=\phi^\alpha$ and $u_n=\STO_\alpha^n\phi^\beta$.
By Proposition~\ref{prop:contraction} the first term is at most
$D\eta^n\|\phi^\alpha-\phi^\beta\|_{k}$. For the second, telescope and apply
Proposition~\ref{prop:contraction} once more, noting
$u_{m+1}-u_m=\STO^m_\alpha(\STO_\alpha\phi^\beta)-\STO^m_\alpha(\phi^\beta)$:
\[
\big\|u_n-u_0\big\|_{k}\le\sum_{m=0}^{n-1}\|u_{m+1}-u_m\|_{k}
\le\frac{D}{1-\eta}\ \big\|\STO_\alpha\phi^\beta-\phi^\beta\big\|_{k} .
\]
Since $\phi^\beta=\STO_\beta\phi^\beta$, estimate \eqref{eq:alphaLip} of
Lemma~\ref{lem:Kphi}(i) gives
\[
\big\|\STO_\alpha\phi^\beta-\phi^\beta\big\|_{k}
=\big\|\STO_\alpha\phi^\beta-\STO_\beta\phi^\beta\big\|_{k}
\le C|\alpha-\beta|\,\|\phi^\beta\|_{k+1}\le CM_{k+1}|\alpha-\beta| .
\]
Here Lemma~\ref{lem:Kphi}(i) applies because $k\le r-2$, which is exactly the range of
the present theorem; this is the point at which the one-derivative estimate
\eqref{eq:est1} is needed rather than the two-derivative estimate \eqref{eq:est2}.
Inserting both bounds in \eqref{eq:stabA} and choosing
$\bar n$ with $D\eta^{\bar n}\le\tfrac12$,
\[
\big\|\phi^\alpha-\phi^\beta\big\|_{k}\le\tfrac12
\big\|\phi^\alpha-\phi^\beta\big\|_{k}+\frac{DCM_{k+1}}{1-\eta}|\alpha-\beta| ,
\]
and the left-hand side is finite by Proposition~\ref{prop:apriori}, so it may be
absorbed. This gives the claim with $C=2DCM_{k+1}(1-\eta)^{-1}$.
\end{proof}

\subsection{Linear response}

\begin{proof}[Proof of Theorem~\ref{thm:LR}]
Let $1\le k\le r-3$, fix $\hat\alpha$ with $|\hat\alpha|<\alpha_2$ and write
$\hat\phi=\phi^\star(\hat\alpha)$, $\phi^\alpha=\phi^\star(\alpha)$,
\[
u_\alpha:=\frac{\phi^\alpha-\hat\phi}{\alpha-\hat\alpha}\in\W^{k+1}_0 .
\]
Part (i) is Proposition~\ref{prop:resolvent}, applicable since $k\le r-2$ and
$\|\hat\phi\|_{k+1}\le M_{k+1}$.

\emph{A priori bound.} By Theorem~\ref{thm:stability} applied at level $k+1$
(legitimate because $k+1\le r-2$), $\|u_\alpha\|_{k+1}\le C$ uniformly in $\alpha$.

\emph{Decomposition.} Since $\phi^\alpha$ and $\hat\phi$ are fixed points of
$\STO_\alpha$ and $\STO_{\hat\alpha}$ respectively,
\[
\phi^\alpha-\hat\phi
=\underbrace{\big[\STO_\alpha\phi^\alpha-\STO_\alpha\hat\phi\big]}_{(\mathrm I)}
+\underbrace{\big[\STO_\alpha\hat\phi-\STO_{\hat\alpha}\hat\phi\big]}_{(\mathrm{II})} .
\]
By Theorem~\ref{thm:diff} at $\phi=\hat\phi$ with $v=\phi^\alpha-\hat\phi$ (both lie in
$\Delta^{k+2}$ with norms $\le M_{k+2}$, and $k+2\le r-1$),
\[
(\mathrm I)=\big(\LL_{\alpha,\hat\phi}+\alpha\mc K_{\alpha,\hat\phi}\big)
\big(\phi^\alpha-\hat\phi\big)+R_\alpha,\qquad
\|R_\alpha\|_{k}\le C_M|\alpha|\,\|\phi^\alpha-\hat\phi\|_{k+1}^2 .
\]
By Theorem~\ref{thm:stability} at level $k+1$, $\|R_\alpha\|_k\le
C'|\alpha-\hat\alpha|^2$. Dividing by $\alpha-\hat\alpha$ and rearranging,
\begin{equation}\label{eq:LRmain}
\Big(\Id-\LL_{\alpha,\hat\phi}-\alpha\mc K_{\alpha,\hat\phi}\Big)u_\alpha
=\frac{\STO_\alpha\hat\phi-\STO_{\hat\alpha}\hat\phi}{\alpha-\hat\alpha}
+\frac{R_\alpha}{\alpha-\hat\alpha} ,
\end{equation}
with $\|R_\alpha\|_k/|\alpha-\hat\alpha|\le C'|\alpha-\hat\alpha|\to0$.

\emph{Replacing $\alpha$ by $\hat\alpha$ in the operator.} Since
$F_{\alpha,\hat\phi,z}-F_{\hat\alpha,\hat\phi,z}=(\alpha-\hat\alpha)A_{\hat\phi,z}$
exactly, estimate \eqref{eq:est1} gives
\[
\big\|\big(\LL_{\alpha,\hat\phi}-\LL_{\hat\alpha,\hat\phi}\big)u_\alpha\big\|_{k}
\le C|\alpha-\hat\alpha|\,\|u_\alpha\|_{k+1}\le C''|\alpha-\hat\alpha| ,
\]
and similarly, using \eqref{eq:est1} inside the definition \eqref{eq:KK} of $\mc K$
together with Proposition~\ref{prop:Kbounded},
\[
\big\|\big(\alpha\mc K_{\alpha,\hat\phi}-\hat\alpha\mc K_{\hat\alpha,\hat\phi}\big)
u_\alpha\big\|_{k}\le C|\alpha-\hat\alpha|\,\|u_\alpha\|_{k+1}
\le C''|\alpha-\hat\alpha| .
\]
Here the uniform bound $\|u_\alpha\|_{k+1}\le C$ is what makes the loss of one
derivative harmless.

\emph{Right-hand side.} By Lemma~\ref{lem:Kphi}(ii) (applicable since
$\hat\phi\in\Delta^{k+2}$),
\[
\frac{\STO_\alpha\hat\phi-\STO_{\hat\alpha}\hat\phi}{\alpha-\hat\alpha}
\ \xrightarrow[\alpha\to\hat\alpha]{\W^{k}}\ \mc K_{\hat\alpha,\hat\phi}\hat\phi .
\]

\emph{Conclusion.} Collecting, \eqref{eq:LRmain} becomes
$\mc R_{\hat\alpha,\hat\phi}u_\alpha=\mc K_{\hat\alpha,\hat\phi}\hat\phi+\eps_\alpha$
with $\|\eps_\alpha\|_{k}\to0$, so by Proposition~\ref{prop:resolvent}
\[
u_\alpha=\mc R_{\hat\alpha,\hat\phi}^{-1}\mc K_{\hat\alpha,\hat\phi}\hat\phi
+\mc R^{-1}_{\hat\alpha,\hat\phi}\eps_\alpha
\ \xrightarrow[\alpha\to\hat\alpha]{\ \W^k\ }\
\mc R_{\hat\alpha,\hat\phi}^{-1}\mc K_{\hat\alpha,\hat\phi}\hat\phi ,
\]
which is \eqref{eq:LRgeneral}.

\emph{Part (iii).} It would be circular to invoke continuity of inversion in
\eqref{eq:neumann}: the operators $\mc R_{\gamma,\phi^\star(\gamma)}$ are continuous in
$\gamma$ only in the mixed norm $\mc L(\W^{k+1},\W^k)$, by \eqref{eq:Rincr}, and
inversion is not continuous for that topology. We argue on increments instead. Write
$R_\gamma:=\mc R_{\gamma,\phi^\star(\gamma)}$,
$N_\gamma:=\mc K_{\gamma,\phi^\star(\gamma)}\phi^\star(\gamma)$ and
$D(\gamma):=\partial_\alpha\phi^\star(\gamma)$, so that part (ii) reads
$R_\gamma D(\gamma)=N_\gamma$. Subtracting the identities at $\alpha$ and $\beta$ gives
the exact relation
\begin{equation}\label{eq:Dincrement2}
R_\alpha\big(D(\alpha)-D(\beta)\big)
=\big(N_\alpha-N_\beta\big)+\big(R_\beta-R_\alpha\big)D(\beta) ,
\end{equation}
and we estimate the two terms on the right.

\emph{(a) A uniform bound one level up: $\sup_{|\gamma|<\alpha_2}\|D(\gamma)\|_{k+1}
<\infty$.} By Proposition~\ref{prop:Kbounded} at level $k+1$, legitimate since
$k+1\le r-2$ and $\phi^\star(\gamma)\in\Delta^{k+2}$, the vector $N_\gamma$ lies in
$\W^{k+1}_0$ with $\|N_\gamma\|_{k+1}\le C_{\mc K}M_{k+2}$. By
Proposition~\ref{prop:resolvent} at level $k+1$, also legitimate since $k+1\le r-2$,
the operator $R_\gamma$ is invertible on $\W^{k+1}_0$ with a bound uniform in $\gamma$;
since it is injective on the larger space $\W^{k}_0$, the two inverses agree on
$\W^{k+1}_0$, so $D(\gamma)=R_\gamma^{-1}N_\gamma\in\W^{k+1}_0$ with
$\|D(\gamma)\|_{k+1}\le C$ uniformly. Note that this uses the smoothing property of
$\mc K$ and not a limit of difference quotients: bounded sets of $\W^{k+1}$ are not
closed under convergence in $\W^k$, so the \emph{a priori} bound
$\|u_\alpha\|_{k+1}\le C$ would not by itself place the limit $D(\gamma)$ in
$\W^{k+1}$.

\emph{(b) $\|N_\alpha-N_\beta\|_{k}\le C|\alpha-\beta|$.} Writing
$g_\gamma:=A_{\phi^\star(\gamma),z}\phi^\star_z(\gamma)$, \eqref{eq:KK} gives
\[
\big(N_\alpha-N_\beta\big)_z
=-\Big[\big(L_{\alpha,\phi^\star(\alpha),z}-L_{\beta,\phi^\star(\beta),z}\big)
g_\alpha\Big]'
-\Big[L_{\beta,\phi^\star(\beta),z}\big(g_\alpha-g_\beta\big)\Big]' .
\]
Using $\|[g]'\|_{W^{k,1}}\le\|g\|_{W^{k+1,1}}$, the first term is at most
$C|\alpha-\beta|\,\|g_\alpha\|_{W^{k+2,1}}\le C|\alpha-\beta|$ by \eqref{eq:est1} with
$m=k+1$, admissible since $k+2\le r-1$, together with
Lemma~\ref{lem:Abound} and $\phi^\star(\alpha)\in\Delta^{k+2}$. For the second,
Lemma~\ref{lem:LY} at level $k+1$ and bilinearity give
$\|g_\alpha-g_\beta\|_{W^{k+1,1}}
\le C\|\phi^\star(\alpha)-\phi^\star(\beta)\|_{0}+C\|\phi^\star(\alpha)
-\phi^\star(\beta)\|_{k+1}\le C|\alpha-\beta|$, by Lemma~\ref{lem:Abound} and
Theorem~\ref{thm:stability} at level $k+1$.

\emph{(c) Conclusion.} By Lemma~\ref{lem:Rincr}(i) applied to the vector $D(\beta)$,
which lies in $\W^{k+1}_0$ by (a),
$\|(R_\beta-R_\alpha)D(\beta)\|_{k}\le C|\alpha-\beta|\,\|D(\beta)\|_{k+1}
\le C'|\alpha-\beta|$. Applying $R_\alpha^{-1}$ to \eqref{eq:Dincrement2} and using
\eqref{eq:resolventbound},
\[
\big\|D(\alpha)-D(\beta)\big\|_{k}
\le C_R'\Big(\|N_\alpha-N_\beta\|_{k}
+\big\|(R_\beta-R_\alpha)D(\beta)\big\|_{k}\Big)
\le C''|\alpha-\beta| ,
\]
which is part (iii). No continuity of $\gamma\mapsto R_\gamma^{-1}$ in
$\mc L(\W^k_0)$ has been used: the only inverse that appears is $R_\alpha^{-1}$ at the
single value $\alpha$, through its uniform bound.
\end{proof}

\subsection{The zero-coupling limit}

\begin{proof}[Proof of Theorem~\ref{thm:zerocoupling}]
At $\alpha=0$ we have $F_{0,\phi,z}=f$ for every $\phi$ and every $z$, so
$\LL_{0,\phi}=P$ acting fibrewise and $(\STO_0\phi)_z=P\phi_z$. The unique fixed point
in $\Delta^{k}$ is therefore $\phi^\star_z(0)=\rho$ for a.e.\ $z$, and
$\mc R_{0,\phi^\star(0)}=\Id-P$, invertible on $\W^k_0$ by
Proposition~\ref{prop:gap}.

Theorem~\ref{thm:LR} at $\hat\alpha=0$ gives
$\partial_\alpha\phi^\star|_{0}=(\Id-P)^{-1}\mc K_{0,\rho}\rho$. By \eqref{eq:KK} with
$L_{0,\cdot,z}=P$ and $\phi_z=\rho$,
\[
\big(\mc K_{0,\rho}\,\rho\big)_z=-\Big[P\big(A_{\rho,z}\,\rho\big)\Big]' .
\]
Since $\phi^\star_{z'}(0)=\rho$ for \emph{every} $z'$, the function $H_{\phi_{z'}}=H_\rho$
does not depend on $z'$ and \eqref{eq:A} gives
\[
A_{\rho,z}=\Big(\int_0^1 W(z,z')\,\dd z'\Big)H_{\rho}=\barW(z)\,H_\rho .
\]
Hence $\big(\mc K_{0,\rho}\rho\big)_z=-\barW(z)\big[P(H_\rho\rho)\big]'$ and, the
resolvent being $z$-independent and linear,
\[
\partial_\alpha\phi^\star_z(\alpha)\big|_{0}
=-\barW(z)\,(\Id-P)^{-1}\big[P(H_\rho\rho)\big]'=\barW(z)\,\Up .
\]
Convergence in $\esssup_z\|\cdot\|_{W^{k,1}}$ is convergence in $\W^k$, which is what
Theorem~\ref{thm:LR} asserts. Finally $\int_\T[P(H_\rho\rho)]'=0$, so $\Up$ is well
defined and has zero average; and $H_\rho\rho\in W^{r-1,1}$ by
Remark~\ref{rem:BVW11}(ii) and $\rho\in\mc C^{r-1}$, so $\Up\in W^{r-2,1}$.
\end{proof}

\begin{proof}[Proof of Corollary~\ref{cor:GK}]
The first equality follows from Theorem~\ref{thm:zerocoupling} and
$\big|\int_\T c\,(\phi^\star_z(\alpha)-\rho-\alpha\barW(z)\Up)\big|
\le\|c\|_{\infty}\,o(\alpha)$, uniformly in $z$. For the second, expand
$(\Id-P)^{-1}=\sum_{n\ge0}P^n$ on zero-average elements of $W^{1,1}$ and use
$\int c\,P^n u=\int (c\circ f^n)\,u$ with $u=[P(H_\rho\rho)]'$:
\[
\int_\T c\,\Up=-\sum_{n\ge0}\int_\T c\,P^n\big[P(H_\rho\rho)\big]'
=-\sum_{n\ge0}\int_\T (c\circ f^n)\big[P(H_\rho\rho)\big]' .
\]
Convergence is exponential because
$\|P^nu\|_{L^1}\le\|P^nu\|_{W^{1,1}}\le C_\star\theta^n\|u\|_{W^{1,1}}$ by
Proposition~\ref{prop:gap} with $k=1$. The global statement follows by integrating in
$z$ and Fubini, since $\int_0^1\barW(z)\dd z=\int\int W$.
\end{proof}

\begin{proof}[Proof of Proposition~\ref{prop:zregularity}]
By Theorem~\ref{thm:zerocoupling},
$\dot\phi_z-\dot\phi_{\bar z}=\big(\barW(z)-\barW(\bar z)\big)\Up$, so
$|\dot\phi_z-\dot\phi_{\bar z}|_{\BV^1}\le\|\Up\|_{\BV^1}|\barW(z)-\barW(\bar z)|$.
Since
$|\barW(z)-\barW(\bar z)|\le\|W(z,\cdot)-W(\bar z,\cdot)\|_{L^1}$, taking the essential
supremum over $z,\bar z\in B(\omega,\varrho)$ gives
$\mathtt{osc}_{\BV^1}(\dot\phi,\omega,\varrho)\le\|\Up\|_{\BV^1}
\mathtt{osc}_{L^1}(W,\omega,\varrho)$; integrating in $\omega$, dividing by
$\varrho^p$ and taking the supremum over $\varrho>0$ gives the claim.
\end{proof}

\section{Higher order response and the tree expansion}\label{sec:higher}

\subsection{The order-by-order recursion}

\begin{lemma}[Expansion of $L_{f+\Delta}$]\label{lem:expansionL}
Let $f+\Delta\in\mc U$, let $p_{\max}\ge0$ and let $u\in W^{m+p_{\max}+1,1}(\T)$ with
$m+p_{\max}+1\le r-1$. Then
\begin{equation}\label{eq:expansionL}
L_{f+\Delta}u=\sum_{p=0}^{p_{\max}}\frac{(-1)^p}{p!}
\Big[P\big(\Delta^{p}\,u\big)\Big]^{(p)}+\mc E ,\qquad
\|\mc E\|_{W^{m,1}}\le C\,\|\Delta\|^{p_{\max}+1}_{\mc C^{m+p_{\max}+1}}
\,\|u\|_{W^{m+p_{\max}+1,1}} .
\end{equation}
\end{lemma}

\begin{proof}
See Appendix~\ref{app:calculus}; it is the Taylor expansion of $g\circ(f+\Delta)$ with
integral remainder, transported through the duality \eqref{eq:duality} and
$p$ integrations by parts, exactly as in Lemma~\ref{lem:exact}.
\end{proof}

\begin{proposition}[Any genuine expansion obeys the recursion]\label{prop:recursion}
Suppose that, for some $m\le r-2$, $\phi^\star_z(\alpha)=\rho+\sum_{j=1}^m
\alpha^j\psi_{j,z}+o(\alpha^m)$ in $\W^{k}$ with $k+m+1\le r-1$ and
$\psi_j\in\W^k_0$. Then $\psi_j=\phi_j$ for $1\le j\le m$, the formal
susceptibilities of \eqref{eq:formaldef}.
\end{proposition}

\begin{proof}
The natural argument --- substitute the two expansions into \eqref{eq:expansionL} with
$p_{\max}=m$ and compare coefficients --- is not legitimate as it stands, and we spell
out why before giving a proof that is. The hypothesis controls the remainder of
$\phi^\star(\alpha)$ only in $\W^{k}$, whereas the $p$-th term of
\eqref{eq:expansionL} applies $p$ derivatives to an expression built from $u$: a
remainder that is $o(\alpha^m)$ in $\W^{k}$ contributes
$[P(\Delta^p\cdot o(\alpha^m))]^{(p)}$, which is controlled only in $\W^{k-p}$. To
compare coefficients in $\W^{k}$ one would need the remainder to be $o(\alpha^m)$ in
$\W^{k+m}$, which is neither assumed nor implied. Note also that the mean field and the
density must be expanded with their remainders, not as infinite series whose
convergence is not known. We therefore test the fixed-point identity against a smooth
function, so that \emph{every} derivative falls on the test function and no derivative
is ever taken of a remainder.

\emph{Step 1: expansion of the mean field.} Put $\psi_0:=\rho$. By
Lemma~\ref{lem:Abound}(i) with $j=r$ the map $v\mapsto A_{v,z}$ is linear with
$\sup_z\|A_{v,z}\|_{\mc C^{r}}\le\|h\|_{\mc C^{r}}\|W\|_{\infty,1}\|v\|_{0}$, and
$\|\cdot\|_{0}\le\|\cdot\|_{k}$; the hypothesis therefore gives
\begin{equation}\label{eq:Dexpansion}
\Delta_{\alpha,z}:=\alpha A_{\phi^\star(\alpha),z}
=\sum_{q=0}^{m}\alpha^{q+1}A_{\psi_q,z}+o\big(\alpha^{m+1}\big)
\qquad\text{in }L^\infty\big([0,1],\mc C^{r}(\T)\big) ,
\end{equation}
and in particular $\sup_z\|\Delta_{\alpha,z}\|_{\mc C^{r}}\le C|\alpha|$.

\emph{Step 2: a scalar identity.} Fix $g\in\mc C^{\infty}(\T)$ and $z\in A_W$. Since
$L_F$ is the density of the push-forward under $F$,
\begin{equation}\label{eq:dualityL}
\int_\T g\,\big(L_{F}u\big)\dd x=\int_\T (g\circ F)\,u\dd x ,
\qquad u\in L^1(\T) ,
\end{equation}
so the fixed-point equation
$\phi^\star_z(\alpha)=L_{f+\Delta_{\alpha,z}}\phi^\star_z(\alpha)$ becomes the scalar
identity
\begin{equation}\label{eq:scalarfix}
\int_\T g\,\phi^\star_z(\alpha)\dd x
=\int_\T g\big(f+\Delta_{\alpha,z}\big)\,\phi^\star_z(\alpha)\dd x .
\end{equation}
Taylor's theorem applied to $g$ gives
\[
g\big(f+\Delta_{\alpha,z}\big)=\sum_{p=0}^{m}\frac{1}{p!}
\big(g^{(p)}\circ f\big)\,\Delta^{p}_{\alpha,z}+\mc r_{\alpha,z},
\qquad
\|\mc r_{\alpha,z}\|_{\infty}
\le\frac{\|g^{(m+1)}\|_{\infty}}{(m+1)!}\,
\|\Delta_{\alpha,z}\|^{m+1}_{\infty}=O\big(|\alpha|^{m+1}\big) ,
\]
uniformly in $z$. Every derivative appearing here is a derivative of $g$. Since
$\|\phi^\star_z(\alpha)\|_{L^1}=1$, inserting this into \eqref{eq:scalarfix} and then
substituting \eqref{eq:Dexpansion} and the assumed expansion of
$\phi^\star_z(\alpha)$ --- each product pairing a uniformly bounded factor with an
$L^1$ factor --- expands both sides of \eqref{eq:scalarfix} in powers of $\alpha$ to
order $m$ with remainder $o(\alpha^m)$, uniformly in $z$.

\emph{Step 3: comparison of coefficients.} Two scalar functions of $\alpha$ with
expansions to order $m$ and $o(\alpha^m)$ remainders have the same coefficients, so we
may equate them. On the left the coefficient of $\alpha^j$ is
$\int_\T g\,\psi_{j,z}$. On the right, the term $p=0$ contributes
$\int_\T(g\circ f)\psi_{j,z}=\int_\T g\,P\psi_{j,z}$ by \eqref{eq:dualityL} with
$F=f$, while by \eqref{eq:Dexpansion} the coefficient of $\alpha^{j}$ in
$\Delta^{p}_{\alpha,z}\phi^\star_z(\alpha)$ is
$\sum\prod_{a=1}^{p}A_{\psi_{j_a},z}\,\psi_{\ell,z}$, the sum being over
$\ell,j_1,\dots,j_p\ge0$ with $p+\ell+j_1+\cdots+j_p=j$. For each such term
\[
\int_\T\big(g^{(p)}\circ f\big)\,w\dd x=\int_\T g^{(p)}\,Pw\dd x
=(-1)^{p}\int_\T g\,\big[Pw\big]^{(p)}\dd x ,
\qquad w:=\prod_{a=1}^{p}A_{\psi_{j_a},z}\,\psi_{\ell,z} ,
\]
the first equality by \eqref{eq:dualityL} and the second by integration by parts on
$\T$, understood in the sense of distributions. Equating coefficients therefore gives,
for every $g\in\mc C^{\infty}(\T)$ and every $1\le j\le m$,
\[
\Big\langle(\Id-P)\psi_{j,z}\ -\!\!\!
\sum_{\substack{p\ge1,\ \ell\ge0,\ j_1,\dots,j_p\ge0\\ p+\ell+j_1+\cdots+j_p=j}}\!\!\!
\frac{(-1)^{p}}{p!}\Big[P\Big(\prod_{a=1}^{p}A_{\psi_{j_a},z}\,\psi_{\ell,z}\Big)
\Big]^{(p)},\ g\Big\rangle=0 ,
\]
that is, the bracket vanishes as a distribution on $\T$.

\emph{Step 4: from distributions to $\W^k$, and induction.} Both terms of the bracket
are in fact functions in $W^{k,1}(\T)$: the first because $\psi_j\in\W^k_0$ and $P$
preserves $W^{k,1}(\T)$, the second by Proposition~\ref{prop:formal}, whose budget is
available precisely because $k+m+1\le r-1$. A distribution that vanishes and is given
by a $W^{k,1}$ function vanishes in $W^{k,1}$, and the above holds for a.e.\ $z$, with
all the estimates of Steps 1--3 uniform in $z$; hence
\[
(\Id-P)\psi_{j}
=\sum_{\substack{p\ge1,\ \ell\ge0,\ j_1,\dots,j_p\ge0\\ p+\ell+j_1+\cdots+j_p=j}}
\frac{(-1)^{p}}{p!}\Big[P\Big(\prod_{a=1}^{p}A_{\psi_{j_a}}\,\psi_{\ell}\Big)
\Big]^{(p)}\qquad\text{in }\W^{k} .
\]
The right-hand side involves only $\psi_i$ with $i\le j-1$, since $p\ge1$, and is a sum
of total derivatives, hence of zero fibrewise average; as $\psi_j\in\W^k_0$ and
$\Id-P$ is injective on $\W^k_0$ by Proposition~\ref{prop:gap}, this determines
$\psi_j$ from $\psi_0,\dots,\psi_{j-1}$ by exactly the recursion
\eqref{eq:formaldef}. Since $\psi_0=\rho=\phi_0$, strong induction on $j$ gives
$\psi_j=\phi_j$ for $1\le j\le m$.
\end{proof}

\subsection{Second order}

\begin{proof}[Proof of Proposition~\ref{prop:formalgenuine}]
The case $j=1$ is Theorem~\ref{thm:zerocoupling} combined with
Proposition~\ref{prop:recursion} at $m=1$. For $j=2$ one differentiates the right-hand
side of \eqref{eq:LRgeneral} at $\alpha=0$, which gives the existence of the
second-order term; the computation is carried out in Appendix~\ref{app:longproofs},
and the budget $k\le r-5$ comes from applying Theorem~\ref{thm:LR} at level $k+2$.
That the resulting coefficient is the formal susceptibility $\phi_2$ of
\eqref{eq:formaldef} is Proposition~\ref{prop:recursion} at $m=2$.
\end{proof}

\begin{proof}[Proof of Theorem~\ref{thm:secondorder}]
Apply \eqref{eq:formaldef} with $j=2$. The admissible tuples are
$(p,\ell,j_1)=(1,1,0)$, $(p,\ell,j_1)=(1,0,1)$, and $(p,\ell,j_1,j_2)=(2,0,0,0)$,
yielding
\[
(\Id-P)\phi_{2,z}=-\big[P(A_{\phi_0,z}\phi_{1,z})\big]'
-\big[P(A_{\phi_1,z}\rho)\big]'+\tfrac12\big[P(A_{\phi_0,z}^2\rho)\big]'' .
\]
Now $\phi_{0,z}=\rho$ gives $A_{\phi_0,z}=\barW(z)H_\rho$, while
$\phi_{1,z}=\barW(z)\Up$ gives, by \eqref{eq:A} and linearity,
\[
A_{\phi_1,z}=\int_0^1 W(z,z')\,H_{\barW(z')\Up}\,\dd z'
=\Big(\int_0^1 W(z,z')\,\barW(z')\,\dd z'\Big)H_\Up
=t_z(P_2,W)\,H_\Up .
\]
Therefore
\[
(\Id-P)\phi_{2,z}=\barW(z)^2\Big(\tfrac12\big[P(H_\rho^2\rho)\big]''
-\big[P(H_\rho\Up)\big]'\Big)-t_z(P_2,W)\big[P(H_\Up\rho)\big]' ,
\]
and applying $(\Id-P)^{-1}$, legitimate since both brackets are total derivatives and
hence of zero average, gives \eqref{eq:secondorder} with
\eqref{eq:UpLambda}--\eqref{eq:UpP}. Finally $\barW(z)^2=t_z(\Lambda_2,W)$ by
\eqref{eq:homdensity}.
\end{proof}

\begin{remark}
The second-order coefficient separates two distinct structural quantities: the square
of the degree of $z$, carried by the cherry, and the mean degree of the neighbours of
$z$, carried by the path. The latter is the graphon version of the quantity underlying
the friendship paradox. Two graphons can share a degree function and still be
distinguished dynamically at second order, as Example~\ref{ex:rank1} shows.
\end{remark}

\subsection{The tree expansion}

We isolate the two graph-theoretic identities used in the induction. For rooted trees
$T,T_1,\dots,T_p$ let $T^+$ denote the tree obtained by attaching a new root to the
root of $T$ by one edge, and $T_1\vee\cdots\vee T_p$ the tree obtained by identifying
the roots. Directly from \eqref{eq:homdensity},
\begin{align}
t_z\big(T^+,W\big)&=\int_0^1 W(z,z')\,t_{z'}(T,W)\,\dd z' ,\label{eq:graft}\\
t_z\big(T_1\vee\cdots\vee T_p,W\big)&=\prod_{i=1}^{p}t_z(T_i,W) ,\label{eq:wedge}
\end{align}
with $|E(T^+)|=|E(T)|+1$ and
$|E(T_1\vee\cdots\vee T_p)|=\sum_i|E(T_i)|$. Identity \eqref{eq:wedge} holds because
the integrals over the disjoint sets of non-root vertices factorise.

\begin{proof}[Proof of Theorem~\ref{thm:tree}]
We argue by strong induction on $j$, using only the definition \eqref{eq:formaldef};
no asymptotic expansion is assumed anywhere in this proof. For $j=0$ set
$\mc T_0=\{\bullet\}$, the one-vertex tree, with $t_z(\bullet,W)=1$ and
$\Up_\bullet:=\rho$; note that $\int\Up_\bullet=1$, the only member of the family with
nonzero average. For $j=1$ the only tuple in \eqref{eq:formaldef} is
$(p,\ell,j_1)=(1,0,0)$, giving
$\phi_{1,z}=-(\Id-P)^{-1}[P(A_{\rho,z}\rho)]'=\barW(z)\Up$ exactly as in the proof of
Theorem~\ref{thm:zerocoupling}; since $\mc T_1$ consists of the single edge and
$t_z(\text{edge},W)=\barW(z)$, this is \eqref{eq:tree} with
$\Up_{\text{edge}}=\Up$.

Assume \eqref{eq:tree} holds for all $i<j$. If
$\psi_{z}=\sum_{S\in\mc T_i}t_z(S,W)\Up_S$, then by \eqref{eq:A}, the linearity of
$v\mapsto A_v$ and \eqref{eq:graft},
\begin{equation}\label{eq:Agraft}
A_{\psi,z}=\int_0^1W(z,z')\sum_{S\in\mc T_i}t_{z'}(S,W)\,H_{\Up_S}\,\dd z'
=\sum_{S\in\mc T_i}t_z\big(S^+,W\big)\,H_{\Up_S} ,
\end{equation}
with $|E(S^+)|=i+1$.

Now consider a generic term of the recursion \eqref{eq:formaldef},
\[
\frac{(-1)^p}{p!}\Big[P\big(A_{\phi_{j_1},z}\cdots A_{\phi_{j_p},z}\,\phi_{\ell,z}\big)
\Big]^{(p)},\qquad p+\ell+\textstyle\sum_ij_i=j .
\]
By the induction hypothesis and \eqref{eq:Agraft},
\[
A_{\phi_{j_i},z}=\sum_{S_i\in\mc T_{j_i}}t_z\big(S_i^+,W\big)\,H_{\Up_{S_i}},
\qquad
\phi_{\ell,z}=\sum_{S_0\in\mc T_\ell}t_z(S_0,W)\,\Up_{S_0} .
\]
Multiplying out and using \eqref{eq:wedge}, the term equals
\[
\sum_{S_0,S_1,\dots,S_p}
t_z\big(S_0\vee S_1^+\vee\cdots\vee S_p^+,\,W\big)\cdot
\frac{(-1)^p}{p!}\Big[P\big(H_{\Up_{S_1}}\cdots H_{\Up_{S_p}}\,\Up_{S_0}\big)
\Big]^{(p)} ,
\]
and the tree $S_0\vee S_1^+\vee\cdots\vee S_p^+$ has
$\ell+\sum_i(j_i+1)=j$ edges. The second factor depends only on $f$ and $h$, because
$H_\psi$ depends only on $h$ and $\psi$, and because by induction the $\Up_{S_i}$
depend only on $f$ and $h$. Since $(\Id-P)^{-1}$ acts on the $\T$ variable only and is
independent of $z$, applying it to \eqref{eq:formaldef} and collecting all terms
contributing to a given $T\in\mc T_j$ yields \eqref{eq:tree} with
\begin{equation}\label{eq:UpT}
\Up_T=(\Id-P)^{-1}\!\!\!\!\sum_{\substack{(p,\ell,S_0,\dots,S_p)\in\mc D(T)}}\!\!\!\!
\frac{(-1)^p}{p!}\Big[P\big(H_{\Up_{S_1}}\cdots H_{\Up_{S_p}}\,\Up_{S_0}\big)
\Big]^{(p)} ,
\end{equation}
where $\mc D(T)$ is the set of tuples with $p\ge1$, $\ell\ge0$,
$S_i\in\mc T_{j_i}$, $p+\ell+\sum_ij_i=j$ and
$S_0\vee S_1^+\vee\cdots\vee S_p^+\cong T$ as rooted trees. That $\mc D(T)$ is finite
and non-empty is the content of Lemma~\ref{lem:decomp} below; we stress that the
decomposition is \emph{not} unique, which is why \eqref{eq:UpT} is a sum over
$\mc D(T)$ rather than a single term.

Each summand in \eqref{eq:UpT} is a $p$-th total derivative, hence of zero average, so
$(\Id-P)^{-1}$ applies and $\int_\T\Up_T=0$. For the regularity, note that
$H_{\Up_{S_i}}\in\mc C^{r}$ by Remark~\ref{rem:BVW11}(ii), so the only loss comes from
the $p$ derivatives: if $\Up_{S_0}\in W^{r-1-\ell,1}$ then the summand lies in
$W^{r-1-\ell-p,1}\subset W^{r-1-j,1}$ because $\ell+p\le j$. Since
$(\Id-P)^{-1}$ preserves $W^{m,1}$ for $m\le r-1$ by Proposition~\ref{prop:gap}, the
induction gives $\Up_T\in W^{r-1-j,1}$.
\end{proof}

\begin{lemma}[Root decompositions of a rooted tree]\label{lem:decomp}
Let $j\ge1$ and $T\in\mc T_j$, and let $B_1,\dots,B_q$ be the branches of $T$ at its
root, i.e.\ the connected components of $T$ minus the root, each rooted at its vertex
adjacent to the root of $T$ and re-attached to a new root; thus $q\ge1$ and
$B_i=S_i^+$ for a unique $S_i\in\mc T_{|E(B_i)|-1}$. Then the set $\mc D(T)$ of tuples
$(p,\ell,S_0,\dots,S_p)$ with $p\ge1$, $\ell\ge0$, $S_i\in\mc T_{j_i}$,
$p+\ell+\sum_{i\ge1}j_i=j$ and
$S_0\vee S_1^{+}\vee\cdots\vee S_p^{+}\cong T$ is non-empty and finite, with
\[
1\le |\mc D(T)|\le\sum_{p=1}^{q}\binom{q}{p}\,p!\ .
\]
\end{lemma}

\begin{proof}
A tuple in $\mc D(T)$ is obtained by selecting an ordered list of $p\ge1$ of the
branches $B_{i_1},\dots,B_{i_p}$, with pairwise distinct indices, to play the role of
the factors $S_1^+,\dots,S_p^+$, the remaining branches being assembled into
$S_0:=\bigvee_{i\notin\{i_1,\dots,i_p\}}B_i$ (the one-vertex tree if all branches are
selected). Conversely every tuple in $\mc D(T)$ arises this way, because an
isomorphism $S_0\vee S_1^+\vee\cdots\vee S_p^+\cong T$ must match the branches at the
root bijectively. Taking $p=1$ and $i_1=1$ shows $\mc D(T)\ne\emptyset$, and the
count of ordered selections gives the upper bound. The edge counts are automatic:
$\ell+\sum_{i\ge1}(j_i+1)=\sum_{i=1}^q|E(B_i)|=j$.
\end{proof}

\begin{remark}[On multiplicities]\label{rem:multiplicities}
The recursion \eqref{eq:formaldef} sums over \emph{ordered} tuples
$(j_1,\dots,j_p)$, since $\Delta^p$ is a $p$-th power, while \eqref{eq:tree} groups
terms by \emph{isomorphism classes} of rooted trees, and several ordered tuples may
produce isomorphic trees. No multiplicity can be mis-counted, because $\Up_T$ is not
given by an independent closed formula to be matched against a count: it is
\emph{defined} by \eqref{eq:UpT} as the sum of exactly those contributions whose tree
is isomorphic to $T$. The content of Theorem~\ref{thm:tree} is therefore the
factorisation of the $W$-dependence --- every contribution carrying the same tree
carries the same factor $t_z(T,W)$, because $t_z(\cdot,W)$ depends only on the
isomorphism class --- and not a combinatorial formula for the profiles. What the
theorem does assert about the profiles is that they are independent of $W$ and of $z$,
which is immediate from \eqref{eq:UpT}. For $j=2$ the decomposition is carried out
explicitly in Theorem~\ref{thm:secondorder}, and the resulting formula is confirmed
numerically in Remark~\ref{rem:numerics2}; the latter test separates the two profiles,
since $t_z(\Lambda_2,W)$ and $t_z(P_2,W)$ are not proportional on the two blocks.
\end{remark}

\begin{remark}[A falsifiable test of the tree structure at order three]\label{rem:order3}
At order $j$ the assertion of Theorem~\ref{thm:tree} is, in essence, a \emph{rank
condition} on the dependence of $\phi_{j,z}$ on the fibre, and such a condition can be
tested numerically without knowing any of the profiles $\Up_T$. This is worth doing at
$j=3$, an order at which we have not proved the expansion of $\phi^\star$ to exist.

There are four rooted trees with three edges. For a block graphon with $B$ blocks of
equal size and block matrix $(w_{ab})$, set $M:=(w_{ab}/B)$ and $d:=M\mathbf 1=\barW$;
then \eqref{eq:homdensity} evaluates the four densities, as vectors indexed by the
block, to
\[
d^{3},\qquad d\cdot(Md),\qquad M\big(d^{2}\big),\qquad M^{2}d ,
\]
powers and products being taken entrywise. These correspond, in order, to the star with
three leaves, to the tree whose root has two children one of which carries a child, to
the tree whose root has a single child carrying two children, and to the path. Hence
Theorem~\ref{thm:tree} predicts that for $B>4$ the $B$-tuple of functions
$\{\phi_{3,z_a}\}_{a\le B}$ is annihilated by every covector $c\in\R^{B}$ annihilating
those four vectors.

What follows is a \emph{conditional} consistency test, not a consequence of our
differentiability results, which reach second order only. Suppose that
$\phi^\star$ admits a genuine third-order expansion,
\begin{equation}\label{eq:thirdorderhyp}
\phi^\star_z(\alpha)=\rho+\alpha\phi_{1,z}+\alpha^{2}\phi_{2,z}
+\alpha^{3}\phi_{3,z}+o(\alpha^{3})\qquad\text{in }\W^{k},
\end{equation}
with $\phi_1,\phi_2$ the coefficients of Proposition~\ref{prop:formalgenuine}. Then the
centred odd quotient
\begin{equation}\label{eq:Q3}
Q_{3,\alpha}(z):=\frac{\phi^\star_z(\alpha)-\phi^\star_z(-\alpha)
-2\alpha\,\phi_{1,z}}{2\alpha^{3}}
\end{equation}
satisfies $Q_{3,\alpha}\to\phi_3$ in $\W^k$, with a remainder that is $o(1)$ and no
better: a rate would require a longer expansion. It is $Q_{3,\alpha}$, with
$\phi_{1,z}=\barW(z)\Up$ known exactly, that we compute.

We tested this with $B=6$, the dynamics of Remark~\ref{rem:numerics} and the block
matrix
\[
(w_{ab})=\begin{pmatrix}
0.178&0.660&0.322&0.643&0.796&0.592\\
0.182&0.738&0.638&0.663&0.574&0.447\\
0.116&0.598&0.815&0.720&0.388&0.622\\
0.943&0.674&0.870&0.854&0.402&0.786\\
0.381&0.148&0.220&0.430&0.757&0.708\\
0.851&0.182&0.159&0.415&0.574&0.680
\end{pmatrix},
\qquad
d=\barW=\begin{pmatrix}0.5318\\0.5403\\0.5432\\0.7548\\0.4407\\0.4768\end{pmatrix}^{\!\top},
\]
for which the four density vectors turn out to be linearly independent, so that the
annihilator is two dimensional. With $\phi_3$ replaced by $Q_{3,\alpha}$,
we computed the normalised residual
$\max_{c}\|c\cdot Q_{3,\alpha}\|_{\infty}\big/\|Q_{3,\alpha}\|_{\infty}$, the maximum
being over an orthonormal basis of the annihilator. It takes the values
\[
4.34\cdot10^{-5},\quad 1.07\cdot10^{-5},\quad 2.64\cdot10^{-6},\quad 1.33\cdot10^{-6}
\qquad\text{for }\alpha=0.2,\ 0.1,\ 0.05,\ 0.025 ,
\]
falling by factors $4.06$ and then $4.05$ before departing from that law at the
smallest $\alpha$. This is more than the vanishing we set out to test, and it is the
behaviour the theory predicts: if \eqref{eq:thirdorderhyp} holds to fifth order, with
remainder $o(\alpha^5)$, then
$Q_{3,\alpha}=\phi_3+\alpha^{2}\phi_5+o(\alpha^{2})$, the even orders cancelling in the
odd difference; the annihilator kills $\phi_3$ exactly but need not kill $\phi_5$,
because the rooted trees with five edges do not lie in the span of those with three.
The residual should therefore scale as $\alpha^2$ until a numerical floor is reached,
which is what is observed, the departure at $\alpha=0.025$ marking that floor.

For calibration, consider the density of the triangle rooted at $z$, which for a block
graphon is
\[
t_z(K_3,W)=\sum_{b,c}M_{zb}M_{bc}M_{cz} .
\]
It has relative component $5.6\cdot10^{-2}$ along the
same annihilator; a cycle contribution of magnitude comparable to the tree terms would
therefore have produced a residual three to four orders of magnitude larger than the
one observed, and would not have decayed quadratically.

Three caveats, and one diagnostic that we report because it corrected our own first
reading of these numbers. (i) The test probes a necessary condition: it would not
detect an error in the profiles $\Up_T$ themselves, only a contribution from a graph
that is not a rooted tree with three edges. (ii) It is conditional on
\eqref{eq:thirdorderhyp}, which we have not proved. (iii) The matrix of the four
density vectors, with columns normalised to unit Euclidean norm, has singular values
$1.977$, $0.301$, $1.27\cdot10^{-2}$, $2.6\cdot10^{-3}$ and hence condition number
$762$. That conditioning does not enter the test, because the annihilator is computed
from the exact block matrix displayed above rather than from simulation output; it
would matter for a least-squares recovery of the four profiles, which we do not
attempt. Finally, the diagnostic: the residual turns out to be insensitive to the
spatial grid, moving by less than a factor $1.3$ as $N$ ranges over
$256,\dots,2048$, and equally insensitive to the fixed-point tolerance over
$10^{-8},\dots,10^{-12}$. It is governed instead by the grid used to invert the lift
when assembling the transfer operator: with that grid at $1.5\cdot10^{4}$ points the
residual at $\alpha=0.025$ is $1.45\cdot10^{-3}$, at $6\cdot10^{4}$ points
$1.38\cdot10^{-4}$, and at $9.6\cdot10^{5}$ points the value quoted above. The
figures reported here use $N=512$ and $9.6\cdot10^{5}$ lift points; the script
performing all of these checks, including the resolution and tolerance studies, is
available from the author.
\end{remark}

\begin{remark}[Interpretation]\label{rem:interpretation}
Theorem~\ref{thm:tree} says that the perturbative expansion of the equilibrium state is
indexed by \emph{rooted trees}, with the graphon entering only through rooted
homomorphism densities and the dynamics only through universal profiles. This is a
dynamical counterpart of the classical fact that, for symmetric graphons, dense graph
convergence is equivalent to convergence of all homomorphism densities
\cite{LovaszSzegedy2006,LovaszBook} --- a counterpart in spirit only, since we use no
directed version of that equivalence (see Appendix~\ref{app:graphons}): the
first $m$ orders of the response see only the tree densities up to $m$ edges, so two
graphons with the same rooted tree densities up to order $m$ have the same
susceptibilities up to order $m$. In particular \emph{no formal susceptibility sees cycle
densities}: triangles, and more generally any homomorphism density of a graph
containing a cycle, are absent from \eqref{eq:tree} at every order. We emphasise what
this does and does not say. It says that the coefficients \eqref{eq:formaldef}, which
by Proposition~\ref{prop:formalgenuine} are the genuine first and second order
coefficients of $\phi^\star(\alpha)$, are functions of tree densities alone. It does
\emph{not} say that the equilibrium state $\phi^\star(\alpha)$ itself is a function of
tree densities: that would require the formal series to determine $\phi^\star$, which
is not known. Cycles may well influence the
equilibrium through effects invisible to this expansion.
\end{remark}

\section{Response to perturbations of the graphon}\label{sec:structural}

\begin{remark}[The uniform setting in $\eps$]\label{rem:epsuniform}
Throughout this section $\eps_0>0$ is fixed first, as in the statement of
Theorem~\ref{thm:structural}, and all the objects of
Sections~\ref{sec:spectral}--\ref{sec:differential} are formed once and for all with
$\|W\|_{\infty,1}$ replaced by $\|W\|_{\infty,1}+\eps_0\|V\|_{\infty,1}$ and
$\mathtt{var}_{p,L^1}(W)$ by $\mathtt{var}_{p,L^1}(W)+\eps_0\,\mathtt{var}_{p,L^1}(V)$.
In particular $\delta_0$, the auxiliary convex class $\mc V$ of \eqref{eq:Vclass}, the
constants $\xi$ and $K_r$ of Lemma~\ref{lem:uniform}, the bounds $M_k$ of
Proposition~\ref{prop:apriori} and the threshold $\alpha_2'$ are common to all the
kernels $W_\eps$, $|\eps|\le\eps_0$, and every fibre map
$F^{\eps}_{\alpha,\phi,z}$ with $|\eps|\le\eps_0$ and $|\alpha|<\alpha_2'$ lies in that
one class $\mc V$. This is the only reason the order of quantifiers in
Theorem~\ref{thm:structural} matters, and it is why $\eps_0$ is fixed before $\alpha$.
\end{remark}

\begin{proof}[Proof of Theorem~\ref{thm:structural}]
Write $\STO^{\eps}_\alpha$ for the STO associated with $W_\eps=W+\eps V$ and observe
that, by the linearity of \eqref{eq:A} in $W$,
\[
F^{\eps}_{\alpha,\phi,z}=f+\alpha A_{\phi,z}+\alpha\eps A^V_{\phi,z}
=F_{\alpha,\phi,z}+\alpha\eps A^V_{\phi,z},
\]
where $A^V$ is \eqref{eq:A} with $W$ replaced by $V$. We first check that $W_\eps$ is admissible in the sense of (A3), uniformly in $\eps$.
Boundedness is clear, and since $\mathtt{osc}_{L^1}(\,\cdot\,,\omega,r)$ is subadditive
in the kernel,
\begin{equation}\label{eq:varWeps}
\mathtt{var}_{p,L^1}(W+\eps V)\ \le\ \mathtt{var}_{p,L^1}(W)
+|\eps|\,\mathtt{var}_{p,L^1}(V)\ <\ \infty
\end{equation}
for every $|\eps|\le1$, with a bound independent of $\eps$. This is precisely why
$\mathtt{var}_{p,L^1}(V)<\infty$ is assumed: without it the perturbed kernel need not
lie in the class for which Theorem~\ref{thm:MPT} provides a fixed point, and the object
$\phi^{\star,\eps}$ would not be known to exist. Granting \eqref{eq:varWeps}, for
$|\eps|\le1$ and $\alpha$ fixed with $|\alpha|<\alpha_2$ all the hypotheses and
conclusions of Sections~\ref{sec:spectral}--\ref{sec:differential} hold for $W_\eps$
with constants uniform in $\eps$, upon replacing $\|W\|_{\infty,1}$ by
$\|W\|_{\infty,1}+\|V\|_{\infty,1}$ and shrinking $\alpha_0$ accordingly; in particular
$\phi^{\star,\eps}$ exists, is unique in the relevant class, and satisfies uniform a
priori bounds.

We now run the arguments of Theorems~\ref{thm:stability} and \ref{thm:LR} in the
variable $\eps$. Three facts are needed, and each is the exact analogue of what was
used there.

\emph{(a) Affine dependence on $\eps$ and the frozen derivative.} For $\phi$ fixed,
$\eps\mapsto F^\eps_{\alpha,\phi,z}$ is affine with
$F^{\eps}_{\alpha,\phi,z}-F^{\eps'}_{\alpha,\phi,z}=(\eps-\eps')\alpha A^V_{\phi,z}$.
Hence \eqref{eq:est1} and Lemma~\ref{lem:derivTO}, applied with
$\Delta=(\eps-\eps')\alpha A^V_{\phi,z}$ and $X=\alpha A^V_{\phi,z}$, give for
$\phi\in\Delta^{k+2}$
\begin{equation}\label{eq:epsLip}
\big\|\STO^{\eps}_\alpha\phi-\STO^{\eps'}_\alpha\phi\big\|_{k}
\le C|\alpha|\,|\eps-\eps'|\,\|\phi\|_{k+1},
\qquad
\partial_\eps\big(\STO^\eps_\alpha\phi\big)
=\alpha\,\mc K^{V,\eps}_{\alpha,\phi}\phi ,
\end{equation}
where, the frozen derivative in $\eps$ being taken at the \emph{perturbed} fibre map,
\begin{equation}\label{eq:KVeps}
\big(\mc K^{V,\eps}_{\alpha,\phi}v\big)_z
:=-\Big[L_{F^{\eps}_{\alpha,\phi,z}}\big(A^V_{v,z}\,\phi_z\big)\Big]' .
\end{equation}
This is \eqref{eq:KK} with $A$ replaced by $A^V$ in the perturbation direction and the
transfer operator taken at $F^{\eps}_{\alpha,\phi,z}$; at $\eps=0$ one has
$F^{0}_{\alpha,\phi,z}=F_{\alpha,\phi,z}$ and hence
$\mc K^{V,0}_{\alpha,\phi}=\mc K^V_{\alpha,\phi}$, which is the operator appearing in
\eqref{eq:structuralLR}, as that formula is evaluated at $\eps=0$. The distinction
matters only for the intermediate estimates, and for those we record that
$\mc K^{V,\eps}$ satisfies \eqref{eq:Kbound} with $\|W\|_{\infty,1}$ replaced by
$\|V\|_{\infty,1}$, by the proof of Proposition~\ref{prop:Kbounded}, with constants
uniform in $|\eps|\le\eps_0$ and $|\alpha|<\alpha_2'$ because the maps
$F^\eps_{\alpha,\phi,z}$ all lie in the class $\mc V$ formed, as in
Remark~\ref{rem:epsuniform}, with $\|W\|_{\infty,1}$ replaced by
$\|W\|_{\infty,1}+\eps_0\|V\|_{\infty,1}$.

\emph{(b) Contraction, uniformly in $\eps$.} Proposition~\ref{prop:contraction} holds
for $\STO^{\eps}_\alpha$ with constants independent of $|\eps|\le1$, since its proof
uses only Proposition~\ref{prop:gap} and \eqref{eq:est1}, both of which hold for the
family $\mc U_\eps$ with uniform constants. Repeating the proof of
Theorem~\ref{thm:stability} with $\eps,\eps'$ in place of $\alpha,\beta$ and using
the first half of \eqref{eq:epsLip} in place of \eqref{eq:alphaLip} gives
\begin{equation}\label{eq:epsstability}
\big\|\phi^{\star,\eps}-\phi^{\star,\eps'}\big\|_{k}\le C'|\eps-\eps'|
\qquad (1\le k\le r-2).
\end{equation}

\emph{(c) The $\phi$-differential is unchanged at $\eps=0$.} Since $W_0=W$, we have
$\STO^0_\alpha=\STO_\alpha$ and therefore
$D_\phi\STO^{0}_\alpha|_{\phi^\star}=\LL_{\alpha,\phi^\star}
+\alpha\mc K_{\alpha,\phi^\star}$ by Theorem~\ref{thm:diff}; more generally
Theorem~\ref{thm:diff} holds for $\STO^\eps_\alpha$ with
$\LL,\mc K$ built from $W_\eps$, with constants uniform in $\eps$.

With (a)--(c) in hand the proof of Theorem~\ref{thm:LR} transfers line by line. Writing
$u_\eps:=(\phi^{\star,\eps}-\phi^\star)/\eps$, decomposing
\[
\phi^{\star,\eps}-\phi^\star
=\big[\STO^{\eps}_\alpha\phi^{\star,\eps}-\STO^{\eps}_\alpha\phi^\star\big]
+\big[\STO^{\eps}_\alpha\phi^{\star}-\STO^{0}_\alpha\phi^\star\big] ,
\]
applying Theorem~\ref{thm:diff} to the first bracket, \eqref{eq:epsstability} at level
$k+1$ to bound the quadratic remainder by $O(\eps^2)$, \eqref{eq:epsLip} to the second
bracket, and Proposition~\ref{prop:resolvent} to invert, we obtain
\[
\mc R_{\alpha,\phi^\star}\ \partial_\eps\phi^{\star,\eps}\big|_{0}
=\alpha\,\mc K^V_{\alpha,\phi^\star}\phi^\star ,
\]
which is \eqref{eq:structuralLR}.

For \eqref{eq:structuralLRsmall}, let $\alpha\to0$. By Theorem~\ref{thm:stability},
$\|\phi^\star-\rho\|_{k+1}=O(\alpha)$, and by \eqref{eq:est1},
$\|\LL_{\alpha,\phi^\star}-P\|_{\W^{k+1}\to\W^{k}}=O(\alpha)$. Hence
\[
\big(\mc K^V_{\alpha,\phi^\star}\phi^\star\big)_z
=-\Big[P\big(A^V_{\rho,z}\,\rho\big)\Big]'+O(\alpha)
=-\overline V(z)\big[P(H_\rho\rho)\big]'+O(\alpha)\qquad\text{in }\W^{k}.
\]
Write $B_\alpha:=\mc K^V_{\alpha,\phi^\star}\phi^\star$ and
$B_0:=-\overline V(\cdot)\big[P(H_\rho\rho)\big]'=\mc K^V_{0,\rho}\rho$, so that the
display reads $\|B_\alpha-B_0\|_{k}=O(\alpha)$. It would now be illegitimate to
write $\mc R_{\alpha,\phi^\star}^{-1}=(\Id-P)^{-1}+O(\alpha)$ in $\mc L(\W^k_0)$, since
the one-derivative loss in \eqref{eq:est1} does not permit it
(Remark~\ref{rem:nodensity}). It is not needed: the source term is one derivative
better than required. Indeed Proposition~\ref{prop:Kbounded}, applied to $\mc K^V$ at
level $k+1$ --- legitimate because $k+1\le r-2$ and $\phi^\star\in\Delta^{k+2}$ by
Proposition~\ref{prop:apriori} --- gives
\[
\|B_0\|_{k+1}\le C_{\mc K}\,\|\rho\|_{k+2}\,\|\rho\|_{0}<\infty ,
\]
so that $B_0\in\W^{k+1}_0$. Decomposing
\[
\mc R_{\alpha,\phi^\star}^{-1}B_\alpha-(\Id-P)^{-1}B_0
=\mc R_{\alpha,\phi^\star}^{-1}\big(B_\alpha-B_0\big)
+\Big(\mc R_{\alpha,\phi^\star}^{-1}-(\Id-P)^{-1}\Big)B_0 ,
\]
the first term is $O(\alpha)$ in $\W^{k}$ by the uniform bound
\eqref{eq:resolventbound}, and the second is $O(\alpha)$ by
Lemma~\ref{lem:resvector}, applicable precisely because $B_0\in\W^{k+1}_0$. Since
$(\Id-P)^{-1}B_0=\overline V(z)\,\Up$ by \eqref{eq:Upsilon}, multiplying by the
prefactor $\alpha$ in \eqref{eq:structuralLR} gives \eqref{eq:structuralLRsmall}. The
extra derivative demanded by Lemma~\ref{lem:resvector} is thus supplied, at no cost, by
the smoothing property \eqref{eq:Kbound} of the very operator that produces the source
term.
\end{proof}

\begin{remark}[On the admissibility of the perturbed kernel]\label{rem:admissible}
Theorem~\ref{thm:structural} requires of $V$ exactly what (A3) requires of $W$: that it
be a bounded measurable kernel of finite $p$-variation. It does \emph{not} require
$W_\eps=W+\eps V$ to be non-negative, and indeed for a signed $V$ non-negativity
generally fails for every $\eps\neq0$ at points where $W$ vanishes, however small
$|\eps|$ is. This costs nothing here, because none of (A1)--(A3),
Theorem~\ref{thm:stability}, Theorem~\ref{thm:LR} or
\cite[Thm.~3.2]{MPT2025} uses non-negativity: as declared in
Section~\ref{sec:setting}, a graphon in this paper is a bounded measurable real kernel.
What does require $W\ge0$ is the \emph{interpretation} of $\barW$ as a degree and the
monotonicity statement of item (ii) after Theorem~\ref{thm:zerocoupling}; those
readings are therefore not available for a general perturbed kernel $W_\eps$, even
though \eqref{eq:structuralLR} and \eqref{eq:structuralLRsmall} are.
\end{remark}

\begin{remark}
Theorem~\ref{thm:structural} answers (Q2) and provides the graphon analogue of the
optimal coupling problem of \cite{Galatolo2022,GalatoloPollicott2017,MacKay2018}. Its
message is that, in the weak coupling regime, \emph{the only structural quantity a
first-order intervention can act upon is the local degree}. An intervention designed to
change the behaviour of a specific node should either add degree at that node or, if
degree is constrained, exploit the second-order term and act on the degrees of its
neighbours. We note that sensitivity analysis on graphons has been considered in other
settings, for instance for SIS epidemics \cite{VizueteFrascaGarin2020}, where the
relevant quantity is spectral rather than degree-based; the degree-proportionality we
obtain is specific to the weak-coupling expansion of the self-consistent equilibrium.
\end{remark}

\section{A conditional prediction for finite networks}\label{sec:finiteN}

Everything in this section is conditional on Assumption~\ref{ass:exchange} below,
which is not proved here and which is, in the heterogeneous deterministic setting, an
open problem; the all-to-all case with additive noise is settled in
\cite{TenagliaTanzi2026}.

Let $G^{(N)}$ be a sequence of directed graphs whose step graphons $W^{(N)}$ satisfy
hypothesis (H2) of \cite{MPT2025} with limit $W$, and let $\mu^{(N)}$ satisfy the
concentration hypothesis (H1). By \cite[Thm.~3.1]{MPT2025}, for each \emph{fixed}
$t\in\N$ the one-node marginals of $(F^t)_*\mu^{(N)}$ at node $i=\lceil z_*N\rceil$
converge weakly to $(\STO^t_\alpha\nu)_{z_*}$.

\begin{assumption}[Exchange of limits]\label{ass:exchange}
There is a sequence $t_N\to\infty$ such that, for every $z_*$ in the set $E$ of
hypothesis (H2) and $i=\lceil z_*N\rceil$, the one-node distributions
$\Pi_i(F^{t_N})_*\mu^{(N)}$ converge weakly to $\phi^\star_{z_*}(\alpha)\,\dd\Leb$ as
$N\to\infty$, uniformly for $|\alpha|<\alpha_2$.
\end{assumption}

\begin{assumption}[Convergence of the normalised degrees]\label{ass:degrees}
With $z_{i(N)}$ denoting the representative of the $i(N)$-th row of the step graphon
$W^{(N)}$,
\begin{equation}\label{eq:degreeconv}
\int_0^1 W^{(N)}\big(z_{i(N)},y\big)\dd y\ \longrightarrow\
\int_0^1 W(z_*,y)\dd y=\barW(z_*) .
\end{equation}
Equivalently, since $W^{(N)}$ is the step graphon of the adjacency matrix $A^{(N)}$,
the normalised in-degrees of the selected vertices converge:
$k_{i(N)}/N\to\barW(z_*)$.
\end{assumption}

\begin{remark}[Why Assumption~\ref{ass:degrees} is not implied by (H2)]
\label{rem:H2degrees}
It is tempting to read \eqref{eq:degreeconv} off hypothesis (H2) of \cite{MPT2025}, and
earlier versions of this work did. That is wrong. Hypothesis (H2) states that
\[
\int_0^1 W^{(N)}(z_*,z)\int_\T h(x,y)\dd\nu_z(y)\dd z\ \longrightarrow\
\int_0^1 W(z_*,z)\int_\T h(x,y)\dd\nu_z(y)\dd z
\qquad\text{for all }x\in\T ,
\]
That is, it tests the rows $W^{(N)}(z_*,\cdot)$ only against the family of functions
\[
z\ \longmapsto\ \int_\T h(x,y)\,\dd\nu_z(y),
\]
indexed by $x\in\T$. Should that family contain no function of nonzero mean, (H2)
carries no information at all about the row sums.
The extreme case is instructive: for $h\equiv0$ hypothesis (H2) holds for
\emph{every} pair of sequences, so one may take $W\equiv0$ together with
$W^{(N)}\equiv1$, which gives $\barW\equiv0$ while $k_i/N=1$ for every $i$ and every
$N$. This is not a pathology confined to $h\equiv0$:
if $\nu_z=\nu$ does not depend on $z$, then (H2) reduces to
\[
H_\nu(x)\int_0^1W^{(N)}(z_*,z)\dd z\ \longrightarrow\ H_\nu(x)\,\barW(z_*)
\qquad\text{for all }x\in\T,
\]
which is
vacuous exactly when $H_\nu\equiv0$ --- the degenerate situation catalogued in
\eqref{eq:degenerate}, and the one in which the first-order response vanishes anyway.
When $H_\nu\not\equiv0$ for a $z$-independent $\nu$, (H2) does give
\eqref{eq:degreeconv}; this is how the Erd\H{o}s--R\'enyi case is handled in
\cite[\S3]{MPT2025}. We therefore impose \eqref{eq:degreeconv} as a separate
hypothesis. It is in any case the natural one here, being precisely the convergence of
the quantity whose relation to the susceptibility the corollary asserts, and it holds
as soon as $\|W^{(N)}(z_{i(N)},\cdot)-W(z_*,\cdot)\|_{L^1}\to0$, which is the first
sufficient condition for (H2) given in \cite[Lem.~4.1]{MPT2025}. A reader unwilling to
assume it may simply delete the second equality in \eqref{eq:finiteNexp} and read the
corollary as a statement about $\barW(z_*)$.
\end{remark}

\begin{corollary}[Degree-proportional susceptibility of nodes]\label{cor:finiteN}
Let $c\in\mc C^0(\T,\R)$, let $z_*\in E$ be a continuity point of $\barW$, where $E$ is
the set of hypothesis (H2) appearing in Assumption~\ref{ass:exchange}, suppose
Assumption~\ref{ass:degrees} holds at $z_*$, let
$i=i(N):=\lceil z_*N\rceil$ and put
$\langle c\rangle^{(N)}_i(\alpha):=\int_\T c\,\dd\big[\Pi_i(F^{t_N})_*\mu^{(N)}\big]$,
with $t_N$ as in Assumption~\ref{ass:exchange}. Then, with
$k_i=\sum_{j}A^{(N)}_{ij}$,
\begin{equation}\label{eq:finiteNexp}
\lim_{N\to\infty}\langle c\rangle^{(N)}_{i(N)}(\alpha)
=\int_\T c\,\rho\;+\;\alpha\,\barW(z_*)\int_\T c\,\Up\;+\;o(\alpha),
\qquad
\barW(z_*)=\lim_{N\to\infty}\frac{k_{i(N)}}{N} .
\end{equation}
Consequently the limiting stationary expectation at a node is, to first order in
$\alpha$, an affine function of the node's asymptotic normalised in-degree, with a
slope $\int_\T c\,\Up$ that is the same for all nodes and all networks.
\end{corollary}

\begin{proof}
By Assumption~\ref{ass:exchange} the measures $\Pi_{i(N)}(F^{t_N})_*\mu^{(N)}$ converge
weakly to $\phi^\star_{z_*}(\alpha)\,\dd\Leb$, so for $c$ continuous and bounded
$\langle c\rangle^{(N)}_{i(N)}(\alpha)\to\int_\T c\,\phi^\star_{z_*}(\alpha)$. Now
apply Theorem~\ref{thm:zerocoupling} and Corollary~\ref{cor:GK} to the limit. Finally
the identity $\barW(z_*)=\lim_N\frac1N\sum_jA^{(N)}_{i(N)j}$ is exactly
Assumption~\ref{ass:degrees}, in the equivalent form recorded there; it is \emph{not}
deduced from (H2), which by Remark~\ref{rem:H2degrees} does not imply it.
\end{proof}

\begin{remark}[Two caveats]\label{rem:finiteNcaveats}
Both restrictions in Corollary~\ref{cor:finiteN} are needed. First, weak convergence
of measures controls $\int c$ only for continuous $c$; to allow bounded measurable
observables one would have to strengthen Assumption~\ref{ass:exchange} to convergence
in total variation of the one-node marginals. Second, \eqref{eq:finiteNexp} is an
expansion of the \emph{limiting} expectation: weak convergence uniform in $\alpha$
does not by itself give convergence of the $\alpha$-derivatives, so we do not assert
that the finite-$N$ susceptibility $\partial_\alpha\langle c\rangle^{(N)}_i$
converges to the slope in \eqref{eq:finiteNexp}. Establishing that would require a
remainder estimate uniform in $N$ in a topology allowing differentiation in $\alpha$.
\end{remark}

\begin{remark}[On Assumption~\ref{ass:exchange}]\label{rem:exchange}
Assumption~\ref{ass:exchange} is not implied by \cite[Thm.~3.1]{MPT2025}, which
controls finitely many iterates only; it amounts to a uniform-in-time propagation of
chaos statement for heterogeneous deterministic networks. The corresponding problem for \emph{globally coupled} maps has
recently seen substantial progress: finite-time propagation of chaos for deterministic
uniformly expanding globally coupled maps is proved in \cite{Tanzi2023}, Gaussian
fluctuations around the mean-field limit in \cite{Cheng2026}, and uniform-in-time
propagation of chaos --- in the Dobrushin--Wasserstein distance, and in total variation
on fixed marginals --- has been obtained by Tenaglia and Tanzi
\cite{TenagliaTanzi2026} for mean-field coupled maps \emph{with additive noise} whose
density is bounded below, where the noise supplies a dimension-free contraction. Both
results are for all-to-all coupling; extending them to graphon couplings, and removing
the noise in the deterministic case, is in our view the most important problem left
open by the present work. See also \cite{AllmeierGast2024} for quantitative accuracy
of the graphon mean-field approximation in a stochastic setting.

We state Corollary~\ref{cor:finiteN} nonetheless because it is the form in which the
theory becomes testable: it predicts a strictly linear relation between node degree and
node susceptibility, with a universal slope, which can be checked numerically on any
dense heterogeneous network of expanding maps. A natural experiment is the pair
$(f,h)$ of Remark~\ref{rem:numerics}, namely
$f(x)=2x+\tfrac1{10}\sin(2\pi x)$ and
$h(x,y)=\sin(2\pi(y-x))+\tfrac12\cos(2\pi(x+2y))$, for which
$\|\Up\|_{\infty}\approx8.0\cdot10^{-2}$, together with a rank-one graphon with
heavy-tailed $g$ as in Example~\ref{ex:rank1}. The nonlinearity of $f$ is not a
convenience but a necessity: this same coupling satisfies
$\int_\T h(x,y)\dd y=0$, so by \eqref{eq:degenerate} it would produce
$\Up\equiv0$, and hence no measurable effect whatsoever, if $f$ were replaced by the
doubling map.
\end{remark}

\begin{appendices}

\section{Background material}\label{app:basics}

\subsection{Transfer operators}\label{app:TO}

Let $F:\T\to\T$ be $\mc C^1$ with $|F'|>0$. Its transfer operator with respect to
Lebesgue measure is determined by the duality relation
\begin{equation}\label{eq:duality}
\int_\T (g\circ F)\,\psi\;\dd x=\int_\T g\,(L_F\psi)\;\dd x
\qquad \text{for all } g\in L^\infty(\T),\ \psi\in L^1(\T),
\end{equation}
equivalently by $(L_F\psi)(x)=\sum_{y\in F^{-1}(x)}\psi(y)/|F'(y)|$. On measures $L_F$
is the push-forward. It is a Markov operator, so $\|L_F\psi\|_{L^1}\le\|\psi\|_{L^1}$,
and $\psi$ is an invariant density iff $L_F\psi=\psi$. If $F\in\mc C^{k+1}$ is
uniformly expanding, $L_F$ preserves $W^{k,1}(\T)$ and satisfies \eqref{eq:LY2}; if in
addition $F$ is topologically mixing, which is automatic for expanding circle maps of
degree $\ge2$, then $L_F$ has a spectral gap on $W^{k,1}$ for $k\ge1$ and
$(\Id-L_F)^{-1}=\sum_{n\ge0}L_F^n$ is bounded on the zero-average subspace. See
\cite{BoyarskyGora,Baladi2000,Liverani1995,Viana1997} and \cite[\S7]{Galatolo2022}. The
regularity of the invariant density, $\rho_F\in\mc C^{r-1}$ for $F\in\mc C^{r}$, is due
to Krzy\.zewski and Szlenk \cite{KrzyzewskiSzlenk1969,Krzyzewski1977}.

We shall use repeatedly the identity, valid for $F\in\mc C^2$ expanding and
$\psi\in W^{1,1}$,
\begin{equation}\label{eq:LYidentity}
(L_F\psi)'=L_F\Big(\frac{\psi'}{F'}\Big)-L_F\Big(\frac{F''}{(F')^2}\,\psi\Big) ,
\end{equation}
which follows from \eqref{eq:duality} applied to $g'$ together with
$\big((g\circ F)/F'\big)'=(g'\circ F)-(g\circ F)F''/(F')^2$.

\subsection{Cones and the Hilbert metric}\label{app:cones}

For $a>0$ let $\mc V_a=\{\psi\in\mc C^0(\T,\R^+):\ \psi(x)/\psi(y)\le e^{a|x-y|}\}$ and
let $\theta_a$ be the associated Hilbert projective metric. A linear operator mapping a
cone with finite diameter strictly inside itself is a contraction for $\theta_a$
(Birkhoff); for transfer operators of uniformly expanding maps this yields exponential
decay of correlations with rates depending only on $\inf|F'|$ and on the distortion
$\sup|F''/(F')^2|$ \cite{Liverani1995,Viana1997}. This is the mechanism behind
\cite[Prop.~4.5 and Lem.~4.5]{MPT2025}; in the present paper we use the resulting
memory-loss estimate (Lemma~\ref{lem:MLMPT}) as a black box, as the base case of
Proposition~\ref{prop:gap}.

\subsection{Graphons, degrees and homomorphism densities}\label{app:graphons}

A graphon is a measurable $W:[0,1]^2\to\R$, not assumed symmetric, so that directed
graphs are covered. A graph $G$ on $\{1,\dots,N\}$ with adjacency matrix $A$ is
represented by the step graphon $W^{(N)}|_{I_i\times I_j}\equiv A_{ij}$ with
$I_i=((i-1)/N,i/N]$ (and $I_1=[0,1/N]$); then
$\frac1NA_{ij}=\int_{I_j}W^{(N)}(z,z')\dd z'$ for $z\in I_i$, which is the identity used
to rewrite \eqref{eq:finite} in graphon form, see \cite[\S2.2]{MPT2025}. The degree
function is $\barW(z)=\int_0^1W(z,z')\dd z'$, so that $\barW^{(N)}|_{I_i}=k_i/N$; the
edge density is $\int\!\!\int W$.

Since $W$ is not assumed symmetric, homomorphism densities must be understood in the
\emph{directed} sense: each edge carries an orientation and contributes the factor
$W(z_u,z_v)$ in that order. Thus for a finite \emph{digraph} $H$ on $q$ vertices,
with arc set $E^{\rightarrow}(H)$,
\begin{equation}\label{eq:homdirected}
t(H,W)=\int_{[0,1]^q}\ \prod_{(u,v)\in E^{\rightarrow}(H)}W(z_u,z_v)
\ \prod_{v\in V(H)}\dd z_v ,
\end{equation}
and the rooted densities \eqref{eq:homdensity} are the local versions of this, the
edges of a rooted tree being oriented away from the root --- which is the convention
fixed in \eqref{eq:homdensity} by writing $(u,v)\in E(T)$ with $u$ the parent of $v$.
For an undirected $H$ and a symmetric $W$, \eqref{eq:homdirected} reduces to the
classical quantity, and in that symmetric setting dense graph convergence in the cut
metric is equivalent to the convergence of $t(H,\cdot)$ for all $H$
\cite{LovaszSzegedy2006,LovaszBook,BorgsChayes2008}. We invoke that equivalence only
as motivation, and only for symmetric kernels; no directed analogue of it is used or
asserted anywhere in this paper, where the only graph theory needed is the pair of
identities \eqref{eq:graft} and \eqref{eq:wedge} for explicitly oriented rooted
trees.

\subsection{Disintegration and the fibred setting}\label{app:disint}

For $\nu\in\M_{1,\Leb}([0,1]\times\T)$ the disintegration $\{\nu_z\}$ with respect to
the partition $\{\{z\}\times\T\}$ satisfies: $z\mapsto\nu_z(A)$ is measurable and
$\nu(A)=\int_0^1\nu_z(A)\dd z$, and $\nu_z$ is a.s.\ a probability measure because the
base marginal is $\Leb$. When $\nu$ is absolutely continuous we write
$\phi_z=\dd\nu_z/\dd\Leb$. The STO acts fibrewise on the disintegration,
$(\STO\nu)_z=(F_{\nu,z})_*\nu_z$, but the maps $F_{\nu,z}$ depend on \emph{all} fibres
through \eqref{eq:A}; this is the entire difficulty, and the entire interest, of the
heterogeneous case.

\subsection{The $p$-variation seminorms of \texorpdfstring{\cite{MPT2025}}{[MPT]}}
\label{app:pvar}

For $\phi\in L^1([0,1],L^1(\T))$, $\omega\in[0,1]$ and $\varrho>0$ put
$\mathtt{osc}_{\BV^1}(\phi,\omega,\varrho)=\esssup_{z,\bar z\in B(\omega,\varrho)}
|\phi_{\bar z}-\phi_z|_{\BV^1}$ and, for $p\in(0,1]$,
\[
\mathtt{var}_{p,\BV^1}(\phi)=\sup_{\varrho>0}\frac{1}{\varrho^p}\int_0^1
\mathtt{osc}_{\BV^1}(\phi,\omega,\varrho)\,\dd\omega ,
\]
with the analogous definition of $\mathtt{var}_{p,L^1}(W)$ for graphons. These
seminorms are used in \cite{MPT2025}, following \cite{KellerPV,Galatolo2018}, to obtain
compactness of the admissible set of densities; they are finite for H\"older or
Lipschitz kernels, block graphons, spatially decaying kernels and constant graphons.
We use them only in Proposition~\ref{prop:zregularity}.

\section{Exact perturbative identities and calculus lemmas}\label{app:calculus}

\begin{lemma}[Products and compositions]\label{lem:calculus}
There is $M>0$ such that for $1\le k\le r$, $\Phi\in\mc C^k(\T)$ and
$\Psi\in W^{k,1}(\T)$,
\[
\|\Phi\Psi\|_{W^{k,1}}\le M\|\Phi\|_{\mc C^k}\|\Psi\|_{W^{k,1}} ,
\]
so that $W^{k,1}(\T)$ is a module over $\mc C^{k}(\T)$. If $\Theta$ is a $\mc C^{k}$
diffeomorphism of $\T$ with $\inf|\Theta'|>0$ then
$\|\Psi\circ\Theta\|_{W^{k,1}}\le M\|\Psi\|_{W^{k,1}}\max_{1\le\ell\le k}
\|\Theta\|^{\ell}_{\mc C^k}(\inf|\Theta'|)^{-1}$.
\end{lemma}

\begin{proof}
Leibniz, Fa\`a di Bruno and the change of variables formula.
\end{proof}

\begin{proof}[Proof of Lemma~\ref{lem:exact}]
Let $g\in\mc C^{m+2}(\T,\R)$. Since $F+t\Delta\in\mc V$ for $t\in[0,1]$, by the
convexity of $\mc V$ established in Lemma~\ref{lem:uniform} --- note that $\mc U$
itself would not do, being non-convex by Remark~\ref{rem:notconvex} ---
(Lemma~\ref{lem:uniform}), the fundamental theorem of calculus gives
\[
g\circ(F+\Delta)-g\circ F=\int_0^1\big(g'\circ(F+t\Delta)\big)\,\Delta\,\dd t .
\]
Multiplying by $u$, integrating over $\T$ and using \eqref{eq:duality} twice,
\[
\int_\T g\,\big(L_{F+\Delta}u-L_Fu\big)
=\int_0^1\!\!\int_\T g'\,L_{F+t\Delta}(\Delta u)\,\dd t
=-\int_0^1\!\!\int_\T g\,\big[L_{F+t\Delta}(\Delta u)\big]'\,\dd t ,
\]
the last step by integration by parts on $\T$, which has no boundary. As $g$ is
arbitrary this is \eqref{eq:exact1}. Likewise, Taylor's formula with integral remainder
gives
\[
g\circ(F+\Delta)=g\circ F+(g'\circ F)\Delta
+\int_0^1(1-t)\big(g''\circ(F+t\Delta)\big)\Delta^2\,\dd t ,
\]
and the same manipulation, now with two integrations by parts in the remainder, yields
\eqref{eq:exact2}.

For the estimates, recall that for $v\in L^1(\T)$ and $0\le i$,
$\|v^{(i)}\|_{L^1}=\sup\{\int_\T g^{(i)}v:\ g\in\mc C^i,\ \|g\|_\infty\le1\}$, so that
$\|v'\|_{W^{m,1}}\le\|v\|_{W^{m+1,1}}$ and $\|v''\|_{W^{m,1}}\le\|v\|_{W^{m+2,1}}$.
Hence, by \eqref{eq:exact1}, Lemma~\ref{lem:LY} and Lemma~\ref{lem:calculus},
\[
\big\|L_{F+\Delta}u-L_Fu\big\|_{W^{m,1}}
\le\sup_{t\in[0,1]}\big\|L_{F+t\Delta}(\Delta u)\big\|_{W^{m+1,1}}
\le C\,\|\Delta\|_{\mc C^{m+1}}\|u\|_{W^{m+1,1}} ,
\]
which is \eqref{eq:est1}; this uses Lemma~\ref{lem:LY} at level $m+1$, whence the
hypothesis $m+1\le r-1$ in part (i). Similarly \eqref{eq:exact2} gives \eqref{eq:est2}
with $\|\Delta^2\|_{\mc C^{m+2}}\le\|\Delta\|^2_{\mc C^{m+2}}$ up to a combinatorial
constant, now using Lemma~\ref{lem:LY} at level $m+2$, whence $m+2\le r-1$ in part
(ii).
\end{proof}

\begin{proof}[Proof of Lemma~\ref{lem:expansionL}]
Identical, expanding $g\circ(f+\Delta)$ to order $p_{\max}$ with integral remainder and
performing $p$ integrations by parts in the term of order $p$, and $p_{\max}+1$ in the
remainder.
\end{proof}

\section{Long proofs}\label{app:longproofs}

\subsection{Proof of Lemma~\ref{lem:LY}}

Let $F\in\mc U$ and $\lambda=\xi^{-1}$. By \eqref{eq:LYidentity} and
$\|L_Fv\|_{L^1}\le\|v\|_{L^1}$,
\[
\|(L_F\psi)'\|_{L^1}\le\Big\|\frac{\psi'}{F'}\Big\|_{L^1}
+\Big\|\frac{F''}{(F')^2}\Big\|_\infty\|\psi\|_{L^1}
\le\lambda\|\psi'\|_{L^1}+K'\|\psi\|_{L^1} ,
\]
with $K'=\sup_{\mc U}\|F''/(F')^2\|_\infty<\infty$ by Lemma~\ref{lem:uniform}. This is
the case $k=1$, $n=1$. We stress that no factor equal to the degree $d$ appears; this
is the advantage of the $L^1$-based scale, see Remark~\ref{rem:nodegree}.

For general $k$ we use the following explicit form of the iterated identity.

\medskip
\noindent\emph{Claim.} For $1\le k\le r-1$ and $\psi\in W^{k,1}(\T)$,
\begin{equation}\label{eq:LYiterated}
\big(L_F\psi\big)^{(k)}
=L_F\!\left(\frac{\psi^{(k)}}{(F')^{k}}\right)
+\sum_{i=0}^{k-1}L_F\big(q_{k,i}\,\psi^{(i)}\big) ,
\end{equation}
where each $q_{k,i}$ is a universal polynomial in $1/F'$ and
$F'',F''',\ldots,F^{(k+1)}$, the same polynomial for every $F$.

\smallskip
\noindent\emph{Proof of the Claim.} For $k=1$, \eqref{eq:LYidentity} is
\eqref{eq:LYiterated} with $q_{1,0}=-F''/(F')^2$. Assume \eqref{eq:LYiterated} for
$k$ and apply \eqref{eq:LYidentity} to each summand, in the form
$\big[L_F v\big]'=L_F\big(v'/F'\big)-L_F\big(\theta_F\,v\big)$ with
$\theta_F:=F''/(F')^2$. For the leading summand, with
$v=\psi^{(k)}/(F')^{k}$,
\[
\frac{v'}{F'}=\frac{\psi^{(k+1)}}{(F')^{k+1}}
+\psi^{(k)}\,\frac{1}{F'}\Big(\frac{1}{(F')^{k}}\Big)' ,
\]
whose first term is the new leading term and whose second term has a coefficient of the
required form; the term $-L_F(\theta_Fv)$ likewise contributes a coefficient times
$\psi^{(k)}$. For a summand with $v=q_{k,i}\psi^{(i)}$, $i\le k-1$, one gets
$v'/F'=(q_{k,i}'\psi^{(i)}+q_{k,i}\psi^{(i+1)})/F'$ with $i+1\le k$, and
$-\theta_Fv$, again of the required form. This proves the Claim, and shows that the
$q_{k,i}$ involve derivatives of $F$ of order at most $k+1$. \hfill$\diamond$

\medskip
Since $L_F$ is Markov, \eqref{eq:LYiterated} gives
\[
\|(L_F\psi)^{(k)}\|_{L^1}\le\Big\|\frac{1}{(F')^{k}}\Big\|_\infty
\|\psi^{(k)}\|_{L^1}+\sum_{i=0}^{k-1}\|q_{k,i}\|_\infty\|\psi^{(i)}\|_{L^1}
\le\lambda^{k}\|\psi^{(k)}\|_{L^1}+\tilde B_k\|\psi\|_{W^{k-1,1}} ,
\]
where $\tilde B_k:=\sum_{i<k}\sup_{\mc U}\|q_{k,i}\|_\infty$, which is finite because
each $q_{k,i}$ is a polynomial in $1/F'$ and in $F'',\ldots,F^{(k+1)}$, all of them
bounded on $\mc U$ by Lemma~\ref{lem:uniform} (here $k+1\le r$).
Together with the inductive hypothesis
$\|L_F\psi\|_{W^{k-1,1}}\le A_k\|\psi\|_{W^{k-1,1}}$ this gives
\[
\|L_F\psi\|_{W^{k,1}}\le\lambda^{k}\|\psi\|_{W^{k,1}}
+B_k\|\psi\|_{W^{k-1,1}},\qquad B_k=\tilde B_k+A_k .
\]
It remains to compose $n$ such inequalities. We write this out, because the maps are
allowed to differ at each step and the lower-order term must be seen to stay bounded
uniformly over all admissible compositions. The argument is an induction on $k$.

For $k=1$ the lower-order bound is the $L^1$ contraction
$\|L(n)\psi\|_{L^1}\le\|\psi\|_{L^1}$, valid for every composition because each
$L_{F_i}$ is Markov; this is \eqref{eq:LY1} with $A_1=1$. Suppose then that
\eqref{eq:LY1} holds at level $k$, i.e.\ that
\[
\sup_{n\ge0}\ \sup_{F_1,\dots,F_n\in\mc U}\ \|L(n)\psi\|_{W^{k-1,1}}
\le A_k\|\psi\|_{W^{k-1,1}}
\qquad\text{for all }\psi\in W^{k-1,1}(\T),
\]
and write $\psi_j:=L(j)\psi$, with $\psi_0=\psi$. Applying the one-step inequality
above to $\psi_{j}$ --- with the map $F_{j+1}$, whichever element of $\mc U$ it is, and
with the constants $\lambda^k$ and $B_k$, which do not depend on that choice --- and
iterating from $j=n-1$ down to $j=0$,
\[
\|\psi_n\|_{W^{k,1}}
\le\lambda^{k}\|\psi_{n-1}\|_{W^{k,1}}+B_k\|\psi_{n-1}\|_{W^{k-1,1}}
\le\lambda^{kn}\|\psi\|_{W^{k,1}}
+B_k\sum_{j=0}^{n-1}\lambda^{k(n-1-j)}\,\|\psi_j\|_{W^{k-1,1}} .
\]
Every intermediate term $\psi_j$ is itself a composition of maps in $\mc U$, so the
inductive hypothesis bounds $\|\psi_j\|_{W^{k-1,1}}\le A_k\|\psi\|_{W^{k-1,1}}$
uniformly in $j$ and in the maps used, and summing the geometric series gives
\begin{equation}\label{eq:LYcomposed}
\|\psi_n\|_{W^{k,1}}\le\lambda^{kn}\|\psi\|_{W^{k,1}}
+\frac{A_kB_k}{1-\lambda^{k}}\,\|\psi\|_{W^{k-1,1}} ,
\end{equation}
which is \eqref{eq:LY2} after renaming the constant
$B_k\rightsquigarrow A_kB_k(1-\lambda^k)^{-1}$. Taking the supremum over $n$ in
\eqref{eq:LYcomposed} and using $\|\psi\|_{W^{k-1,1}}\le\|\psi\|_{W^{k,1}}$ produces a
finite $A_{k+1}:=1+A_kB_k(1-\lambda^{k})^{-1}$ for which \eqref{eq:LY1} holds at level
$k+1$, closing the induction. All constants so obtained depend only on $K_r,\xi,k$;
in particular they depend neither on $n$ nor on which elements of $\mc U$ occur in the
composition, which is what is required for the non-stationary compositions used
throughout. \hfill$\square$

\subsection{Proof of Proposition~\ref{prop:gap}}

\emph{Base case $k=1$.} Let $\psi\in W^{1,1}(\T)$ with $\int_\T\psi=0$, and let
$F_i=F_{\alpha,\nu_i,z}\in\mc U$. By Remark~\ref{rem:BVW11}(i),
$\|\cdot\|_{\BV^1}=\|\cdot\|_{W^{1,1}}$ on $W^{1,1}(\T)$, and $L(n)\psi\in W^{1,1}$ by
Lemma~\ref{lem:LY}. Hence Lemma~\ref{lem:MLMPT} reads
\[
\|L(n)\psi\|_{W^{1,1}}=\|L(n)\psi\|_{\BV^1}\le C_1e^{-\lambda n}\|\psi\|_{\BV^1}
=C_1e^{-\lambda n}\|\psi\|_{W^{1,1}} ,
\]
which is \eqref{eq:gap} with $C_\star=C_1$, $\theta=e^{-\lambda}$. Here it is
essential that all the maps $F_i$ lie in the \emph{same} family $\mc U_z(\alpha)$: by
\eqref{eq:Uzalpha} they are then of the form $F_{\alpha,\nu_i,z}$ with a common $z$
\emph{and a common} $\alpha$, and some sequence $\{\nu_i\}\subset\M_{1,\Leb}$, which is
precisely the situation covered by Lemma~\ref{lem:MLMPT} (compare
Remark~\ref{rem:fixedfibre}). The constants $C_1,\lambda$ furnished by that lemma do
not depend on $z$, on $\alpha$ with $|\alpha|<\alpha_0$, nor on the sequence
$\{\nu_i\}$, so the estimate is uniform over the fibres and over the coupling
strengths.

\emph{Induction $k-1\Rightarrow k$, $k\ge2$.} Fix $z\in A_W$ and suppose
\eqref{eq:gap} holds at level $k-1$ with constants $C_{k-1},\theta_{k-1}$ independent
of $z$. Fix $n_1\in\N$, let $\psi\in W^{k,1}$ with $\int\psi=0$ and set
$g_i:=L(in_1)\psi$, where now all the maps are elements of the single family
$\mc U_z(\alpha)$. All $g_i$ have zero average. By Lemma~\ref{lem:LY},
\[
\|g_{i+1}\|_{W^{k,1}}\le\lambda^{kn_1}\|g_i\|_{W^{k,1}}+B_k\|g_i\|_{W^{k-1,1}},
\]
and by the induction hypothesis,
\[
\|g_{i+1}\|_{W^{k-1,1}}\le C_{k-1}\theta_{k-1}^{n_1}\|g_i\|_{W^{k-1,1}}
\le C_{k-1}\theta_{k-1}^{n_1}\|g_i\|_{W^{k,1}} .
\]
Writing $v_i=\big(\|g_i\|_{W^{k,1}},\|g_i\|_{W^{k-1,1}}\big)^{\!\top}$ we obtain
$v_{i+1}\preceq Mv_i$ componentwise, with
\[
M=\begin{pmatrix}\lambda^{kn_1}& B_k\\[2pt] C_{k-1}\theta_{k-1}^{n_1}&0\end{pmatrix}
\succeq0 .
\]
Write $c:=C_{k-1}\theta_{k-1}^{n_1}>0$ for the lower-left entry. The characteristic
polynomial of $M$ is $\mu^2-\lambda^{kn_1}\mu-B_kc$, so the spectral radius is
\begin{equation}\label{eq:specradius}
\varrho=\varrho(n_1)=\tfrac12\lambda^{kn_1}
+\tfrac12\sqrt{\lambda^{2kn_1}+4B_kC_{k-1}\theta_{k-1}^{n_1}} .
\end{equation}
Both $\lambda^{kn_1}$ and $\theta_{k-1}^{n_1}$ tend to $0$ as $n_1\to\infty$ while
$B_k$ and $C_{k-1}$ do not depend on $n_1$, so $\varrho(n_1)\to0$; fix once and for
all an $n_1$ with $\varrho(n_1)\le\varrho_k<1$. All entries of $M$ are uniform over
$\mc V$, hence over every $\mc U_z(\alpha)$, so this choice is uniform in $z$ and in
$\alpha$ as well.

Both off-diagonal entries of $M$ are strictly positive, so $M$ is irreducible, and
Perron--Frobenius provides a \emph{strictly positive} left eigenvector
$(a,b)$, $a,b>0$, with $(a,b)M=\varrho\,(a,b)$; normalise $a+b=1$. Define the balanced
norm $\|g\|_{(a,b)}:=a\|g\|_{W^{k,1}}+b\|g\|_{W^{k-1,1}}$, so that
$\|g\|_{(a,b)}=(a,b)\cdot\big(\|g\|_{W^{k,1}},\|g\|_{W^{k-1,1}}\big)^{\!\top}$.
Since $M$ has nonnegative entries and $v_{i+1}\preceq Mv_i$ componentwise with
$v_i\succeq0$,
\[
\|g_{i+1}\|_{(a,b)}=(a,b)\cdot v_{i+1}\le (a,b)\cdot Mv_i
=\varrho\,(a,b)\cdot v_i=\varrho\,\|g_{i}\|_{(a,b)} ,
\]
whence $\|g_i\|_{(a,b)}\le\varrho^{i}\|\psi\|_{(a,b)}$ and, since $a>0$ and
$\|\cdot\|_{(a,b)}\le\|\cdot\|_{W^{k,1}}$,
\begin{equation}\label{eq:alongmultiples}
\big\|L(in_1)\psi\big\|_{W^{k,1}}\le a^{-1}\varrho^{i}\,\|\psi\|_{W^{k,1}} .
\end{equation}
For a general $n$ write $n=in_1+s$ with $0\le s<n_1$. Applying \eqref{eq:LY2} to the
last $s$ factors and then \eqref{eq:alongmultiples},
\[
\|L(n)\psi\|_{W^{k,1}}\le\lambda^{ks}\|L(in_1)\psi\|_{W^{k,1}}
+B_k\|L(in_1)\psi\|_{W^{k-1,1}}
\le\frac{1+B_k}{\min(a,b)}\,\varrho^{i}\,\|\psi\|_{W^{k,1}} ,
\]
and $\varrho^{i}\le\varrho^{(n-n_1)/n_1}=\varrho^{-1}\big(\varrho^{1/n_1}\big)^{n}$.
This is \eqref{eq:gap} at level $k$ with
$C_\star=(1+B_k)\varrho^{-1}\min(a,b)^{-1}$ and
$\theta=\varrho^{1/n_1}<1$. This is the standard two-norm argument; compare
\cite[Prop.~9]{Galatolo2022} and \cite{Liverani1995}.

\emph{Fibred statement.} For $v\in\W^k_0$ we have
$\big(\LL_{\alpha,\phi^{(n)}}\cdots\LL_{\alpha,\phi^{(1)}}v\big)_z
=L_{\alpha,\phi^{(n)},z}\cdots L_{\alpha,\phi^{(1)},z}v_z$ with $\int_\T v_z=0$ and all
maps in $\mc U_z(\alpha)$ for a.e.\ $z$ --- the same $z$ and the same $\alpha$ in
every factor, as required ---;
applying \eqref{eq:gap} fibrewise and taking $\esssup_z$ gives the claim, and \eqref{eq:resolvent0} follows by summing the Neumann
series. \hfill$\square$

\subsection{Existence of the second-order term (Proposition~\ref{prop:formalgenuine})}

Let $1\le k\le r-5$. By Theorem~\ref{thm:LR} applied at level $k+2$, which
is legitimate since $k+2\le r-3$, the map $\alpha\mapsto\phi^\star(\alpha)$ is $\mc C^1$
into $\W^{k+2}$ and
\[
D(\alpha):=\partial_\alpha\phi^\star(\alpha)
=\mc R^{-1}_{\alpha,\phi^\star(\alpha)}\ \mc K_{\alpha,\phi^\star(\alpha)}
\phi^\star(\alpha) .
\]
We show that $\alpha\mapsto D(\alpha)$ is differentiable at $0$ in $\W^{k}$, which
gives \eqref{eq:expansion} with the second-order coefficient $\tfrac12 D'(0)$. The
right-hand side is a composition of three maps, each differentiable at $\alpha=0$.

\emph{(i)} $\alpha\mapsto\phi^\star(\alpha)\in\W^{k+2}$ is differentiable at $0$ with
derivative $\barW\,\Up$, by Theorems~\ref{thm:LR} and \ref{thm:zerocoupling} at level
$k+2$.

\emph{(ii)} $\alpha\mapsto N(\alpha):=\mc K_{\alpha,\phi^\star(\alpha)}
\phi^\star(\alpha)\in\W^{k}$ is differentiable at $0$. Write
$N(\alpha)_z=-\partial_x G_z(\alpha)$ with
$G_z(\alpha):=L_{\alpha,\phi^\star(\alpha),z}\big(A_{\phi^\star(\alpha),z}\,
\phi^\star_z(\alpha)\big)$, and split
\[
G_z(\alpha)-G_z(0)
=\underbrace{\Big[L_{\alpha,\phi^\star(\alpha),z}-P\Big]
\big(A_{\phi^\star(\alpha),z}\phi^\star_z(\alpha)\big)}_{(\mathrm I_\alpha)}
+\underbrace{P\Big(A_{\phi^\star(\alpha),z}\phi^\star_z(\alpha)
-A_{\rho,z}\rho\Big)}_{(\mathrm{II}_\alpha)} .
\]
For $(\mathrm I_\alpha)$ apply \eqref{eq:exact1} with $F=f$ and
$\Delta=\alpha A_{\phi^\star(\alpha),z}$, which is legitimate since
$k+3\le r-2$; using $\|\Delta-\alpha A_{\rho,z}\|_{\mc C^{k+2}}
=O(\alpha\|\phi^\star(\alpha)-\rho\|_0)=O(\alpha^2)$ by
Lemma~\ref{lem:Abound}(ii) and Theorem~\ref{thm:stability}, and
$\|A_{\phi^\star(\alpha),z}\phi^\star_z(\alpha)-A_{\rho,z}\rho\|_{W^{k+2,1}}
=O(\alpha)$, one gets
\[
(\mathrm I_\alpha)=-\alpha\Big[P\big(A_{\rho,z}^2\,\rho\big)\Big]'
+o(\alpha)\qquad\text{in } W^{k+1,1}(\T),
\]
uniformly in $z$. For $(\mathrm{II}_\alpha)$, bilinearity gives
\[
A_{\phi^\star(\alpha),z}\phi^\star_z(\alpha)-A_{\rho,z}\rho
=A_{\phi^\star(\alpha)-\rho,z}\,\phi^\star_z(\alpha)
+A_{\rho,z}\big(\phi^\star_z(\alpha)-\rho\big)
=\alpha\big(A_{\dot\phi,z}\,\rho+A_{\rho,z}\,\dot\phi_z\big)+o(\alpha),
\]
by (i), Lemma~\ref{lem:Abound}(ii) and Lemma~\ref{lem:calculus}, with
$\dot\phi_z=\barW(z)\Up$. Adding the two and differentiating in $x$,
\begin{equation}\label{eq:Ndot}
N'(0)_z=\big[P\big(A_{\rho,z}^2\,\rho\big)\big]''
-\big[P\big(A_{\dot\phi,z}\,\rho\big)\big]'
-\big[P\big(A_{\rho,z}\,\dot\phi_z\big)\big]' ,
\end{equation}
all three terms lying in $\W^{k}$ by Lemmas~\ref{lem:Abound} and \ref{lem:LY}.

\emph{(iii)} The family $R_\alpha:=\mc R_{\alpha,\phi^\star(\alpha)}$ satisfies the
hypotheses of Lemma~\ref{lem:Rincr}: $\phi^\star(\alpha)\in\Delta^{k+2}$ with
$\|\phi^\star(\alpha)\|_{k+2}\le M_{k+2}$ by Proposition~\ref{prop:apriori}, the
Lipschitz bound $\|\phi^\star(\alpha)-\phi^\star(\beta)\|_{k+1}\le C|\alpha-\beta|$
holds by Theorem~\ref{thm:stability} at level $k+1$, and $\phi^\star(0)=\rho$.
Consequently, for every \emph{fixed} $v\in\W^{k+2}_0$,
\begin{equation}\label{eq:Rdotapp}
\frac{R_\alpha-R_0}{\alpha}\,v\ \longrightarrow\ \dot R\,v
=\big[P\big(A_{\rho,z}\,v_z\big)\big]'-\mc K_{0,\rho}v
\qquad\text{in }\W^{k}_0 ,
\end{equation}
with the quantitative rate \eqref{eq:Rdot}. We stress that we do \emph{not} assert
differentiability of $\alpha\mapsto R_\alpha$ in $\mc L(\W^{k+1},\W^{k})$: that
assertion is false even for the doubling map with $h\equiv1$ and $W\equiv1$, as
Remark~\ref{rem:nodiff} shows. Everything below uses \eqref{eq:Rdotapp} for one
specific vector, whose membership in $\W^{k+2}_0$ is verified presently. Finally
$R_0=\Id-P$ is invertible on $\W^{k}_0$ by Proposition~\ref{prop:gap}, and
$R_\alpha^{-1}$ is bounded on $\W^{k}_0$ uniformly in $\alpha$ by
\eqref{eq:resolventbound}.

To combine (i)--(iii) we must not differentiate $\alpha\mapsto\mc R_\alpha^{-1}$ in the
operator norm of $\mc L(\W^k_0)$, which is unavailable
(Remark~\ref{rem:nodensity}); we argue instead on the single vector produced by the
increment quotient. Write $R_\alpha:=\mc R_{\alpha,\phi^\star(\alpha)}$,
$R_0=\Id-P$, $\dot R$ for the operator of Lemma~\ref{lem:Rincr}(ii) and
$u_0:=R_0^{-1}N(0)$, so that $D(\alpha)=R_\alpha^{-1}N(\alpha)$. The
resolvent identity $R_\alpha^{-1}-R_0^{-1}=R_\alpha^{-1}(R_0-R_\alpha)R_0^{-1}$ gives
\begin{equation}\label{eq:Dincrement}
D(\alpha)-D(0)=R_\alpha^{-1}\big(N(\alpha)-N(0)\big)
+R_\alpha^{-1}\big(R_0-R_\alpha\big)u_0 ,
\end{equation}
so that, dividing \eqref{eq:Dincrement} by $\alpha$ and introducing the increment
quotients,
\begin{equation}\label{eq:quotients}
a_\alpha:=\frac{N(\alpha)-N(0)}{\alpha},\qquad
b_\alpha:=\frac{R_0-R_\alpha}{\alpha}\,u_0,\qquad
\frac{D(\alpha)-D(0)}{\alpha}=R_\alpha^{-1}\big(a_\alpha+b_\alpha\big) .
\end{equation}
Here $u_0\in\W^{k+2}_0$: indeed $N(0)=\mc K_{0,\rho}\rho\in\W^{k+2}_0$ by
Proposition~\ref{prop:Kbounded} at level $k+2$, legitimate since $k+2\le r-2$, and
$R_0^{-1}=(\Id-P)^{-1}$ preserves $\W^{k+2}_0$ by Proposition~\ref{prop:gap} at level
$k+2$. These two spare derivatives on $u_0$ are exactly what
Lemma~\ref{lem:Rincr}(ii) demands, and they are the reason the second quotient
converges at all: by Remark~\ref{rem:nodiff} no amount of uniform control on
$\|u_0\|_{k+1}$ alone would suffice. By (ii) we have $a_\alpha\to N'(0)$ in
$\W^{k}_0$, and by (iii), that is by Lemma~\ref{lem:Rincr}(ii) applied to the fixed
vector $u_0\in\W^{k+2}_0$, $b_\alpha\to-\dot R\,u_0$ in $\W^{k}_0$.

\emph{(iv) The limit vector has one derivative to spare.} Set
\[
w:=N'(0)-\dot R\,u_0
=N'(0)-\Big[P\big(A_{\rho,z}\,(u_0)_z\big)\Big]'+\mc K_{0,\rho}u_0 ,
\]
using the formula for $\dot R$ in Lemma~\ref{lem:Rincr}(ii). We claim that
$w\in\W^{k+1}_0$, not merely in $\W^{k}_0$. For the third term this is
Proposition~\ref{prop:Kbounded} at level $k+1$, which costs nothing precisely because
$\mc K$ is smoothing: $\|\mc K_{0,\rho}u_0\|_{k+1}\le C_{\mc K}\|\rho\|_{k+2}\|u_0\|_0$,
admissible since $k+1\le r-2$. For the second term, $u_0\in\W^{k+2}_0$ by the
paragraph above, $\|A_{\rho,z}\|_{\mc C^{k+2}}$ is bounded uniformly in $z$ by
Lemma~\ref{lem:Abound}(i), so $A_{\rho,z}(u_0)_z\in W^{k+2,1}(\T)$ by
Lemma~\ref{lem:calculus}, $P$ preserves $W^{k+2,1}(\T)$ by Lemma~\ref{lem:LY} at level
$k+2\le r-1$, and one derivative leaves $W^{k+1,1}(\T)$. For the first term, inspect
the three summands of \eqref{eq:Ndot}: $A_{\rho,z}^2\rho\in W^{k+3,1}(\T)$ since
$\rho\in\mc C^{r-1}$ and $k+3\le r-1$, while $A_{\dot\phi,z}\rho$ and
$A_{\rho,z}\dot\phi_z$ lie in $W^{k+2,1}(\T)$, the latter because
$\dot\phi_z=\barW(z)\Up$ with $\Up\in W^{r-2,1}(\T)$ and $k+2\le r-2$; applying $P$ and
then two, respectively one, derivatives leaves $W^{k+1,1}(\T)$ in all three cases.
Finally each summand of $w$ has zero fibrewise average, being either a total derivative
or in the range of $\mc K$. Every index constraint invoked here is implied by
$k\le r-5$, so no hypothesis beyond those of the statement is needed.

\emph{Conclusion.} Decompose
\[
R_\alpha^{-1}\big(a_\alpha+b_\alpha\big)-R_0^{-1}w
=R_\alpha^{-1}\big(a_\alpha+b_\alpha-w\big)
+\big(R_\alpha^{-1}-R_0^{-1}\big)w .
\]
The first term tends to $0$ in $\W^k_0$ because $\|a_\alpha+b_\alpha-w\|_{k}\to0$ by
(ii)--(iii) while $\|R_\alpha^{-1}\|_{\W^k_0\to\W^k_0}$ is bounded uniformly in
$\alpha$ by \eqref{eq:resolventbound}. The second tends to $0$ by
Lemma~\ref{lem:resvector}, which applies exactly because of (iv), with the quantitative
bound $C|\alpha|\,\|w\|_{k+1}$. Hence the quotient \eqref{eq:quotients} converges and
\begin{equation}\label{eq:Dprime}
D'(0)=R_0^{-1}w=R_0^{-1}N'(0)-R_0^{-1}\,\dot R\,R_0^{-1}N(0) ,
\end{equation}
which is the expected formula. The minus sign is the one carried by the identity
$\partial_\alpha(R_\alpha^{-1})
=-R_\alpha^{-1}(\partial_\alpha R_\alpha)R_\alpha^{-1}$
in the formal computation --- a computation we have replaced, rather than justified,
since the operator derivative appearing in it does not exist in the norm that would be
required. This establishes the existence of the second-order
term; by Proposition~\ref{prop:recursion} its
coefficient is the formal susceptibility $\phi_2$ of \eqref{eq:formaldef}, so
\eqref{eq:expansion} holds with $\phi_2=\tfrac12D'(0)$.

Here \eqref{eq:Dprime} serves only to establish \emph{existence} of the second-order
term; the closed formula \eqref{eq:secondorder} comes from the recursion
\eqref{eq:formaldef} at $j=2$, an elementary identification of the coefficient of
$\alpha^2$ in the fixed point equation which does not involve the resolvent of
$\mc R_\alpha$. \hfill$\square$

\end{appendices}


\end{document}